\documentclass[10pt]{article}

\usepackage[english,activeacute]{babel}
\usepackage{epsfig}

\usepackage[latin1]{inputenc}
\usepackage[T1]{fontenc}
\usepackage{stmaryrd}
\usepackage{amsmath,amsfonts,amssymb,mathrsfs,amsthm}
\usepackage{xcolor}
\usepackage{dsfont}
\usepackage{graphicx}
\usepackage[mathscr]{eucal}
\usepackage{makeidx}
\usepackage{verbatim}
\usepackage{graphics,graphicx}
\usepackage{textcomp}
\usepackage{float}
\usepackage[colorlinks=true]{hyperref}

\newtheorem{lemma}{Lemma}[section]
\newtheorem{theorem}[lemma]{Theorem}
\newtheorem{proposition}[lemma]{Proposition}
\newtheorem{corollary}[lemma]{Corollary}
\newtheorem{assumption}[lemma]{Assumption}

\theoremstyle{definition}
\newtheorem{example}[lemma]{Example}
\newtheorem{definition}[lemma]{Definition}
\newtheorem{remark}[lemma]{Remark}

\makeatletter
\def\keywords{
    \vspace{1ex}
    \noindent
    \if@twocolumn
      \small{\bf  Keywords}\/---$\!$    \else
      \begin{center}\small\ {\bf Keywords}\end{center}\quotation\small
    \fi}
\def\endkeywords{\vspace{0.6em}\par\if@twocolumn\else\endquotation\fi
    \normalsize\rm}
\makeatother

\newcommand{\calR}{\ensuremath{\mathcal R}}
\renewcommand{\P}{\ensuremath{\mathcal P}}

\newcommand{\calN}{\ensuremath{\mathcal N}}
\newcommand{\calS}{\ensuremath{\mathcal S}}

\newcommand{\calC}{\ensuremath{\mathcal C}}

\newcommand{\ID}{\ensuremath{\mathbb D}}

\DeclareMathOperator{\Meas}{Meas}

\DeclareMathOperator{\Sp}{Sp}
\DeclareMathOperator{\loc}{loc}

\DeclareMathOperator{\Vol}{Vol}

\DeclareMathOperator{\Int}{Int}

\DeclareMathOperator{\Eucl}{Eucl}
\DeclareMathOperator{\diam}{diam}

\newcommand{\mb}[1]{\ensuremath{\mathbb{#1}}}
\newcommand{\N}{{\mb{N}}}

\newcommand{\R}{{\mb{R}}}
\newcommand{\C}{{\mb{C}}}

\newcommand{\f}{\ensuremath{\mathfrak{f}}}

\newcommand{\de}{\delta}
\newcommand{\eps}{\varepsilon}

\newcommand{\M}{\ensuremath{\mathcal M}}
\newcommand{\T}{\ensuremath{\mathbb T}}

\let \Re \relax
\DeclareMathOperator{\Re}{Re}

\newcommand{\ovl}[1]{\overline{#1}}

\renewcommand{\S}{\ensuremath{\mathbb S}}

\DeclareMathOperator{\supp}{supp}

\DeclareMathOperator{\dist}{dist}

\DeclareMathOperator{\vect}{span}

\renewcommand{\d}{\ensuremath{\partial}}

\let \div \relax
\DeclareMathOperator{\div}{div}

\DeclareMathOperator{\length}{length}

\newcommand{\E}{\mathscr E}
\newcommand{\Z}{\mathbb Z}

\newcommand{\X}{\ensuremath{\mathcal X}}

\newcommand\bna{\begin{eqnarray*}}
\newcommand\ena{\end{eqnarray*}}

\def\bal#1\nal{\begin{align*}#1\end{align*}}
\def\baln#1\naln{\begin{align}#1\end{align}}

\newcommand\bnan{\begin{eqnarray}}
\newcommand\enan{\end{eqnarray}}

\newcommand\bnp{\begin{proof}}
\newcommand\enp{\end{proof}}

\newcommand\bneq{\begin{eqnarray*}\left\lbrace \begin{array}{rcl}}
\newcommand\eneq{\end{array} \right.\end{eqnarray*}}
\newcommand\bneqn{\begin{eqnarray}\left\lbrace \begin{array}{rcl}}
\newcommand\eneqn{\end{array} \right.\end{eqnarray}}

 \numberwithin{equation}{section} 
 
\newcommand\grando[1]{\mathcal{O}\left(#1\right)} 
\newcommand\Lap{\Delta_{g}}
\newcommand\nablag{\nabla_{g}}

\newcommand\nor[2]{\left\|#1\right\|_{#2}}
\newcommand\sgn{\textnormal{sgn}}
\newcommand\e{\varepsilon}

\newcommand{\Ek}{\mathcal{E}_k}

\newcommand{\ce}{c}
\newcommand{\q}{q}
\newcommand{\qf}{q_{\f}}
\newcommand\CO[1]{\mathbf{C}_{#1}}

\begin{document}
\title{On uniform observability of gradient flows in the vanishing viscosity limit}

\author{Camille Laurent\footnote{CNRS UMR 7598 and Sorbonne Universit\'es UPMC Univ Paris 06, Laboratoire Jacques-Louis Lions, F-75005, Paris, France, email: laurent@ann.jussieu.fr} and Matthieu L\'eautaud\footnote{Laboratoire de Math\'ematiques d'Orsay, Universit\'e Paris-Sud, CNRS, Universit\'e Paris-Saclay, B\^atiment 307, 91405 Orsay Cedex France, email: matthieu.leautaud@math.u-psud.fr.}
}

\date{}

\maketitle

\begin{abstract}
We consider a transport equation by a gradient vector field with a small viscous perturbation $-\eps\Delta_g$. We study uniform observability (resp. controllability) properties in the (singular) vanishing viscosity limit $\eps\to 0^+$, that is, the possibility of having a uniformly bounded observation constant (resp. control cost). We prove with a series of examples that in general, the minimal time for uniform observability may be much larger than the minimal time needed for the observability of the limit equation $\eps=0$. 
We also prove that the two minimal times coincide for positive solutions.
The proofs rely on a semiclassical reformulation of the problem together with (a) Agmon estimates concerning the decay of eigenfunctions in the classically forbidden region~\cite{HS:84} (b) fine estimates of the kernel of the semiclassical heat equation~\cite{LY:86}.
\end{abstract}


\setcounter{tocdepth}{2} 
\tableofcontents

\section{Introduction and main results}

\subsection{Introduction}
\label{s:intro-intro}
Given a smooth connected compact manifold $\M$ without boundary (the case of a bounded domain of $\R^n$ is also discussed in Section~\ref{s:main-res} below), a smooth real valued vector field $X$ on $\M$ and a real valued potential $\q(x)$, we consider the question of observability/detectability for the autonomous transport equation
  \begin{equation}
  \label{e:transport-libre}
 \left\{
 \begin{aligned}
&(\d_t - X-\q) u = 0 , & \text{ in } \R\times \M ,  \\
& u|_{t=0} = u_0, & \text{ on } \M ,
 \end{aligned}
 \right.
 \end{equation}
from an observation (open) set $\omega \subset \M$ through the time interval $(0,T)$. More precisely, the question is whether there exists a constant $C_0=C_0(T,\omega)>0$ such that 
 \begin{multline}
 \label{e:transport-obs-ineq}
C_0^2 \int_0^T\int_\omega |u(t,x)|^2 ds(x)dt \geq \|u (T)\|_{L^2(\M)}^2 , \\
 \text{ for all } u_0 \in L^2(\M) \text{ and } u \text{ solution of~\eqref{e:transport-libre}}.
 \end{multline}
 Here, $ds(x)$ denotes any positive density measure\footnote{See e.g.~\cite[Chapter~16 p427]{Lee:book}: given a local chart $(U_\kappa, \kappa)$ of $\M$, we have $\int_{U_\kappa} u \ ds=  \int_{\kappa(U_\kappa)} u\circ \kappa^{-1}(y) \varphi^\kappa(y) dy$ for an appropriate smooth positive function $\varphi^\kappa$, and for any $u \in C^0_c(U_\kappa)$.} on $\M$, and the $L^2$ norm is defined accordingly.
The observability question~\eqref{e:transport-obs-ineq} is naturally solved by introducing an appropriate Geometric Control Condition (recall $\d \M=\emptyset$):  we say that $(\M,X,\omega,T)$ satisfies (GCC) if for all $x\in \M$, there is $t \in (0,T)$ such that $\phi_{-t}(x)\in \omega$, where $(\phi_t)_{t\in\R}$ denotes the flow of $X$ (see Section~\ref{s:GCC} for precise statements and proofs). We also say that $(\M,X,\omega)$ satisfies (GCC) if $(\M,X,\omega,T)$ does for some $T>0$; and if so, we denote by $T_{GCC}(\M,X,\omega)$ the infimum of times for which $(\M,X,\omega,T)$ satisfies (GCC).

\bigskip
On the other hand, endowing $\M$ with a Riemannian metric $g$, one may want to investigate the observability question for the viscously damped transport equation:
  \begin{equation}
  \label{e:transport-viscous}
 \left\{
 \begin{aligned}
&(\d_t - X -\q - \eps \Delta_g ) u = 0 , & \text{ in } \R^+_*\times \M ,  \\
& u|_{t=0} = u_0, & \text{ on } \M ,
 \end{aligned}
 \right.
 \end{equation}
from the same observation set $(0,T)\times\omega$. The question is whether there exists a constant $C_0(T,\eps)>0$ such that 
 \begin{multline}
 \label{e:transport-viscous-obs}
C_0(T,\eps)^2 \int_0^T\int_\omega |u(t,x)|^2 ds(x)dt \geq \|u (T)\|_{L^2(\M)}^2 , \\
 \text{ for all } u_0 \in L^2(\M) \text{ and } u \text{ solution of~\eqref{e:transport-viscous}},
 \end{multline}
(and one may then choose the Riemannian volume density $ds(x) = d\Vol_g(x)$ without changing the problem).
For fixed $\eps>0$, Equation~\eqref{e:transport-viscous} is of parabolic type and the observability inequality~\eqref{e:transport-viscous-obs} is known to hold for any open set $\omega \neq \emptyset$ and $T>0$, see~\cite{FI:96} (see also~\cite{LR:95} and its variant in~\cite{Lea:10}).
Of course, in such results, the observability constant $C_0(T,\eps)$ in~\eqref{e:transport-viscous-obs} depends {\em a priori} on $\eps$. 
For many different reasons (some of them described in Section~\ref{s:motiv} below), it is interesting to investigate the behavior of the observability constant $C_0(T,\eps)$ in the vanishing viscosity limit $\eps \to 0^+$. This problem was first studied in the one dimensional setting by Coron and Guerrero in~\cite{CG:05}, and later extended to any dimension by Guerrero and Lebeau~\cite{GL:07}. Their main result in this direction can be formulated (in the present geometric context, see the remark preceding Proposition~\ref{p:guerrero-lebeau}) as follows.
\begin{theorem}[Guerrero-Lebeau~\cite{GL:07}] 
\label{t:Leb-Guer}
Given an open set $\omega \subset \M$, the following two results hold.
\begin{itemize}
\item \cite[Theorem~1]{GL:07} Assume $(\M,X,\overline{\omega}, T)$ does not satisfy (GCC). Then there is $C, \eps_0>0$ such that any constant $C_0(T,\eps)$ in~\eqref{e:transport-viscous-obs} satisfies $C_0(T,\eps)\geq \exp(C/\eps)$ for $\eps \in (0,\eps_0)$.
\item \cite[Theorem~3]{GL:07} Assume $(\M,X,\omega)$ satisfies (GCC). Then there is $T_{unif}(\omega) \geq T_{GCC}(\M,X,\omega)$ and $K_0>0$ such that for all $T\geq T_{unif}(\omega)$,~\eqref{e:transport-viscous-obs} holds with $C_0(T,\eps)\leq K_0$ for all $\eps \leq 1$.
\end{itemize}
\end{theorem}
Note that the results in~\cite{GL:07} are even more general since time-dependent vector fields are allowed and the boundary-value problem is also considered (with Dirichlet boundary conditions). We also refer to our Proposition \ref{proplowerobse+} below for a more precise lower bound of the constant $C$ when (GCC) is not satisfied.

Note that if~\eqref{e:transport-viscous-obs} holds for some $T_0$ and constant $C_0(T_0,\eps)$, then it also holds for all times $T\geq T_0$ with the same constant $C_0(T_0,\eps)$.
In~\cite{GL:07}, the question of the minimal time $T_{unif}(\omega)$, more precisely defined by
\bal
T_{unif}(\omega) & = \inf \left\{ T >0 \text{ for which there exist }K_0, \eps_0>0 \right. \\
&\left.  \quad  \text {such that~\eqref{e:transport-viscous-obs} holds with }C_0(T,\eps)\leq K_0 \text{ for all }  \eps \in (0,\eps_0) \right\} ,
\nal
and its link with the minimal observation time $T_{GCC}(\M,X,\omega)$ associated to the limit problem~\eqref{e:transport-libre} is left open. 
In particular, the formulation of the results in~\cite{GL:07} (see e.g. Theorem~2 and the discussion thereafter in that reference) suggests the possible existence of a universal constant $\mathfrak{K}\geq 1$ such that 
\begin{equation}
\label{e:universal-K}
T_{unif}(\omega) \leq \mathfrak{K}\ T_{GCC}(\M,X,\omega) .
\end{equation}

\bigskip
The present article investigates this question in a very particular case, namely assuming the vector field $X$ is a gradient vector field, i.e. $X= \nabla_g \f$ for a function $\f\in W^{2,\infty}(\M;\R)$ (note that the gradient is taken with respect to the Riemannian metric $g$). Hence, Equation~\eqref{e:transport-viscous} becomes
  \begin{equation}
  \label{e:transport-viscous-gradient}
 \left\{
 \begin{aligned}
&(\d_t - \nabla_g \f \cdot \nabla_g -\q - \eps \Delta_g ) u = 0 , & \text{ in } \R^+_*\times \M ,  \\
& u|_{t=0} = u_0, & \text{ on } \M ,
 \end{aligned}
 \right.
 \end{equation}

Here, given two vector fields $Y_1$ and $Y_2$, we have denoted $Y_1\cdot Y_2 = g(Y_1,Y_2)$ or  $(Y_1\cdot Y_2)(x) = Y_1(x)\cdot Y_2(x) = g_x(Y_1(x),Y_2(x))$ for all $x\in \M$. We denote similarly $|Y|_g = \sqrt{Y \cdot Y}$ the associated Riemannian norm of a vector (or a vector field).
Note that the vector field $\nabla_g \f$ is canonically identified with the derivation $\nabla_g \f \cdot \nabla_g$.

In this context, the first consequence of our main results can be (loosely) stated as follows. 
\begin{theorem}
\label{t:loose}
\begin{enumerate}
 \item There are geometries $(\M,g)$ such that for all $\Lambda>0$, one can find $\f \in C^\infty(\M)$ and $\omega$ open such that $(\M, \nablag \f , \omega)$ satisfies (GCC) and  $T_{unif}(\omega) \geq \Lambda \ T_{GCC}(\M, \nablag \f , \omega)$.
 \item There are $(\M,\f, X,\omega)$ such that for all $\Lambda>0$, one can find a metric $g_\Lambda$ on $\M$ such that 
 \begin{itemize}
 \item $X = \nabla_{g_\Lambda}\f$,
 \item  $(\M, X , \omega)$ satisfies (GCC), 
 \item $T_{unif}(\omega) \geq \Lambda  \ T_{GCC}(\M, X , \omega)$.
 \end{itemize}
\end{enumerate}
\end{theorem}
In particular, Theorem~\ref{t:loose} states that there is no $\mathfrak{K}$ such that~\eqref{e:universal-K} holds  for all $(\M,X,\omega)$.

The second item in Theorem~\ref{t:loose} stresses the importance of the viscosity one chooses. Namely, with the same vector field $X$, changing the metric $g$, that is the viscous perturbation, may change the minimal uniform observability time.
 We also obtain related results for domains of $\R^{n}$ (see Section~\ref{s:sect-general}).

\bigskip
Our second main result in this setting concerns the uniform observability of {\em positive} solutions to~\eqref{e:transport-viscous-gradient}. Recall that nonnegative data $u_0\geq 0$ give rise to positive solutions to~\eqref{e:transport-viscous-gradient}. We define $C_0^+(T,\eps)$ the observability constant for positive solutions, that is for which~\eqref{e:transport-obs-ineq} holds for all $u_0\geq0$, and accordingly set 
\begin{multline}
\label{defTunif+}
T_{unif}^+(\omega) = \inf \left\{ T >0 \text{ for which there exist }K_0, \eps_0>0 \text{ such that~\eqref{e:transport-viscous-obs} holds}\right. \\
 \left.  \quad\quad \quad \text{for all }u_0\geq 0, \text{ with }C_0^+(T,\eps)\leq K_0 \text{ for all }  \eps \in (0,\eps_0) \right\} .
\end{multline}

\begin{theorem}[Positive solutions]
\label{t:positive-intro}
For all $\f \in C^3(\M;\R)$, and $\omega \subset \M$ such that $(\M,\nabla_g\f,\omega)$ satisfies (GCC), we have $T_{unif}^+(\omega) = T_{GCC}(\M,\nabla_g\f,\omega)$.
\end{theorem}

As usual, these uniform observability/non-observability results can be reformulated in terms of uniform controllability/non-controllability statements for an adjoint controlled equation, see Section~\ref{s:obs-cont} below.

\subsection{Background and motivation}
\label{s:motiv}
Uniform controllability problems for singular perturbations of partial differential equations already appeared in the reference book of Lions~\cite[Chapter 3]{Lio2:88}. In the context of transport/heat equation in vanishing viscosity limit, this study was initiated by Coron and Guerrero on the 1D problem with constant speed in~\cite{CG:05}, where the authors make a conjecture on the minimal time needed to achieve uniform controllability. Then, the estimates on this minimal time have been improved successively in~\cite{Gla:10,Lissy:12,Lissy:14,Lissy:15}. 
We also refer to the articles~\cite{Munch:18,YM:19,YM:19bis} proposing numerical experiments to find the optimal minimal time. 
 Such uniform control properties in singular limits are also addressed for vanishing dispersion in~\cite{GG:08} and for vanishing dispersion and viscosity in~\cite{GG:09}. 
 
Whereas the one dimensional problem with a constant vector field has received a lot of attention in the past fifteen years, there are very few results in higher dimension or for non-constant vector fields. Besides~\cite{GL:07} we are only aware of the results of~\cite{Barcena:20} and~\cite{Barcena:21} for the flat Laplace operator and the vector field $\d_{x_1}$, with several boundary conditions.

Note that controllability problems for nonlinear conservation laws in vanishing viscosity have also been studied in~\cite{GG:07},~\cite{Lea:11}, and~\cite{Marbach:14}.

\bigskip
There are several motivations for studying the vanishing viscosity limit. 
A first motivation comes from the theory of conservation laws, for which the vanishing viscosity criterium is a selection principle for the physical (called entropy) solution, see~\cite{Kru:70} or~\cite[Chapter 6]{Daf:00}.
It is therefore very natural, when considering control problems for conservation laws, to study the cost of the viscosity,
that is, to determine if known controllability properties for the hyperbolic equation are still valid for the model with
small viscosity, and how the size of the control evolves as the viscosity approaches $0$. So far the only known results in this directions seem to be~\cite{GG:07} and~\cite{Lea:11}.

\bigskip
Another important motivation for studying singular limits in control problems is the seek of controllability properties for the perturbated system itself. This is well-illustrated by the papers~\cite{Cor:96,CF:96,ChaNavier:09,CMS:19}, where the authors investigate the Navier-Stokes system with Navier slip or slip-with-friction boundary conditions. They use a global controllability result for the inviscid equation (in this case, the Euler equation) to deduce global approximate controllability of the Navier-Stokes system. 

\bigskip
On the other hand, the study of gradient fields naturally arises as the simplest dynamical situation among all vector fields.
The importance of gradient vector fields with a vanishing viscosity coefficient also appears in theoretical physics and differential topology, through the Witten-Helffer-Sj\"ostrand theory~\cite{Witten:82,HS:85}. See e.g. the monographs~\cite{Helffer:booksemiclassic,CFKS}. 
In that theory, the operator $- \nabla_g \f \cdot \nabla_g - \eps \Delta_g$ (and its analogues acting on forms) is conjugated to a particular semiclassical Schr\"odinger operator, namely 
\baln
\label{e:intro-witten}
P_\eps = -\eps^2 \Delta_g + \frac{|\nabla_g \f|_g^2}{4} + \frac{\eps}{2} \Delta_g \f,
\naln 
sometimes called the Witten Laplacian.
Topological properties of the couple $(\M,\f)$ (e.g. the Morse inequalities, linking the number of critical points of the Morse function $\f$ with the Betti numbers of $\M$) are deduced from spectral properties of the Witten Laplacian.
We also refer to~\cite{DR:20} for the understanding of other links between the spectral theories of the Witten Laplacian and the vector field $\nabla_g \f \cdot \nabla_g$ (in appropriate spaces linked to the dynamics of the gradient flow), in the semiclassical limit $\eps\to 0^+$.

\bigskip
Viscous perturbations of gradient dynamics also arise naturally in molecular dynamics. 
Indeed, in $\R^n$, the operator $- \nabla \f \cdot \nabla - \eps \Delta$ is the infinitesimal generator of the so-called overdamped Langevin process
$$
dX_t = \nabla \f(X_t) dt + \sqrt{2\eps} dB_t ,
$$
where $X_t\in \R^n$ and $(B_t)_{t\geq0}$ is a standard Brownian motion of dimension $n$. This stochastic process is a classical model in statistical physics. It is used in particular for the simulation of molecular dynamics at low temperature (proportional to $\eps$), see~\cite{Chan:43,SM:79}. The possible convergence to equilibrium, as well as the so-called metastability phenomenon are closely related to the low-lying eigenvalues (and associated eigenfunctions) of $- \nabla \f \cdot \nabla - \eps \Delta$, or equivalently of the Witten Laplacian~\eqref{e:intro-witten}.
For a very precise asymptotic study of the exponentially small eigenvalues of this operator, we refer e.g. to \cite{HKN:04,Michel:19} in the case of a compact manifold and \cite{HN:06,LePeutrec:10,DLLN:19} in the case with boundary (see also the references therein).

We finally remark that the above-mentioned works concerning the Witten Laplacian mostly study the behavior of the bottom of the spectrum $P_\eps$ (thus linked to critical points of $\f$). 
In the present work, we rely on a similar conjugation.

\subsection{Main results}
\label{s:main-res}
As already seen in the end of Section~\ref{s:intro-intro}, the results of this article go in two different directions. 

In a first part (Section~\ref{s:general-bounds}), we prove some general lower bounds on the time $T_{unif}(\omega)$ for a general class of domains and vector fields. This implies in particular that the quite natural idea that $T_{unif}(\omega)$ is $T_{GCC}(\M,X,\omega)$ or even $\mathfrak{K}T_{GCC}(\M,X,\omega)$ for a universal constant $\mathfrak{K}$ is false in general. This might be interpreted by the fact that in the vanishing viscosity limit, some strong oscillations can be responsible for concentration phenomena. The latter are not only described by the flow of $X$, but other parameters where an Agmon distance plays a crucial role. We also study (in Section~\ref{s:revolution}) the particular case of surfaces of revolutions where we obtain refined lower bounds. This analysis also shows that the global geometry of the Riemannian manifold we consider has an effect on the vanishing viscosity limit. In particular, with a fixed vector field, we show that the choice of the Laplacian $\Lap$ can change drastically the time $T_{unif}(\omega)$ of uniform controllability. This shows definitely that the flow of the vector field is not the unique parameter defining $T_{unif}(\omega)$.

In a second part (Section~\ref{s:positive}), our results go exactly in the opposite direction, but for a specific class of solutions, namely positive solutions. As announced in Theorem \ref{t:positive-intro}, we prove that $T_{unif}(\omega)^+=T_{GCC}(\M,\nabla_g\f,\omega)$ for positive solutions.  
This shows that the dynamics of positive solutions are actually well represented by the sole flow of $\nabla_g\f$.

In both cases, using the change of unknown $v= e^{\frac{\f}{2\eps}} u$, see Section \ref{sectconjug}, the problem is reduced (modulo lower order terms, and in weighted spaces) to observability of solutions of a semiclassical heat equation 
\bal
\e\partial_t v- \eps^2 \Lap v+Vv =0 , 
\nal
where $V= \frac{|\nablag \f|_g^2}{4}$. Note that most of the results we obtain are of interest for this particular question as well.

\subsubsection{A general lower bound} 
\label{s:sect-general}
The first family of results in this paper concern the general setting of a compact connected Riemannian manifold $(\M,g)$, with or without boundary $\d\M$, and the associated internal/boundary observability question. Namely, we consider the parabolic--transport problem with small viscosity $\eps>0$ and Dirichlet boundary conditions:
  \begin{equation}
  \label{e:transport-viscous-boundary}
 \left\{
 \begin{aligned}
&(\d_t - X -\q - \eps \Delta_g ) u = 0 , & \text{ in } (0,T)\times\Int(\M) ,  \\
& u = 0 & \text{ on } (0,T)\times \d \M ,  \\
& u|_{t=0} = u_0, & \text{ in } \M .
 \end{aligned}
 \right.
 \end{equation}
 Moreover, we assume that the vector field $X$ is a gradient vector field for the metric $g$, that is: there is $\f \in W^{2,\infty}(\M)$ (at least) such that 
 $$
 X = \nabla_g \f \cdot \nabla_g .
 $$
For the Dirichlet problem~\eqref{e:transport-viscous-boundary} as well as for the case $\d\M =\emptyset$ discussed in Section~\ref{s:intro-intro}, one may discuss the behavior of the observability constant (and in particular its possible uniform boundedness in the limit $\eps\to 0^+$) in the internal observability inequality~\eqref{e:transport-viscous-obs}. 
Also, a boundary observability problem for~\eqref{e:transport-viscous-boundary} can be formulated as follows (see~\cite{GL:07} and Section~\ref{s:boundary-shit} below). Given $\theta \in C^\infty(\d\M)$, there exist a constant $C_0(T,\eps)>0$ such that 
 \begin{multline}
 \label{e:transport-viscous-obs-boundary}
C_0(T,\eps)^2 \int_0^T\nor{ \theta \eps  \d_{\nu} u|_{\d\M}(t)}{H^{1/2}(\d\M)}^2  dt \geq \|u (T)\|_{H^1_0(\M)}^2 , \\
 \text{ for all } u_0 \in H^1_0(M) \text{ and } u \text{ solution of~\eqref{e:transport-viscous-boundary}}.
 \end{multline}
  Here $\d_\nu$ denotes a unit normal (for the metric $g$) vector field to $\d\M$.
 Then, the uniform observability question is whether $C_0(T,\eps)$ remains uniformly bounded in the limit $\eps\to 0^+$, and the associated minimal uniform observation time is defined again by 
 \begin{multline*}
T_{unif}(\theta)  = \inf \left\{ T >0 \text{ for which there exist }K_0, \eps_0>0 \text{ such that~\eqref{e:transport-viscous-obs-boundary} holds }\right. \\
 \left.  \quad  \text { for all }u_0\in H^1_0(\M) , \text{ with }C_0(T,\eps)\leq K_0 \text{ for all } \eps \in (0,\eps_0) \right\} .
\end{multline*}
Before going further, let us first give the definition of an analogue of the condition (GCC) in case of a manifold with boundary $\d\M \neq \emptyset$ (called here Flushing Condition), as used in the Guerrero-Lebeau result~\cite{GL:07}. For this, we need to extend $(\M,g)$ in a slightly bigger Riemannian manifold $(\tilde{\M}, \tilde{g})$, i.e. such that $\M\subset \Int(\tilde{\M})$ and $\tilde{g}$ any Riemannian metric on $\tilde{\M}$ such that $\tilde{g}|_{\M}=g$. In the case of a bounded domain of $\R^n$, one may choose $\tilde{\M} = \R^n$. We also extend $\f \in W^{2,\infty}(\M)$ as a {\em compactly supported } function $\tilde{\f}\in  W^{2,\infty}(\tilde{\M})$ such that $\tilde{\f}|_{\M} = \f$. 
\begin{definition} 
\label{d:def-GCC_bdry}
For $x\in \M$, denote by $\gamma_x$ the maximal solution to 
$$\dot{\gamma}_x(t) = -\nabla_{\tilde{g}}\tilde{\f}(\gamma_x(t)) , \quad \gamma_x(0)=x . $$
Note that this solution is defined globally in time since $\tilde{\f}$ is compactly supported.

We say that $(\M,\nabla_g \f,T)$ (resp. $(\M,\nabla_g \f,\omega,T)$) satisfies the Flushing Condition (FC) if for all $x\in \M$ there is $t\in (0,T)$ such that $\gamma_x(t)\notin \M$ (resp. $\gamma_x(t)\notin \M$ or $\gamma_x(t)\in \omega$). We also say that $(\M,\nabla_g \f)$ (resp. $(\M,\nabla_g \f,\omega)$) satisfies (FC) if there is $T>0$ such that $(\M,\nabla_g \f,T)$ (resp. $(\M,\nabla_g \f,\omega,T)$) does. We then define accordingly the time $T_{FC}(\M,\nabla_g \f)$ (resp.  $T_{FC}(\M,\omega,\nabla_g \f)$) as the infimum of $T>0$ for which this property holds. 

Remark finally that these definitions do not depend on the extensions $(\tilde{\M}, \tilde{g})$ and~$\tilde{\f}$.
\end{definition}
Guerrero-Lebeau~\cite{GL:07} prove an analogue of Theorem~\ref{t:Leb-Guer} in the present setting (and for general vector fields), namely: if $(\M,\nabla_g \f)$ (resp. $(\M,\nabla_g \f,\omega)$) satisfies the Flushing Condition, given $\theta$ any nontrivial observation function, then there is $T_{u}>0$ and $K_0>0$ such that for all $T\geq T_{u}$, \eqref{e:transport-viscous-obs-boundary} holds (resp.~\eqref{e:transport-viscous-obs} holds for all solutions to~\eqref{e:transport-viscous-boundary}), with $C_0(T,\eps)\leq K_0$ for all $\eps \leq 1$.

\bigskip
Two important geometric quantities in our results are the potential associated to the function $\f$, defined by 
\baln
\label{e:def:V}
V(x) := \frac{| \nabla_g \f(x)|_g^2}{4}  ,
\naln
and the Agmon distance (see e.g.~\cite[Chapter~3]{Helffer:booksemiclassic}) to the minimum of this potential, namely, with $E_0=\min_\M V$, 
\baln
\label{defAgmon-0}
d_A(x,y) &= \inf \left\{ \int_0^1\sqrt{\left(V(\gamma(t)) - E_0 \right)_+} |\dot{\gamma}(t)|_g  dt,\  \gamma \in \mathbf{U}_1(x,y)  \right\} ,\nonumber \\
\mathbf{U}_1(x,y) & =  \left\{ \gamma \in W^{1,\infty} ([0,1] ; \M ), \gamma(0)=x , \gamma(1)= y \right\} , \nonumber \\
d_{A}(x) & = d_{A}(x,V^{-1}(E_0)) = \inf_{y \in V^{-1}(E_0)} d_A(x,y) .
\naln
Here $V^{-1}(E_0)$ is the classically allowed region at the potential minimum, $$\left(V(x) - E_0 \right)_+ = \max \left(V(x) - E_0 , 0\right),$$ and $d_{A}(x)$ is the Agmon distance of $x$ to the set $V^{-1}(E_0)$ for the (pseudo-)metric $(V-E_0)_+ g$. Remark that the index $(\cdot)_+$ is not needed at the bottom energy $E_0=\min_\M V$; however, we keep it here since the definition~\eqref{defAgmon-0} will also be useful for a general energy level.

Our main result in this general setting formulates as follows.
\begin{theorem}
\label{thmtmpsminomega}
We assume that $\f \in W^{2,\infty}(\M)$ (hence $V \in W^{1,\infty}(\M)$) and $q \in L^\infty(\M)$. 
We let $E_0=\min_\M V$, set 
$$
W_{E_0}(x) = \frac{\f(x)}{2} + d_{A}(x)  ,
$$
and fix $\omega \subset \M$ (resp. $\theta \in C^\infty(\d \M)$), and in the case of {\em boundary estimates}, we further assume $\f,q \in C^\infty(\M)$. For any $\delta >0$, there is $\eps_0>0$ such that for all $\eps\in (0,\eps_0)$ the observability inequality~\eqref{e:transport-viscous-obs} (resp.~\eqref{e:transport-viscous-obs-boundary}) with constant $C_0(T,\eps)$ implies  
\bal
C_0(T, \eps) & \geq  \exp \frac{1}{\eps} \left( \min_{\ovl{\omega}}W_{E_0} - \max_{V^{-1}(E_0)}  \frac{\f}{2}- E_0 T -\delta \right)  , \\
 \text{resp. }C_0(T, \eps) & \geq  \exp \frac{1}{\eps} \left( \min_{\supp \theta}W_{E_0} -  \max_{V^{-1}(E_0)}\frac{\f}{2}- E_0 T-\delta \right)  . 
\nal
In particular, we have
\begin{align}
\label{e:estiminf-Tunifomega}
E_ 0 T_{unif}(\omega)\geq\min_{\ovl{\omega}}W_{E_0}-\max_{V^{-1}(E_0)}\frac{\f}{2} ,\\
\label{e:estiminf-Tunifomega-bord}
( \text{resp. } E_ 0 T_{unif}(\theta) \geq \min_{\supp\theta}W_{E_0}  - \max_{V^{-1}(E_0)}  \frac{\f}{2} ).
\end{align}
\end{theorem}
Note that the quantity in the right hand-side of \eqref{e:estiminf-Tunifomega}-\eqref{e:estiminf-Tunifomega-bord} as well as $E_0$ are invariant under the change $\f \to \f + C$ for $C$ constant on $\M$. This is consistent with the fact that the equations remain unchanged by such a modification of $\f$. 
Note also that if $E_0=0$ and $V^{-1}(E_0) \cap \omega = \emptyset$ (resp. $V^{-1}(E_0) \cap \supp(\theta) = \emptyset$), a precised version of this result (see Theorem~\ref{thm-general} below) actually shows that $C_0(T,\eps) \geq e^{c/\eps}$ for one $c>0$ and all time $T>0$ (in particular, uniform observability never holds).
This is consistent with (and a particular case of) the Guerrero-Lebeau~\cite{GL:07} result (first part of Theorem~\ref{t:Leb-Guer} above) for in this case, $(\M,\nabla_g \f,\omega)$ does not satisfy (GCC). Indeed, a point $x_0 \in V^{-1}(E_0)$ satisfies $\nabla_g \f(x_0)=0$ and is thus a stationary point of the gradient dynamics.

We refer to Remark~\ref{r:Cinfty-reg-bdry} concerning the additional smoothness assumption for the boundary estimate.

\begin{theorem}
\label{thm:counterexamplesomega}
Assume $\M = \ovl{\Omega}$ where $\Omega \subset \R^n$ is any smooth bounded connected open set endowed with $g=\Eucl$ the Euclidean metric and $q \in L^\infty(\M)$. There exist $\omega\subset \Omega$ and constants $c_\omega,\delta>0$ such that for any $\lambda>0$, there is a function $\f_\lambda \in C^\infty(\ovl{\Omega})$ such that: 
\begin{itemize}
\item  $(\ovl{\Omega},\nabla \f_\lambda,\omega)$ satisfies (FC) and $T_{FC}(\ovl{\Omega},\nabla \f_\lambda,\omega) \leq \diam(\Omega)$;
\item $T_{unif}(\omega)\geq c_\omega \lambda$, 
\item $ \lambda^2\delta\leq \nor{\nabla \f_\lambda}{L^\infty(\Omega)}^2  \leq \lambda^2 \diam(\Omega)^{2}+n$. 
\end{itemize}
In particular, for all $\Lambda >0$, there is $\f \in C^\infty(\ovl{\Omega})$ such that $$T_{unif}(\omega) \geq \Lambda T_{FC}(\ovl{\Omega},\nabla \f,\omega).$$
\end{theorem}
The result of Theorem~\ref{thm:counterexamplesomega} is already of interest in dimension one. In this case $\Omega = (-L,L)$, the vector field we consider is $\f'(x)\d_x$ with $\f'\geq 1$ on $[-L,L]$ and $\f'(0)=1$ and the observation set $\omega$ is a neighborhood of the boundary point $L$ (note that this would correspond to the case $M<0$ in the Coron-Guerrero problem with the notation of~\cite{CG:05}).
Note that the function $\f_\lambda$ in this result satisfies $\max \f_\lambda - \min \f_\lambda \simeq \lambda$.
As a consequence, one cannot even hope to have existence of a constant $\mathfrak{K}>0$ depending {\em only} on $\min_{x\in [0,L]}\f'(x)$ (a uniform flushing time) such that $T_{unif}(\omega) \leq \mathfrak{K} \ T_{FC}(\ovl{\Omega},\nabla \f_\lambda,\omega)$.
However, at this point, it does not seem hopeless that such a constant $\mathfrak{K}$ depends only on $\nor{\nabla \f_\lambda}{L^\infty(\Omega)}$, at least for a fixed metric.

\begin{remark}
In the case $\d \M = \emptyset$, Theorem~\ref{thmtmpsminomega} does not seem to suffice to construct functions $\f, \omega$ having $\frac{T_{unif}(\omega)}{T_{GCC}(\M,\nablag\f,\omega)}$ arbitrarily large. In a domain of $\R^n$, Theorem~\ref{thmtmpsminomega} is however enough to provide counterexamples.

Another drawback of Theorem~\ref{thmtmpsminomega} is that it does not produce any useful lower bound in case $\omega$ is a whole neighborhood of $\d\M$ (or in the boundary observation case from the whole boundary $\d \M$).

We remedy these issues in the next section on surfaces of revolution.
\end{remark}

\subsubsection{Lower bounds on surfaces of revolution}
\label{s:revol-intro}
In Theorems~\ref{thmtmpsminomega} and~\ref{thm:counterexamplesomega} above, the lower bound of the minimal uniform observability time is essentially due to the contribution of the potential $V(x) = \frac{|\nabla_g \f(x)|_g^2}{4}$ (and the difference between its maximal and minimal values on $\M$). 
In this section, we consider a family of geometric settings, namely surfaces of revolution, for which the contribution of the geometry of $(\M,g)$ plays an important role. 
This leads in particular to explicit versions of Theorem~\ref{t:loose}.

The precise description of the geometry of the surfaces we consider is given in Section~\ref{s:revol} and we only describe here features required to state the result.
We may consider either:
\begin{enumerate}
\item $\M =\calS \subset \R^3$ a smooth compact surface diffeomorphic to the sphere $\S^2$;
\item $\M =\calS \subset \R^3$ a smooth compact surface diffeomorphic to the disk $\ID$;
\item $\M =\calS \subset \R^3$ a smooth compact surface diffeomorphic to the cylinder $[0,1]\times \S^1$;
\item $\M =\calS$ a smooth compact surface diffeomorphic to the torus $\T^2 = \S^1\times \S^1$. 
\end{enumerate}
We assume moreover that it has revolution invariance around an axis. In particular, the axis may intersect $\calS$ in two points (sphere), one point (disk) or no point (cylinder or torus). Except near these points, $\calS$ has a global coordinate chart $(s,\theta)\in(0,L)\times \S^1$ for some $L>0$. 
In the first three cases, the surface is endowed with the metric $g$ inherited from the Euclidean metric on $\R^3$ which writes
\baln
\label{e:def-metric-R}
g = ds^2  +R(s)^2 d\theta^2  ,
\naln
on account to the rotation invariance. 
Here the function $R:(0,L)\to (0,\infty)$ describes the shape of $\calS$ (distance to the revolution axis).
In the torus case, we simply endow $\T^2$ with the metric~\eqref{e:def-metric-R}.

We further assume that the function $\f$ and $q$ are themselves rotationally invariant, that is $\f=\f(s)$ and $q=q(s)$ in these coordinates.
In this setting (and as opposed to results presented in above Section~\ref{s:sect-general}), the relevant quantities for our analysis are the following.

We define for any $\ce>0$ (that can be chosen) the ($\theta$-invariant) effective potential
\baln
\label{e:def-Vce}
V_\ce (s) := \frac{\ce^2}{R(s)^2} + \frac{| \f'(s)|^2}{4} .
\naln
Note that, as opposed to the potential appearing in~\eqref{e:def:V}, this potential $V_\ce$ is different from $\frac{|\nabla_g\f|_g^2}{4} =  \frac{| \f'(s)|^2}{4}$.
Moreover, it depends explicitly on the geometry (namely, on $R$).
We shall make the simplifying assumption that 
\baln
\label{s:asspt-uniq-min}
V_{\ce}^{-1}(\min V_\ce) = \{s_{\min}\} \text{ consists in a single point }s_{\min} .
\naln
Note that in case $\calS$ has a boundary, one may have $s_{\min}$ at the boundary (see Section~\ref{s:revolution} for more precise statements). 
As in the previous section, we introduce the associated Agmon distance, which simply writes in the first three cases:
\baln
\label{e:defbisdA}
d_{A}^{\ce}(s) = \left| \int_{s_{\min}}^s \sqrt{ V_\ce(y)- V_\ce(s_{\min})} dy \right|  .
\naln
In the fourth case $\calS = \S^1_L\times \S^1$, an analog of \eqref{defAgmon-0} still makes sense on $\S^1_L$ when $V_{\ce}$ is defined on $\S^1_L$. We may choose a representation of $\S^1_L=\R/L\Z$ in which $s_{\min}= 0 + L \Z$, a definition of the Agmon distance then reads
\baln
\label{e:defbisdAT2}
d_{A}^{\ce}(s) = \min_{s_{\min}\in L \Z}\left[  \left| \int_{s_{\min}}^s \sqrt{ V_\ce(y)- V_\ce(s_{\min})} dy \right|\right].
\naln
We also set
\baln
\label{e:def-Wce}
W^{\ce}(s)=d_A^{\ce}(s)+\frac{\f(s)}{2} ,\quad \text{for } s \in (0,L) .
\naln
Then, our main result in this geometric context can be (loosely) stated as follows.
\begin{theorem}
\label{thmtmpsmin}
Let $\ce>0$, assume that $V_\ce$ satisfies~\eqref{s:asspt-uniq-min} and that, in the coordinates $(s,\theta)\in (0,L)\times \S^1$, we have $\omega = (0,\delta)\times \S^1 \cup  (L-\delta,L)\times \S^1$ (this can be rewritten in an intrinsic way on $\calS$). 
Assuming the observability inequality~\eqref{e:transport-viscous-obs} (resp.~\eqref{e:transport-viscous-obs-boundary}) with constant $C_0(T,\eps, \omega)$ (resp. $C_0(T,\eps,  \{L\}\times \S^1)$), 
there is a sequence $\eps_k \to 0^+$ such that  for any $\delta>0$, there is $k_0(\delta)>0$ such that 
\bal
C_0(T,\eps_k, \omega)  \geq e^{\frac{1}{\eps_k}\big(W^{\ce}_\omega-W^{\ce}_m-V_{\ce}(s_{\min})T-\delta\big)} \quad \text{ for all } k\geq k_0(\delta) ,\quad W^{\ce}_\omega = \min_{\bar{\omega}}W^{\ce},\\
\text{resp.}\quad C_0(T,\eps_k,  \{L\}\times \S^1)  \geq e^{\frac{1}{\eps_k}\big(W^{\ce}(L)-W^{\ce}_m-V_{\ce}(s_{\min})T-\delta\big)} \quad \text{ for all } k\geq k_0(\delta) ,
\nal
where $W^{\ce}_m = \inf_{(0,L)}W^{\ce}$.
In particular, we have 
\baln
\label{e:estiminf-Tunif}
V_{\ce}(s_{\min})T_{unif}(\omega)& \geq W^{\ce}_{\omega}- W^{\ce}_m  , \\
\label{e:estiminf-Tunif-bord}
V_{\ce}(s_{\min})T_{unif}(\{L\}\times \S^1)&\geq W^{\ce}(L)- W^{\ce}_m   .
\naln
\end{theorem}

Theorem~\ref{thmtmpsmin} differs from Theorem~\ref{thmtmpsminomega} in several respects. First notice that the potential appearing in Theorem~\ref{thmtmpsminomega} is $\frac{| \f'(s)|^2}{4}$, that is $V_0(s)$ with the definition of $V_{\ce}$ in~\eqref{e:def-Vce}. In particular, it does not depend on $R$: neither does its minimal value, nor the associated Agmon distance and function $W^0$.
Therefore, in this very particular geometric context, the results of Theorem~\ref{thmtmpsminomega} do not depend on the geometry of $R$, and hence only formulate as a one dimensional result in the $s$ variable. As such, they do not care about the ``transverse dynamics'' in the $\theta$-variable. Theorem~\ref{thmtmpsmin} overcomes this lack and shows that both have to be taken into account.

Another difference with the estimates of Theorem~\ref{thmtmpsminomega} is that $- \max_{V^{-1}(E_0)}\frac{\f}{2} = - \max_{V^{-1}(E_0)}W$ is here replaced by $- \inf_{(0,L)}W^{\ce} $. This improvement is due to the ``one dimensional'' underlying framework (in which localization properties of eigenfunctions are better understood).

Again, we remark that the initial problem is invariant by changing $\f$ by a constant $\f+C_0$. In Estimate~\eqref{e:estiminf-Tunif}, both the potential $V_{c}$ and the quantity $W_{\omega}^{\ce}-W_{m}^{\ce}$ are as well invariant by this change of the function $\f$.

We now state three particular examples of application of Theorem~\ref{thmtmpsmin}. 
The latter imply Theorem~\ref{t:loose}.

\begin{corollary}
\label{c:sphere-disk}
Assume $\calS$  is a surface of revolution in $\R^3$ diffeomorphic to $\S^2$ (resp. $\ID$) and such that $R^{-1}(\max R)$ is a single point ($R$ has a unique maximum). Denote by $N,S$ (resp. $N$ only) the north/south poles of $\calS$, which are the only two (resp. the unique) invariant points under the revolution symmetry.
Then, for any $\delta>0$, there exists $\f_\delta  \in C^\infty(\calS)$ invariant by rotation such that
with $\omega_\delta = B_g(N, \delta) \cup B_g(S, \delta)$ (resp. $\omega_\delta = B_g(N, \delta)$) we have
\begin{enumerate}
\item $(\calS,\nabla_g\f_\de,\omega_\delta)$ satisfies (GCC) (resp. satisfies (FC)) and $T_{GCC}(\calS,\nabla_g\f_\de,\omega_\delta) = L-2\delta \leq L= \dist_g(N,S)$ (resp. $T_{FC}(\calS,\nabla_g\f_\de,\omega_\delta) = L-\delta \leq L= \dist_g(N,\d \calS)$);
\item For all $\ce>0$, there is $\delta_0,C >0$ such that for all $\delta \in (0,\delta_0)$ 
$$
T_{unif}(\omega_\delta) \geq  \left( \frac{\ce^2}{R(s_{\min})^2} + \frac14 \right)^{-1} \left(   \ce \log(\frac1{\delta}) - C \right)  .
$$
\end{enumerate}
\end{corollary}
This result proves the first item in Theorem~\ref{t:loose}.
In particular, the limit $\delta\to 0^+$ prevents from the existence of a universal constant $\mathfrak{K}\geq 1$ such that $T_{unif}(\omega) \leq \mathfrak{K} \ T_{GCC}(\M,X,\omega)$.
Note that in this construction, the vector fields $\nabla_g\f_\de \cdot \nabla_g$ are rotationally invariant and independent of the metric $g$. Moreover, for $\delta < \delta'$, the two functions $\f_\de$ and $\f_{\de'}$ coincide on $\M\setminus \omega_{\delta'}$.

In our second result, the geometry is close to that of Corollary~\ref{c:sphere-disk}. However, we consider fixed $\omega$ and $\f$ (and even a fixed vector field), but let the metrics $g$ vary.
We denote $\S^1_L = \R/L\Z$ and $\S^1 =\S^1_{2\pi}$.
\begin{corollary}
\label{c:example-torus}
Assume $\calS=\S^1_L \times \S^1$ with coordinates $(s,\theta)$ and let $\f\in C^{\infty}(\S^1_L)$. 
Let $\omega = I_\omega \times \S^1$ with $I_\omega$ a nonempty interval such that $\ovl{I_\omega} \neq \S^1_L$.

Then, there is a constant $C>0$ such that for any $\delta \in (0,1)$, there exists a function $R_\delta \in C^\infty(\S^1_L ;\R^+_*)$ such that 
\begin{enumerate}
\item the vector field $\nabla_{g_\delta}\f \cdot \nabla_{g_\delta}= \f'(s) \d_s$ (defined by \eqref{e:def-metric-R} associated to $R_\delta$) does not depend on the metric $g_\delta$, the triple $(\calS,\nabla_{g_\delta}\f,\omega)$ satisfies (GCC) if and only if $\f'\neq 0$ on $\S^1_L \setminus I_\omega$, and, in this case, we have
$$
T_{GCC}(\calS , \nabla_{g_\delta}\f ,I_\omega \times \S^1) =  \left| \int_{\S^1_L \setminus I_\omega} \frac{ds}{\f'(s)} \right| < +\infty  ,
$$
\item for the transport equation~\eqref{e:transport-viscous-gradient} with viscosity given by the metric $g_\delta$ under the form \eqref{e:def-metric-R} associated to $R_\delta$, we have
$$
T_{unif}(I_\omega \times\S^1) \geq C\delta^{-1/2} ,
$$
\item  $\delta^{1/2} (1-C\delta) \leq \min_{\S^1}R_\delta \leq \delta^{1/2}$ for all $\delta \in (0,1)$.
\end{enumerate}
\end{corollary}
This result implies the second item in Theorem~\ref{t:loose}.

Another application is given by the following result, which is an analogue of Corollary~\ref{c:example-torus} for the boundary observability problem in the cylindrical geometry.
\begin{corollary}
\label{c:example-cylinder}
Assume $\calS=[0,L]\times\S^1$ (whence $\d \calS = (\{0\}\cup\{L\})\times\S^1$), and let $\f\in C^{\infty}([0,L])$. Then, for any $\gamma >2$, $\delta_0>0$ there is a constant $C>0$ such that for all $\delta \in (0,\delta_0]$, there exists a function $R_\delta: [0,L]\to\R^+_*$ such that 
\begin{enumerate}
\item  the vector field $X= \nabla_{g_\delta}\f \cdot \nabla_{g_\delta} = \f'(s) \d_s$ does not depend on the metric $g_\delta$ (defined by~\eqref{e:def-metric-R} with $R_\delta$);
\item $(\calS,\nabla_{g_\delta}\f)$ satisfies (FC) (in the sense of Definition~\ref{d:def-GCC_bdry}) if and only if $\f' \neq 0$ on $[0,L]$ and in this case,
$$
T_{FC}(\M,\nabla_{g_\delta}\f) =  \left| \int_0^L \frac{ds}{\f'(s)} \right| <+\infty ;
$$
\item for the transport equation with viscosity and Dirichlet boundary conditions~\eqref{e:transport-viscous-gradient}, and with metric $g_\delta$ under the form~\eqref{e:def-metric-R} associated to $R_\delta$, we have
$$
T_{unif}(\d\calS) \geq \left(\max_{[0,L]}\frac{|\f'|^2}{4} \right)^{-1} \left(  \frac{1}{\delta^{\gamma/2-1}}\frac{1}{\gamma/2-1}  - C \right) .
$$
\item \label{itemdelta}$(s+\delta)^{\frac{\gamma}{2}} (1-C (s+\delta)) \leq  R_\delta(s) \leq  (s+\delta)^{\frac{\gamma}{2}}$ for $s\in [0,L/4]$ and $(L-s+\delta)^{\frac{\gamma}{2}} (1-C (L-s+\delta)) \leq  R_\delta(s) \leq  (L-s+\delta)^{\frac{\gamma}{2}}$ for $s\in [3L/4,L]$; 
\end{enumerate}
\end{corollary}
Note the link between the asymptotic singularity of the metric $R_\delta(s) \sim (s+\delta)^{\frac{\gamma}{2}}$ (i.e. $\calS$ becomes close to a ``conical geometry'' for small $\delta$) and the blowup of the minimal time $T_{unif} \gtrsim\frac{1}{\delta^{\gamma/2-1}}$. Note also that the limit case $\gamma =2$, all calculations can be done as well and lead to $R_\delta(s) \sim s+\delta$ together with $T_{unif}\gtrsim -\log \delta$. This corresponds to the case where the geometry of the cylinder degenerates towards that of the disk, leading to the same blowup estimate as in Corollary~\ref{c:sphere-disk}.

\subsubsection{Observability for positive solutions}
As already mentioned, our last result concerns the uniform observability question for positive solutions of~\eqref{e:transport-viscous-gradient}, and is restricted to the case $\d\M=\emptyset$.  We also assume $\f\in C^3$.
Note that if $u_0 \in L^2(\M; \R^+)$, then the associated solution $u$ to~\eqref{e:transport-viscous-gradient} satisfies $u(t,x) \geq 0$ for a.e. $(t,x) \in \R^+\times\M$ (see e.g. Theorem~9 in Section~7.1 p369 together with Problem~7 in Section~7.5 in~\cite{Evans:98}, or Chapter~III, Theorem 7.1, p181 in~\cite{LSU:68}).

We consider the observability inequality for nonnegative solutions:
 \begin{multline}
 \label{e:transport-viscous-obsplus}
C_0^+(T,\eps)^2 \int_0^T\int_\omega |u(t,x)|^2 ds(x)dt \geq \|u (T)\|_{L^2(\M)}^2 , \\
 \text{ for all } u_0 \in L^2(\M; \R^+)\text{ and } u \text{ solution of~\eqref{e:transport-viscous-gradient}}.
 \end{multline}
and the associated minimal time $T_{unif}^+(\omega)$ of uniform observability for positive solution, already defined in~\eqref{defTunif+}. The main result we obtain in this context is the above Theorem~\ref{t:positive-intro}, stating that $T_{unif}^+(\omega)=T_{GCC}(\M, \nabla_g\f,\omega)$.
As a byproduct of our analysis, we also obtain a lower estimate on the blow up of the control cost when the Geometric Control Condition is not satisfied. It involves the definition of a quantity that roughly speaking, measures how two points are far from being the image of a trajectory at time $t$, namely
\begin{align}
\label{defdXintro} d_{\nabla_g\f}(x,y,t)&=\frac{1}{4} \inf \left\{ \int_0^t \left|\dot{\gamma}(s)-\nabla_g \f(\gamma(s))\right|_g^{2} ds , \ \gamma \in \mathbf{U}_t(x,y)
\right\} , \\
 \mathbf{U}_t(x,y) & = \left\{\gamma \in W^{1,\infty} ([0,t] ; \M ), \gamma(0)=x , \gamma(t)= y \right\} . \nonumber
\end{align}
Note that we have $ d_{\nabla_g\f}(x,y,t)\approx d(x,\phi_{-t}(y))$ for bounded $t$, where $d$ denotes the Riemannian distance and $\phi_t$ the flow of the vector field $\nabla_g \f$ (see Lemma~\ref{lmdist} where $d_{\nabla_g\f}(x,y,t)$ is interpreted as a control cost from $x$ to $\phi_{-t}(y)$ with time varying metric).
\begin{proposition}
\label{proplowerobse+}
Assume that $(\M,\nabla_g \f,\overline{\omega})$ does not satisfy (GCC). Then, we have
\bal
d_{([0,T], \overline{\omega})} : = \sup_{y\in \M}\inf_{x\in\overline{\omega},t\in[0,T]}d_{\nabla_g\f}(x,y,t)>0.
\nal
Moreover, for any $\delta>0$, there is $\e_0>0$ such that for all $\eps \in (0,\eps_0)$, we have
\baln
C_{0}(T,\e)\geq C_{0}^{+}(T,\e)\geq e^{\frac{d_{([0,T], \overline{\omega})}-\delta}{\e}}\label{lowercost+}.
\naln
\end{proposition}
Note that this exponential blowup is a refinement of the Guerrero-Lebeau~\cite{GL:07} result (first part of Theorem~\ref{t:Leb-Guer} above). However, we provide here with a precise geometric rate (namely $d_{([0,T], \overline{\omega})}$) quantifying this blowup phenomenon.

The proofs of Theorem~\ref{t:positive-intro} and Proposition~\ref{proplowerobse+} rely on estimates on the kernel of the associated equation. Note that kernel estimates have already been used in~\cite{Miller:04} to prove lower bounds for the cost of controllability of the usual heat equation in the short time asymptotics, and in~\cite{LL:18} to prove observability of positive solutions to the heat equation with optimal constants.

\subsection{Further remarks}
In this section, we collect several remarks and comments related to our results. 
\begin{enumerate}
\item The principal interest of working with gradient vector fields $X=\nablag\f$ is that the associated operator $-\nablag\f-\eps\Delta_g$ can be conjugated to a {\em selfadjoint} Sch\"odinger operator~\eqref{e:intro-witten}. 
And the limit $\eps\to 0^+$ then corresponds to the semiclassical limit, which has been the object of many studies (see e.g.~\cite{Witten:82,Simon:83,HS:84,HS:85,Helffer:booksemiclassic,CFKS,Allibert:98,HKN:04,HN:06}.
This conjugation does no longer hold in case $X$ is not a gradient vector field.
One could also consider that giving counter-examples with gradient flows is ``stronger'' than general counter-examples. 
We do not know wether an analogue of Theorem~\ref{t:positive-intro} for positive solutions remains true for general vector fields $X$. This seems to be an interesting open problem.

\item In the context of surfaces of revolution, as presented in Section~\ref{s:revol-intro}, we only provide with lower bounds of $T_{unif}$. It would of course be interesting to obtain related upper bounds on this uniform control time. This would require to provide a very precise description of several spectral quantities (spectral gaps, localization of eigenfunctions at all energy levels...) for the semiclassical Schr\"odinger operator $P_\eps$ in~\eqref{e:intro-witten}, and seems to be a difficult question.
See the companion paper~\cite{LL:20-1D} for an upper bound in a related one dimensional situation.

\item The one dimensional one well problem is considered in~\cite{LL:20-1D}. In this very particular situation, we are able to provide with 
\begin{itemize}
\item improved lower bounds on the minimal time when compared to Theorem~\ref{thmtmpsminomega};
\item an upper bound on the minimal time.
\end{itemize}
This requires the knowledge of precise information on the spectral gap and the localization of eigenfunctions at all energy levels $E\geq E_0$ (whereas Theorem~\ref{thmtmpsminomega} is only concerned with the bottom energy level $E_0$). See also the discussion at the beginning of Section~\ref{s:general-bounds}.

\item Notice that if one is not interested in null-controllability (i.e. driving the solution exactly to zero at time $T$), but rather in approximate controllability with a reasonable cost (and a precision depending on the viscosity $\eps$), one might be satisfied by the following statement. 
\begin{proposition}
\label{p:guerrero-lebeau-control}
Suppose $\d\M= \emptyset$ and $(\M, X,\omega,T)$ satisfies (GCC) (resp. $\d\M \neq \emptyset$ and $(\M, X,\omega,T)$ satisfies (FC)). Then, there exist $C,C_0>0$ such that for all $y_0\in L^2(\M)$, $\eps \in (0,1]$ there is $h=h_\eps \in L^2((0,T)\times \omega)$ with 
$$
\nor{h_\eps}{L^2((0,T)\times \omega)} \leq C \nor{y_0}{L^2(\M)} ,
$$
such that the associated solution to~\eqref{e:transport-viscous-control-internal} satisfies
$$
\nor{y(T)}{L^2(\M)} \leq Ce^{-\frac{C_{0}}{\e}}  \nor{y_0}{L^2(\M)} .
$$
\end{proposition}
That is to say, one can drive the solutions $e^{-\frac{C_{0}}{\e}} $ close to zero with a uniformly bounded cost.
This result follows from Proposition~\ref{p:guerrero-lebeau} below (a particular case of \cite[Proposition~3]{GL:07}) together with~\cite[Appendix]{LL:17approx}.
This can be particularly useful for numerical purposes, since $e^{-\frac{C_{0}}{\e}} = 0$ numerically for $\eps$ small enough.

In the situation of Theorems~\ref{thmtmpsminomega} or~\ref{thmtmpsmin}, this means that for intermediate times $T \in \big(T_{GCC}(\M , \nabla_g\f ,\omega), T_{unif}(\omega) \big)$ (resp. for $T \in \big(T_{FC}(\M , \nabla_g\f ,\omega), T_{unif}(\omega) \big)$ if $\d\M \neq \emptyset$), controlling the solution $e^{-\frac{C_{0}}{\e}}$ close to zero costs $\approx 1$, whereas controlling the solution exactly to zero costs $\approx e^{\frac{C}{\eps}}$.

\item Note that in the context of revolution surfaces of Section~\ref{s:revol-intro}, we prove a complementary result compared to~\cite[Theorem~1.9/Corollary~1.10]{LL:18}. We prove in Proposition~\ref{propinfdA} that in all cases of Section~\ref{s:revol-intro}, for any rotationally invariant set $\omega$, we have (with the notation of~\cite{LL:18}, the constant $\mathfrak{K}_{eig}(\omega)$ being the smallest constant $\mathfrak{K}$ in the inequality $\nor{\psi_\lambda}{L^2(\M)} \leq Ce^{\mathfrak{K}\sqrt{\lambda}}\nor{\psi_\lambda}{L^2(\omega)}$ where $-\Delta_g \psi_\lambda = \lambda \psi_\lambda$)
$$
\mathfrak{K}_{eig}(\omega) = R(s_{\min}) d_A(\omega) , \quad \text{with } d_A(\omega) = \inf_{x\in \omega} d_A(x) .
$$
In~\cite{LL:18}, we only proved $\mathfrak{K}_{eig}(\omega) \geq R(s_{\min}) d_A(\omega)$ (and only in case $\calS$ is diffeomorphic to a sphere). This result is close to that of Allibert~\cite{Allibert:98}, which already proves this in case $\calS$ is diffeomorphic to a cylinder and the function $R$ has a single local maximum which is non-degenerate.
\end{enumerate}

\bigskip
\noindent
{\em Acknowledgements.} 
The authors are partially supported by the Agence Nationale de la Recherche under grants SRGI ANR-15-CE40-0018 (for C.L.), SALVE ANR-19-CE40-0004 (for M.L.) and ISDEEC ANR-16-CE40-0013 (for C.L. and M.L.).
Part of this research was done when the second author was in \'Ecole Polytechnique, Centre de Math\'ematiques Laurent Schwartz UMR7640,  91128 Palaiseau cedex France.

We would like to thank Olivier Glass (who suggested Proposition \ref{p:guerrero-lebeau-control}), Franck Sueur and Dorian Le~Peutrec for interesting discussions related to this work. 
We are also grateful to Jon Asier B\'{a}rcena-Petisco for his comments on a first version of the article. Finally we would like to thank both referees for their careful reading of the manuscript and their constructive remarks that helped to improve the presentation of the paper.

\section{Preliminaries: duality, conjugation of gradient flows and (GCC)}
\label{s:preliminaires}

\subsection{Uniform controllability problems and dual formulation}

In this section, we reformulate the above uniform observability questions in terms of uniform controllability statements. This follows essentially the classical duality approach, see~\cite{DR:77} or~\cite[Chapter~2.3]{Cor:book}.

\subsubsection{Duality between internal control and observation problems}
\label{s:obs-cont}
 In this section, we present the controllability problems associated to the above observability questions, and we briefly describe the duality between the control and the observation problems. 
 We introduce the internal control problem
    \begin{equation}
  \label{e:transport-viscous-control-internal}
 \left\{
 \begin{aligned}
&(\d_t  + X + \div_g(X) -\q - \eps \Delta_g ) y = \mathds{1}_\omega h , & \text{ in } (0,T)\times \M ,  \\
& y = 0, & \text{ on } (0,T)\times \d \M , \\
& y|_{t=0} = y_0, & \text{ on } \M .
 \end{aligned}
 \right.
 \end{equation}
Notice that, as opposed to~\eqref{e:transport-viscous}, the operator appearing in these control problems is $X^* = -X - \div_{g}(X)$, where the adjoint is taken in the space $L^2(\M, d\Vol_g)$.

That the appropriate dual observation problem is~\eqref{e:transport-viscous} is a consequence of the following lemma.
 \begin{lemma}[Duality equation]
 \label{l:duality-equation-1}
 For all solutions $u \in C^0([0,T]; L^2(\M))$ of \eqref{e:transport-viscous} on $[0,T]$ with Dirichlet boundary conditions and all $y\in C^0([0,T]; L^2(\M))$ solution to~\eqref{e:transport-viscous-control-internal} with $h \in L^2((0,T)\times \M)$, we have
 \begin{align}
 \label{e:duality-equation-1}
 (u(T) , y_0)_{L^2(\M)} - (u_0, y(T))_{L^2(\M)} +\int_0^T\big(\mathds{1}_\omega u(t), h(T-t) \big)_{L^2(\M)} dt =0 .
\end{align}
 \end{lemma}
 
 Notice that one passes from the observed evolution to the controlled evolution by changing $(X, q)$ into $(-X , q-\div_{g}(X))$. The interest of adding a potential term $q(x)$ in these equations is that the free equation~\eqref{e:transport-viscous-boundary} and the controlled equation~\eqref{e:transport-viscous-control-internal} then have the same form (i.e. the adjoint of a vector field is not a vector field but the adjoint of a first order differential operator is a first order differential operator).
 
 \begin{definition}[Controllability and cost]
  Given $(\omega, \eps, T)$, we say that~\eqref{e:transport-viscous-control-internal} is null-controllable from $(\omega, T)$ if for any $y_0 \in L^2(\M)$, there is $h \in L^2((0,T)\times \M)$ such that the associated solution to~\eqref{e:transport-viscous-control-internal} satisfies $y(T) =0$. 
  If~\eqref{e:transport-viscous-control-internal} is null-controllable from $(\omega, T)$, we define for $y_0 \in L^2(\M)$ the set $U(y_0) \neq \emptyset$ of all such controls $h \in L^2((0,T)\times \M)$, and the cost function 
  $$
  \mathcal{C}_0(T,\eps) := \sup_{y_0 \in L^2(\M) , \nor{y_0}{L^2(\M)}\leq 1} \left\{\inf_{h \in U(y_0)} \nor{h}{L^2((0,T)\times \M)} \right\} .
  $$
 \end{definition}
  As a corollary of Lemma~\ref{l:duality-equation-1}, together with classical arguments (see e.g.~\cite{DR:77} or~\cite[Chapter~2.3]{Cor:book}) we deduce the following statement.
 \begin{corollary}[Observability constant = control cost]
\label{c:obscontrol}
 Given $(\omega, \eps, T)$, Equation~\eqref{e:transport-viscous-control-internal} is null-controllable from $(\omega, T)$ if and only if the observability inequality~\eqref{e:transport-viscous-obs} holds. 
 Moreover, we then have $\mathcal{C}_0(T,\eps)= C_0(T,\eps)$.
  \end{corollary}

As a consequence, all lower bounds on $C_0(T,\eps)$ formulated in Theorems~\ref{thmtmpsminomega} and~\ref{thmtmpsmin} translate into lower bounds on $\calC_0(T,\eps)$. The time $T_{unif}(\omega)$ is equal to the minimal time of uniform {\em controllability}, and all lower bounds on the time $T_{unif}(\omega)$ obtained in Theorems~\ref{thmtmpsminomega} and~\ref{thmtmpsmin} and their corollaries apply.

\medskip
The uniform observability result for positive solutions of the heat equation in Theorem~\ref{t:positive-intro} also has a controllability counterpart. 
This fact was indeed proved by Le~Balc'H~\cite[Theorem 4.1]{L:18} for the classical heat equation. 
In the present context, the uniform observability estimate for positive solutions, associated to Theorem~\ref{t:positive-intro}, implies the following controllability result.
\begin{corollary}
\label{corcontrol+}
Let $\M$ be a compact Riemannian manifold with $\d\M=\emptyset$, $X=\nabla_{g} \f$ where $\f\in C^{3}(\M)$, and $\omega\subset \M$ an open subset. Assume that $(\M, \nabla_{g} \f,\omega)$ satisfies (GCC), and $T>T_{GCC}(\M, \nabla_{g} \f,\omega)$. Then, there exist $C,\e_0>0$ so that for any $y_{0}\in L^{2}(\M)$ and $0<\e\leq \e_{0}$, there exists a control $h\in L^2([0,T],L^2(\omega))$ with 
$$\nor{h}{L^2([0,T],L^2(\omega))}\leq C\nor{y_{0}}{L^{2}(\M)}$$
such that the solution of \eqref{e:transport-viscous-control-internal} satisfies $y(T)\geq0$.
\end{corollary}
We refer to Section~\ref{s:cont-pos} for a proof.
 
\subsubsection{Duality between boundary control and observation problems}
\label{s:boundary-shit}
We now briefly discuss the boundary case and we refer to~\cite{GL:07} for the details.
The boundary control problem under interest is 
     \begin{equation}
  \label{e:transport-viscous-control-boundary}
 \left\{
 \begin{aligned}
&(\d_t  + X + \div_{g}(X) -\q - \eps \Delta_g ) y = 0, & \text{ in } (0,T)\times \Int( \M) ,  \\
& y =  \theta h , & \text{ on } (0,T)\times \d \M , \\
& y|_{t=0} = y_0, & \text{ on } \M ,
 \end{aligned}
 \right.
 \end{equation}
 where $\theta  \in C^\infty(\d\M ; \R)$ is meant to be a smooth version of $\mathds{1}_\Gamma, \Gamma \subset \d \M$. 
Solutions of~\eqref{e:transport-viscous-control-boundary} are defined in the sense of transposition, and a well-posedness statement can be written as follows.

\begin{lemma}[Guerrero-Lebeau~\cite{GL:07} pp~1814-1815]
Assume $X$ is a $L^\infty$ vector field on $\M$ with $\div_{g}(X) \in L^\infty(\M)$, $q \in L^\infty(\M)$, and let $T>0$. Then, there exists $C>0$ such that for all $y_0 \in H^{-1}(\M)$, all $h \in L^2(0,T ; H^{-1/2}(\d \M))$ and all $\eps>0$, there is a unique solution $y$ to~\eqref{e:transport-viscous-control-boundary} in the sense of transposition, which satisfies $y \in L^2((0,T)\times \M) \cap  C^0([0,T]; H^{-1}(\M)) \cap H^1(0,T ; H^{-2}(\M))$ with 
\begin{multline*}
 \nor{y}{L^2((0,T)\times \M)} +  \nor{y}{L^\infty(0,T; H^{-1}(\M))} +  \nor{\d_t y}{L^2(0,T ; H^{-2}(\M))} \\
 \leq \frac{C}{\eps} \left( \nor{y_0}{H^{-1}(\M)}+ \eps^{1/2} \nor{h}{L^2(0,T ; H^{-1/2}(\d \M))}\right) .
\end{multline*}
\end{lemma}
Such solutions in particular solve the first equation of~\eqref{e:transport-viscous-control-boundary} in the sense of distributions.
 \begin{definition}[Controllability and cost]
  Given $(\theta, \eps, T)$, we say that~\eqref{e:transport-viscous-control-boundary} is null-controllable from $(\theta, T)$ if for any $y_0 \in H^{-1}(\M)$, there is $h \in L^2(0,T; H^{-1/2}(\d\M))$ such that the associated solution to~\eqref{e:transport-viscous-control-boundary} satisfies $y(T) =0$. 
  If~\eqref{e:transport-viscous-control-boundary} is null-controllable from $(\theta, T)$, we define for $y_0 \in H^{-1}(\M)$ the set $U(y_0) \neq \emptyset$ of all such controls $h \in L^2(0,T; H^{-1/2}(\d\M))$, and the cost function 
  $$
  \mathcal{C}_0(T,\eps) := \sup_{y_0 \in H^{-1}(\M) , \nor{y_0}{H^{-1}(\M)}\leq 1} \left\{\inf_{h \in U(y_0)} \nor{h}{L^2(0,T; H^{-1/2}(\d\M))} \right\} .
  $$
 \end{definition}

We now describe the link with the boundary observation problem~\eqref{e:transport-viscous-boundary}. 
 We start with the duality identity.
\begin{lemma}[Duality equation]
 For all solutions $u \in C^0([0,T]; H^1_0(\M))$ of \eqref{e:transport-viscous-boundary} on $[0,T]$ and all $y\in C^0([0,T]; H^{-1}(\M))$ solution to~\eqref{e:transport-viscous-control-internal} with $h \in L^2(0,T; H^{-1/2}(\M))$, we have
 \begin{multline*}
 \langle u(T) , y_0\rangle_{H^1_0 , H^{-1}} - \langle u_0, y(T)\rangle_{H^1_0 , H^{-1}} \\
  - \int_0^T \left\langle \theta \eps \d_{\nu}u|_{\d\M}(t), h(T-t) \right\rangle_{H^{1/2}(\d\M) ,H^{-1/2}(\d\M)} dt =0 .
\end{multline*}
 \end{lemma}
 The proof is omitted here and only consists in an integration by parts for smooth solutions, and then a density argument.
 As in the internal case, classical duality arguments (see~\cite{DR:77} or~\cite[Chapter~2.3]{Cor:book}) yield the following statement.
 
 \begin{corollary}[Observability constant = control cost]
 Given $(\theta, \eps, T)$, Equation~\eqref{e:transport-viscous-control-boundary} is null-controllable from $(\theta, T)$ if and only if the observability inequality~\eqref{e:transport-viscous-obs-boundary} holds. 
 Moreover, we then have $\mathcal{C}_0(T,\eps)= C_0(T,\eps)$.
 \end{corollary}
Again, all lower bounds on $C_0(T,\eps)$ formulated in Theorems~\ref{thmtmpsminomega} and~\ref{thmtmpsmin} translate into lower bounds on $\calC_0(T,\eps)$. The time $T_{unif}(\theta)$ is equal to the minimal time of uniform {\em controllability}, and all lower bounds on the time $T_{unif}(\theta)$ obtained in Theorems~\ref{thmtmpsminomega} and~\ref{thmtmpsmin} and their corollaries apply.

\subsection{The vanishing viscosity limit for gradient flows. Conjugation and reformulation}
\label{sectconjug}

We focus in this article on the very particular case~\eqref{e:transport-viscous-gradient} where $X$ is a gradient vector field (with respect to the same metric $g$ defining the viscous perturbation $\eps \Lap$) of a weight function $\f : \M\to \R$, that is $X = \nablag \f \cdot \nablag$. In this case, it is classical (see e.g.~\cite{Witten:82,HS:85}) that the operator $-\eps \Lap - \nablag \f \cdot \nablag-\q$ can be conjugated to a ``semiclassical selfadjoint operator''. Here, $X\cdot Y$ is the inner product of the two vector fields $X$ and $Y$ given by the metric $g$ in $T\M$.

The first basic computation is the following:
$$
e^{-\frac{\f}{2\eps}} \Lap e^{\frac{\f}{2\eps}} = \Lap + \frac{1}{\eps}\nablag \f \cdot \nablag + \frac{|\nablag \f|_g^2}{4\eps^2} + \frac{\Lap \f}{2\eps} .
$$
We denote by
 \begin{equation}
 \label{e:def-Peps}
 \frac{1}{\eps^2}P_\eps := - \Lap + \frac{|\nablag \f|_g^2}{4\eps^2} + \frac{\Lap \f}{2\eps}-\frac{ \q}{\eps} , \quad \text{ that is } \quad P_\eps := - \eps^2 \Lap + \frac{|\nablag \f|_g^2}{4} + \eps  \qf  ,
 \end{equation}
 where $\qf= \frac{\Lap \f}{2}-\q$. The above computation implies that 
\baln
\label{e:trans-eps}
e^{-\frac{\f}{2\eps}} \big( \frac{1}{\eps^2}P_\eps  \big)e^{\frac{\f}{2\eps}} = - \Lap - \frac{1}{\eps}\nablag \f \cdot \nablag-\frac{\q}{\eps} .
\naln
The interest of this conjugation is that the operator $P_\eps$ is selfadjoint in $L^2(\M,d\Vol_g)$ endowed with domain $D(P_\eps) = H^2(\M)\cap H^1_0(\M)$. Henceforth, the operator $- \Lap - \frac{1}{\eps}\nablag \f \cdot \nablag-\frac{\q}{\eps} $ is also selfadjoint in $L^2(\M, e^{\frac{\f}{2\eps}} d\Vol_g)$.
Let us now reformulate the uniform observability problem~\eqref{e:transport-viscous-obs} in terms of the heat equation involving the operator $P_\eps$ defined in~\eqref{e:def-Peps}.

\medskip
Note that the constant coefficient one dimensional problem introduced in~\cite{CG:05} enters the ``gradient flow'' setting with $\M =(0,L)\subset \R$, $g=1$, $\Lap=\d_x^2$, $\q=0$, $\f = M x$ for $M\in \R$, and thus $\nablag \f \cdot \nablag = M \d_x$. In that context, this form together with its formulation~\eqref{e:def-Peps} have already been used in~\cite{CG:05,Gla:10,Lissy:14,Lissy:15}.
\begin{lemma}
\label{lemequiveunif}
Given $T_0, C_0, \eps >0$ and a function $u$, the following statements are equivalent. 
\begin{enumerate}
\item The function $u$ solves 
  \begin{equation}
  \label{e:heat-transp-eps}
 \left\{
 \begin{aligned}
&(\d_t - \nablag \f\cdot \nablag -\q - \eps \Delta_g ) u = 0 , & \text{ in } (0,T_0)\times\Int(\M) ,  \\
& u = 0 & \text{ on } (0,T_0)\times \d \M ,
 \end{aligned}
 \right.
 \end{equation}
\begin{align}
\label{e:obs-heat-transp-eps}
\text{resp. }\quad \nor{u(T_0)}{L^2(\M)}^2 \leq C_0^2\int_0^{T_0}\nor{u}{L^2(\omega)}^2 dt , \\
\label{e:obs-heat-transp-epsboundary}
\text{resp. }\quad \nor{u(T_0)}{H^1_0(\M)}^2 \leq C_0^2\int_0^{T_0}\nor{\theta \eps \partial_{\nu}u|_{\d\M}}{H^{1/2}(\d\M)}^2 dt .
\end{align}
\item The function $v(t,x)= e^{\f(x)/2\eps}u(t,x)$ solves 
  \begin{equation}
  \label{e:heat-transp-eps-Witt}
 \left\{
 \begin{aligned}
&\eps \partial_t v + P_\eps v =0  , & \text{ in } (0,T_0)\times\Int(\M) ,  \\
& v = 0 & \text{ on } (0,T_0)\times \d \M ,
 \end{aligned}
 \right.
 \end{equation}
\begin{align}
\label{e:obs-heat-transp-eps-Witt}
\text{resp. }\quad \nor{e^{-\frac{\f}{2\eps}}v(T_0)}{L^2(\M)}^2 \leq C_0^2\int_0^{T_0}\nor{e^{-\frac{\f}{2\eps}}v }{L^2(\omega)}^2 dt , \\
\label{e:obs-heat-transp-eps-Wittboundary}
\text{resp. }\quad \nor{e^{-\frac{\f}{2\eps}}v(T_0)}{H^1_0(\M)}^2 \leq C_0^2\int_0^{T_0}\nor{\theta e^{-\frac{\f}{2\eps}} \eps \partial_{\nu}v|_{\d\M} }{H^{1/2}(\d\M)}^2 dt .
\end{align}
\item The function $w(t,x) = v(t/\eps, x)= e^{\f(x)/2\eps}u(t/\eps,x)$ solves
  \begin{equation}
  \label{e:heat-transp-eps-Witt-eps}
 \left\{
 \begin{aligned}
& \partial_t w + \frac{1}{\eps^2}P_\eps w =0  , & \text{ in } (0,T_0)\times\Int(\M) ,  \\
& w = 0 & \text{ on } (0,T_0)\times \d \M ,
 \end{aligned}
 \right.
 \end{equation}
\begin{align}
\text{resp. }\quad \nor{e^{-\frac{\f}{2\eps}}w(\eps T_0)}{L^2(\M)}^2 \leq C_0^2\int_0^{\eps T_0}\nor{e^{-\frac{\f}{2\eps}}w}{L^2(\omega)}^2 dt ,\nonumber \\
\label{e:obs-heat-transp-eps-Witt-epsboundary}
\text{resp. } \quad \nor{e^{-\frac{\f}{2\eps}}w(\eps T_0)}{L^2(\M)}^2 \leq C_0^2\int_0^{\eps T_0}\nor{\theta e^{-\frac{\f}{2\eps}} \eps \partial_{\nu}w}{H^{1/2}(\d\M)}^2 dt .
\end{align}
\end{enumerate}
\end{lemma}
\bnp
Start e.g. with $u$ satisfying~\eqref{e:heat-transp-eps}.
Using~\eqref{e:trans-eps} and the definition of $P_\eps$, Equation~\eqref{e:heat-transp-eps} rewrites equivalently as 
\bal
\partial_t u + e^{-\frac{\f}{2\eps}} \frac{1}{\eps}P_\eps e^{\frac{\f}{2\eps}}  u=0 , \quad t \in [0,T_0] .
\nal
The function $v= e^{\frac{\f}{2\eps}}  u$ then satisfies~\eqref{e:heat-transp-eps-Witt} (and conversely). Setting $w(t,x) = v(t/\eps , x)$ then satisfies~\eqref{e:heat-transp-eps-Witt-eps}, and conversely.

The proof that \eqref{e:obs-heat-transp-epsboundary}$\Leftrightarrow$\eqref{e:obs-heat-transp-eps-Wittboundary}$\Leftrightarrow$\eqref{e:obs-heat-transp-eps-Witt-epsboundary} uses additionally that, on account to the Dirichlet boundary condition, we have $\d_\nu(e^{\frac{\f}{2\eps}} u) |_{\d\M} =e^{\frac{\f}{2\eps}} \d_\nu u |_{\d\M}$.
\enp

\subsection{(GCC) and controllability of the limit equation $\eps=0$}
\label{s:GCC}
In this section, we characterize the observability inequality~\eqref{e:transport-obs-ineq} for solutions of~\eqref{e:transport-libre} in terms of the Geometric Control Condition (GCC). In this section, $\M$ is always assumed without boundary.

We denote by $(\phi_t)_{t\in \R}$ the flow of the vector field $X$, namely
\baln
\label{e:def-flot}
\dot{\phi}_t(x) = X(\phi_t(x)) ,\quad  \phi_0(x) = x \in \M.
\naln
This flow is globally defined on account to the compactness of $\M$.
We consider the following definition of the geometric control condition in the manifold $\M$ for the vector field $X$ and the set $\omega$, which we denote (GCC).

\begin{definition}
\label{d:def-GCC}
Let $\M$ be a compact manifold without boundary, $X$ a Lipschitz vector field on $\M$, $\omega \subset \M$, $\chi \in C^0(\M)$, $I\subset \R$ and interval and $T>0$. We say that 
\begin{itemize}
\item $(\M,X, \omega,I)$ satisfies (GCC) if for all $x\in \M$, there is $t \in I$ such that $\phi_{-t}(x)\in \omega$;
\item $(\M,X, \omega,T)$ satisfies (GCC) if $(\M,X, \omega,(0,T))$ satisfies (GCC);
\item $(\M,X, \omega)$ satisfies (GCC) if there is $T>0$ such that $(\omega,T)$ satisfies (GCC);
\item $(\M,X, \chi,T)$ satisfies (GCC) if $(\M,X, \{\chi\neq 0\},T)$ satisfies (GCC);
\item $(\M,X, \chi)$ satisfies (GCC) if $(\M,X, \{\chi\neq 0\})$  satisfies (GCC).
\end{itemize}
\end{definition}
In this section, the manifold $\M$ is fixed. To lighten notation, we omit the dependence on $\M$ in $(\M,X, \omega,I)$ and we simply write $(X,\omega,I)$ instead of $(\M,X, \omega,I)$ (with a similar notation for the other definitions). 

Note in particular that an open set $\omega$ satisfying (GCC) must contain all singular points of the vector field $X$ (i.e. all points $x\in \M$ such that $X(x)=0$).
We now provide with different reformulations of this property. 
\begin{lemma}
\label{l:equiv-GCC}
Let $\M$ be a compact manifold and $X$ a Lipschitz vector field on $\M$. 
Given $\omega \subset \M$ and $T>0$, the following properties are equivalent:
\begin{enumerate}
\item \label{i:om-gcc} $(X,\omega,T)$ satisfies (GCC);
\item \label{i:+om} $\bigcup_{t\in (0,T)}\phi_t(\omega) \supset \M$;
\item \label{i:-om} $\bigcup_{t\in (0,T)}\phi_{-t}(\omega) \supset \M$;
\item \label{i:-GCC-X} $(-X,\omega,T)$ satisfies (GCC).
\end{enumerate}
\end{lemma}
\bnp
The definition of $(\X,\omega,T)$ satisfying (GCC) is equivalent to: for all $x\in \M$, there is $t \in (0,T)$ such that $x \in \phi_t(\omega)$.
Equivalence between the Items~\ref{i:om-gcc} and~\ref{i:+om} follows. Item~\ref{i:-om} is equivalent to $\bigcup_{t\in (0,T)}\phi_{T-t}(\omega) \supset \phi_{T} (\M) = \M$ after having applied $\phi_T$, which itself is equivalent to Item~\ref{i:+om}. Equivalence between Item~\ref{i:-GCC-X} and Item~\ref{i:-om} finally follows from the fact that the flow of $-X$ is $(\phi_{-t})_{t\in \R}$.
\enp

\begin{proposition}
Assume $\M$ is a compact manifold, $X$ is a Lipschitz vector field on $\M$, $ds$ is a positive density on $\M$, and $q \in L^\infty(\M)$.
Given $\omega \subset \M$, $\chi \in C^0(\M)$ and $T>0$, the following statements hold true:
\begin{enumerate}
\item \label{i:GCC-obs} If $(X,\omega ,T)$ satisfies (GCC), then the observability inequality~\eqref{e:transport-obs-ineq} for solutions of~\eqref{e:transport-libre} is true;
\item \label{i:obs-GCC} The observability inequality~\eqref{e:transport-obs-ineq} for solutions of~\eqref{e:transport-libre} implies that $(X,\ovl{\omega},[0,T])$ satisfies (GCC);
\item  \label{i:obs-GCC-cont} The observability inequality 
\begin{multline}
\label{e:obs-fct-cont}
C_0^2 \int_0^T\int_\M |\chi(x) u(t,x)|^2 ds(x)dt \geq \|u (T)\|_{L^2(\M)}^2 , \\
 \text{ for all } u_0 \in L^2(\M) \text{ and } u \text{ solution of~\eqref{e:transport-libre}}.
\end{multline}
holds true if and only if $(X,\chi,T)$ satisfies (GCC);
\item \label{i:equiv-0-T} In all the above observability statements, $\|u (T)\|_{L^2(\M)}^2$ can be equivalently replaced by $\|u (0)\|_{L^2(\M)}^2$.
\end{enumerate}
\end{proposition}
The proof below is inspired by that in~\cite{DL:09,LL:16} for the wave equation. It is constructive and would also yield a characterisation of the HUM control operator (see e.g.~\cite{HKL:15} or \cite[Section~1.2]{Lea:18} for more on controllability/stabilization properties for transport equations).
\bnp
First notice that for $u_0 \in L^2(\M)$, the unique solution to~\eqref{e:transport-libre} is explicitly given by
$$
u(t,x) = e^{\int_0^tq \circ\phi_{t-\tau}(x)d\tau}u_0\circ \phi_t(x)  \quad \in C^0(\R ; L^2(\M)) .
$$
A first direct consequence is the existence of a constant $C_{T,q}>1$ such that 
$$
C_{T,q}^{-1} \|u (0)\|_{L^2(\M)}^2\leq \|u (T)\|_{L^2(\M)}^2 \leq C_{T,q} \|u (0)\|_{L^2(\M)}^2 , \quad \text{ for all solutions to~\eqref{e:transport-libre}},
$$
which proves Item~\ref{i:equiv-0-T}. Next, we write the observation term in~\eqref{e:transport-obs-ineq} (the same holds for~\eqref{e:obs-fct-cont} if we replace $\mathds{1}_\omega$ by $\chi$) as 
\begin{align*}
\int_0^T\int_\omega |u(t,x)|^2 ds(x)dt  & = \int_0^T\int_\M \mathds{1}_\omega(x) |u(t,x)|^2 ds(x)dt \\
&  =  \int_0^T\int_\M \mathds{1}_\omega(x)  u_0^2\circ \phi_t(x) e^{2\int_0^tq \circ\phi_{t-\tau}(x)d\tau} ds(x)dt .
\end{align*}
Using the change of variable $y=\phi_t(x)$ (see e.g.~\cite[Proposition~16.42 p432]{Lee:book}), we obtain 
\begin{align*}
\int_0^T\int_\omega |u(t,x)|^2 ds(x)dt  
 =  \int_0^T\int_\M \mathds{1}_\omega(\phi_{-t}(y))  u_0^2(y)e^{2\int_0^tq \circ\phi_{-\tau}(y)d\tau} (\phi_{-t}^*ds)(y)dt
\end{align*}
(note that $\div_{ds}(X)$ is defined by $\frac{d}{dt} (\phi_{t}^*ds)|_{t=0} = \div_{ds}(X) ds$, so that this expression simplifies slightly in case $\div_{ds}(X)=0$).
Using that the density is positive on the compact $\M$, we get the existence of $C_T>1$ such that
$$
C_T^{-1} ds(y) \leq e^{2\int_0^tq \circ\phi_{-\tau}(y)d\tau} (\phi_{-t}^*ds)(y) \leq C_T ds(y) \quad \text{ uniformly for } (t,y) \in [0,T]\times \M .
$$
As a consequence, we obtain
\begin{align*}
C_T^{-1}  \int_0^T\int_\M \mathds{1}_\omega(\phi_{-t}(y))  u_0^2(y) ds(y)dt 
& \leq \int_0^T\int_\M \mathds{1}_\omega(x) |u(t,x)|^2 ds(x)dt \\
&   \leq C_T \int_0^T\int_\M \mathds{1}_\omega(\phi_{-t}(y))   u_0^2(y) ds(y)dt .
\end{align*}
Hence setting 
\begin{align*}
g_{\omega,T}(y) &:= \int_0^T  \mathds{1}_\omega(\phi_{-t}(y)) dt \in L^\infty(\M) \\
 \text{ resp. } \quad g_{\chi,T}(y) &:= \int_0^T  \chi^2 (\phi_{-t}(y)) dt \in C^0(\M) ,
\end{align*}
we deduce
\begin{align}
\label{e:equiv-g-t-obs}
C_T^{-1}  \int_\M g_{\omega,T}(y) u_0^2(y) ds(y)
& \leq \int_0^T\int_{\M} \mathds{1}_\omega(x) |u(t,x)|^2 ds(x)dt \nonumber \\
&  \leq C_T \int_\M g_{\omega,T}(y)  u_0^2(y) ds(y) . 
\end{align}
Recalling that $\omega$ is an open set and $\M$ compact, together with Definition~\ref{d:def-GCC}, we deduce that if $(X,\omega,T)$ satisfies (GCC), then we have the existence of $c>0$ such that $g_{\omega,T}(y)\geq c$ for a.e. $y\in\M$. The lower bound in~\eqref{e:equiv-g-t-obs} then implies the observability inequality~\eqref{e:transport-obs-ineq}, and Item~\ref{i:GCC-obs} follows. 

Concerning Item~\ref{i:obs-GCC}, if $(X,\ovl{\omega},[0,T])$ does not satisfy (GCC), then there is a point $x_0\in \M$ such that $\phi_{-t}(x_0) \cap \ovl{\omega} = \emptyset$ for all $t\in [0,T]$. The set $\ovl{\omega}\times [0,T]$ being compact, there is a neighborhood $U$ of $x_0$ such that $\phi_{-t}(U) \cap \ovl{\omega} = \emptyset$ for all $t\in [0,T]$. Setting $u_0 = \mathds{1}_U$, we have on the one hand that $\nor{u_0}{L^2(\M)} >0$. On the other hand, we have $\mathds{1}_{\omega}\circ \phi_{-t}(y)= 0$ for all $y\in U$ and $t\in (0,T)$. This implies that $g_{\omega,T} u_0 = g_{\omega,T}(y) \mathds{1}_U = 0$ a.e. and, according to the upper bound in~\eqref{e:equiv-g-t-obs}, that $\int_0^T\int_\omega |u(t,x)|^2 ds(x)dt  =0$. This contradicts~\eqref{e:equiv-g-t-obs} and concludes the proof of Item~\ref{i:obs-GCC}.

 Finally, the proof of Item~\ref{i:obs-GCC-cont} is split in two parts. That $(X,\chi,T)$ satisfies (GCC) implies the observability inequality~\eqref{e:obs-fct-cont} follows as in the proof of Item~\ref{i:GCC-obs}. Now assume that $(X,\chi,T)$ does not satisfy (GCC). Then there is a point $x_0\in \M$ such that $\phi_{-t}(x_0) \cap \{\chi \neq 0\} = \emptyset$ for all $t\in (0,T)$. Hence, we have $g_{\chi,T}(x_0)= \int_0^T  \chi^2 (\phi_{-t}(x_0)) dt =0$.
We now choose a sequence of continuous real-valued initial data $(u_0^n)_{n\in \N}$ such that $\nor{u_0^n}{L^2(\M)}=1$ and $(u_0^n)^2(x)ds(x) \rightharpoonup \delta_{x_0}$ in the sense of measures on $\M$. The fact that $g_{\chi,T}$ is continuous on $\M$ together with the upper bound in~\eqref{e:equiv-g-t-obs}
implies that, denoting by $u^n$ the solution of~\eqref{e:transport-libre} with initial datum $u_0^n$, we have
\begin{align*}
\int_0^T\int_{\M} |\chi(x)u^n(t,x)|^2 ds(x)dt  \leq&  C_T \int_\M g_{\chi,T}(y)  (u_0^n)^2(y) ds(y)\\
 &  \to C_T \langle \delta_{x_0} , g_{\chi,T} \rangle_{\Meas, C^0} = 0 ,
\end{align*}
which contradicts the observability inequality~\eqref{e:obs-fct-cont}, and concludes the proof of Item~\ref{i:obs-GCC-cont}.
\enp

\section{General lower bounds without geometric assumption}
\label{s:general-bounds}
In this section, we consider a general manifold (with or without boundary) $\M$, and prove the lower bound for the minimal time of uniform controllability provided in Theorem~\ref{thmtmpsminomega}. We also give a proof of Theorem~\ref{thm:counterexamplesomega} as a corollary.
To do this, we use the semiclassical reformulation~\eqref{e:heat-transp-eps-Witt}-\eqref{e:obs-heat-transp-eps-Witt} of the problem in Lemma~\ref{lemequiveunif}, as well as exponential decay properties of eigenfunctions of the operator $P_\eps$. We rely on the Helffer-Sj\"ostrand theory as developed in~\cite{HS:84,Helffer:booksemiclassic,DS:book}. All results presented in this section apply as well for the semiclassical heat equation.

The result of Theorem~\ref{thmtmpsminomega} is stated for a potential minimum. However, we shall prove a seemingly more general result, at any energy level in $V(\M)$. We shall then explain why this latter result is not more general, and how it can be improved in dimension one.
We recall the definition of $V$ in terms of $\f$ in~\eqref{e:def:V}, and define the classically allowed region at any energy level $E$: 
$$
 K_E = \{x \in \M , V(x) \leq E \} .
$$
We then define the Agmon distance (see e.g.~\cite[Chapter~3]{Helffer:booksemiclassic}) to the set $K_E$ at the energy level $E$: 
\baln
\label{defAgmon}
d_{A,E}(x ,y) &=  \inf \left\{ \int_0^1\sqrt{\left(V(\gamma(t)) - E \right)_+}  |\dot{\gamma}(t)|_g dt, \ \gamma \in \mathbf{U}_1(x,y) \right\} , \nonumber\\
 \mathbf{U}_1(x,y) & = \left\{ \gamma \in W^{1,\infty} ([0,1] ; \M ), \gamma(0)=x , \gamma(1)= y \right\} ,\nonumber\\
d_{A,E}(x)& = \inf_{y \in K_E}d_{A,E}(x ,y) .
\naln
That is to say, $d_{A,E}(x)$ is the distance of $x$ to the set $K_E$ for the (pseudo-)metric $(V-E)_+ g$. Here again $\left(V(x) - E \right)_+ = \max \left(V(x) - E , 0\right)$. We will use, as in~\eqref{defAgmon-0} the notation $d_A=d_{A,E_0}$ where $E_0=\min_{\M}V$ for the Agmon distance at the bottom energy.
Note that $d_{A,E}$ vanishes identically on $K_E$ (and only on this set).
Finally, an important function in the estimates below is given by
\baln
\label{e:def-WE}
W_E (x)=d_{A,E}(x)+\frac{\f(x)}{2}.
\naln
We shall prove in this section the following result.
\begin{theorem}
\label{thm-general}
Assume the observability estimate~\eqref{e:obs-heat-transp-eps} (resp. the boundary observability inequality~\eqref{e:obs-heat-transp-epsboundary}) for all solutions to~\eqref{e:heat-transp-eps} with constant $C_0 = C_0(T_0, \eps)$. 
Then, for all $E \in V(\M) = [\min_\M V, \max_\M V]$ and all $\delta >0$, there is $\eps_0>0$ such that we have for all $\eps\in (0,\eps_0)$
$$
C_0(T_0, \eps) \geq  \exp \frac{1}{\eps} \left( \min_{\ovl{\omega}}W_E - \max_{K_E} W_E  - \delta - E T_0 \right)  , 
$$
in the internal observation case, and 
$$
C_0(T_0, \eps) \geq  \exp \frac{1}{\eps} \left( \min_{\ovl{\Gamma}}W_E - \max_{K_E} W_E  - \delta - E T_0 \right)  , 
$$
in the boundary observation case.

In particular, we have for all $E \in V(\M)$, for each respective case,
\bal
T_{unif}(\omega) &\geq \frac{1}{E}\left(\min_{\ovl{\omega}}W_E - \max_{K_E} W_E  \right)  ,\\
T_{unif}(\Gamma) &\geq \frac{1}{E}\left(\min_{\ovl{\Gamma}}W_E - \max_{K_E} W_E  \right) .
\nal
\end{theorem}
Theorem~\ref{thmtmpsminomega} is then the particular case $E=E_0=\min_\M V$ in Theorem~\ref{thm-general}.
Unfortunately, the function $E \mapsto \frac{1}{E} \left(\min_{\ovl{\omega}}W_E - \max_{K_E} W_E \right)$ is a decreasing function of $E$. Indeed, 
\begin{itemize}
\item the sets $K_E$ are increasing in $E$, hence the function $E \mapsto \max_{K_E} W_E= \max_{K_E} \frac{\f}{2}$ increases; 
\item $E \mapsto d_{A,E}(x)$ is decreasing in $E$, hence the function $E \mapsto\min_{\ovl{\omega}}W_E$ decreases. 
\end{itemize}
Therefore, the estimate of Theorem~\ref{thm-general} simply reduces to that for $E=E_0$, that is Theorem~\ref{thmtmpsminomega} in the introduction.
This comes from the fact that the estimate involving the term $\max_{K_E} W_E$ is very rough (see Section \ref{sectrough} below for a more precise discussion). This can be improved in the one dimensional one well case, see~\cite{LL:20-1D}.
Moreover, it is interesting to notice that the proof of Theorem~\ref{thm-general} is not more involved than the direct proof of Theorem~\ref{thmtmpsminomega}, and we shall re-use part of it in case of revolution surfaces in Section~\ref{s:revolution}.

\subsection{Eigenfunctions of semiclassical Schr\"odinger operators}
\label{s:sub-eig-schro}
In this section, we collect classical results concerning eigenfunctions of semiclassical Schr\"odinger operators, and some of their decay properties.
Recall that $P_\eps = - \eps^2 \Lap + \frac{|\nablag \f|_g^2}{4} + \eps  \qf =  - \eps^2 \Lap + V+ \eps  \qf $ is defined in~\eqref{e:def-Peps}.
We first need to prove existence of eigenfunctions near any energy level. 
\begin{lemma}[Existence of eigenfunctions]
\label{l:exist-eig}
Assume $V \in W^{1,\infty}(\M)$ and $\qf \in L^\infty(\M)$  are both real valued. 
For all $E \in V(\M) = [\min_\M V, \max_\M V]$ and all $\eps \in (0,1]$, there is $E_\eps = E + \grando{\eps^{2/3}}$ and $\psi_\eps \in H^2(\M) \cap H^1_0(\M)$ such that $P_\eps \psi_\eps = E_\eps \psi_\eps$.
\end{lemma}
Note that the $\grando{\eps^{2/3}}$ precision is relatively poor, and can be improved in different situations (e.g. if there is a critical point of $V$ at energy $E$). These refinements are however not needed here.

\bnp
The proof consists in constructing a (very rough) quasimode. 
Assume first that $E$ is reached by an interior point, i.e. there is $x_0 \in \Int(\M)$ such that $V(x_0) = E$. 
 We then only work in a local chart near $x_0$, centered at $x_0$ (hence we work in $\R^n$ in a neighborhood of $0$).
 
 We take a cutoff function $\chi \in C^\infty_c(\R)$ such that $\chi = 1$ in a neighborhood of $0$. We set 
$u^\eps (x) = \eps^{-\frac{n}{3}} \chi \left( \eps^{-\frac23} |x| \right)$, so that $u^\eps$ is smooth and moreover supported in the chart for all $\eps<\eps_0$ with $\eps_0$ sufficiently small, and thus $u^\eps \in H^2(\M) \cap H^1_0(\M)$. Notice also that 
\bal
\nor{u^\eps}{L^2(\M)}^2 = \int |u^\eps(x)|^2 \sqrt{|g|}(x) dx = \int \chi^2(|y|) \sqrt{|g|}(\eps^{2/3} y) dy = c_0+ \grando{\eps^{2/3}} ,
\nal
with $c_0>0$.
We now estimate $(P_\eps - E)u^\eps$.
For this, we first have $\nor{\eps \qf u^\eps }{L^2} \leq C \eps$. Second, we always have the rough estimate $V(x)-E = V(x)-V(0) = \grando{|x|}$ so that we have 
\begin{align}
\label{e:quasim-est-1}
\nor{ (V-E) u^\eps }{L^2}^2  \leq C \int |x|^2 \chi^2 \left( \eps^{-\frac23} |x| \right)  \frac{dx}{ \eps^{\frac{2n}{3}}} \leq C \eps^{4/3}.
\end{align}
Third, we have 
\begin{align}
\label{e:quasim-est-2}
\nor{ \eps^2 \Delta_g u^\eps }{L^2}^2& = \nor{\frac{\eps^2}{\sqrt{|g|}}\sum_{i,j} \d_i \left(g^{ij} \sqrt{|g|} \d_j u^\eps \right) }{L^2}^2 \nonumber \\
& \leq C\eps^4 \sum_{i,j} \int \Big[\d_i \left(g^{ij} \sqrt{|g|}  \frac{x_j}{|x|} \right)  \eps^{-\frac23}\chi' \left( \eps^{-\frac23}|x| \right) \nonumber \\
& \quad + \left(g^{ij} \sqrt{|g|} \frac{x_j}{|x|} \right) \frac{x_i}{|x|} \eps^{-4/3} \chi'' \left( \eps^{-\frac23} |x| \right) \Big]^{2}  \frac{dx}{ \eps^{\frac{2n}{3}}} \nonumber \\
& \leq C \eps^{\frac43}.
\end{align}
Combining the above four estimates yields the existence of $D,\eps_0>0$ such that for all $\eps<\eps_0$, we have, 
$$
\nor{(P_\eps - E) u^\eps}{L^2}  \leq C \eps^{\frac23} \leq D\eps^{\frac23} \nor{u^\eps}{L^2}  .
$$
Hence, if $E \notin \Sp(P_\eps)$, this implies $\nor{(P_\eps - E)^{-1}}{L^2\to L^2}\geq (D\eps^{\frac23})^{-1}$.
Finally, the operator $P_\eps$ being selfadjoint, we have, for $z \in \C \setminus \Sp(P_\eps)$, $\|(P_\eps-z)^{-1}\|_{L^2\to L^2} = \frac{1}{d(z, \Sp(P_\eps))}$, so that, if $E \notin \Sp(P_\eps)$, 
$$
\frac{1}{d(E, \Sp(P_\eps))} \geq (D\eps^{\frac23})^{-1}. 
$$
In any case, this implies $d(E, \Sp(P_\eps))\leq D\eps^{\frac23}$, and using that the spectrum of $P_\eps$ is purely pointwise, this proves the sought result.

Assume now that $E$ is not reached by an interior point, i.e. $E \notin V(\Int(\M))$. This means in particular that 
 there is $x_0 \in \d \M$ such that $V(x_0) = E$. 
 Then, we again work in a local chart near $x_0$, centered at $x_0$. In this chart, $\M$ is given by $\R^{n-1}\times\R^-$ and $x_0$ by $0$. We denote $(x', x_n) \in \R^{n-1}\times\R^-$ local coordinates.
We then take $\chi$ as above and further define $\check{\chi} \in C^\infty_c(\R)$, non-identically vanishing, such that $\supp(\check{\chi}) \subset (-1,0)$. We define $u^\eps (x) = \eps^{-\frac{n}{3}}\check{\chi}(\eps^{-\frac23} x_n) \chi ( \eps^{-\frac23} |x'| )$. One can check that all above properties of $u^\eps$ are still satisfied, and in particular~\eqref{e:quasim-est-1}-\eqref{e:quasim-est-2}.
In addition, we have $\supp u^\eps \subset \R^{n-1}\times\R^-_*$, and thus $u^\eps \in H^2(\M) \cap H^1_0(\M)$. The remainder of the proof then follows the same as in the first case.
\enp

\begin{remark}
Note that near a noncritical value of $V$, or near the boundary of $\M$, the appropriate local model is $-\eps^2\partial_x^2\pm x$. Considering concentrating quasimodes of the form $\chi\left(\frac{x}{\eps^{\alpha}}\right)$ leads to 
\bal
\eps^2\partial_x^2\left( \chi\left(\frac{x}{\eps^{\alpha}}\right)\right)=\frac{\eps^2}{\eps^{2\alpha}}\chi''\left(\frac{x}{\eps^{\alpha}}\right), 
\quad  \text{ and } \quad 
x\chi\left(\frac{x}{\eps^{\alpha}}\right)=\eps^{\alpha} \left(\frac{x}{\eps^{\alpha}}\right)\chi\left(\frac{x}{\eps^{\alpha}}\right) .
\nal
Henceforth, the right scaling is given by $2-2\alpha =\alpha$, that is $\alpha=2/3$.
The quasimode we construct are then $\grando{\eps^{2/3}}$. If one wants to obtain a better remainder, one could replace $\chi$ by an Airy function, as one should replace $\chi$ by a Hermite function in the case of bottom of potential (in which case the precision of the quasimode is improved). Also, the remainder $\e^{2/3}$ is actually the worst possible case. 
\end{remark}

The next result states the decay estimates for eigenfunctions in the classically forbidden region, and is a consequence of so-called Agmon estimates (see~\cite{HS:84,Helffer:booksemiclassic,DS:book}). Here, it is a particular case of~\cite[Propositions~3.3.1 and~3.3.4]{Helffer:booksemiclassic}. Note that with respect to~\cite[Proposition~3.3.1]{Helffer:booksemiclassic}, our operator $P_\eps$ contains an additional term, namely multiplication by $\eps \qf$. However, this contribution is of lower order and can be absorbed in the proof of the Agmon estimates.

\begin{theorem}[Decay of eigenfunctions in the classically forbidden region]
\label{t:decay-agmon}
Assume $V \in W^{1,\infty}(\M)$ and $\qf \in L^\infty(\M)$. Let $$E \in V(\M) = [\min_\M V, \max_\M V]$$ and assume $\psi_\eps \in H^2(\M) \cap H^1_0(\M)$ and $E_\eps$ satisfy
\begin{align}
\label{e:requirements-E}
P_\eps \psi_\eps = E_\eps \psi_\eps , \quad \nor{\psi_\eps}{L^2(\M)} = 1, \quad \max(E_\eps -E , 0)  = o(1) \text{ as } \eps \to 0^+ .
\end{align}
Then for all $\delta>0$, there exist $C=C(\delta) ,\eps_0=\eps_0(\delta)>0$ such that, for all $\eps<\eps_0$, we have 
\begin{align}
\label{e:estimee-L2-psi}
\nor{  e^{\frac{1}{\eps} d_{A,E}} \psi_\eps}{L^2(\M)}\leq C e^{\frac{\delta}{\eps}} .
\end{align}
Assuming further that $V, \qf \in C^\infty(\M)$, we have $\psi_\eps \in C^\infty(\M)$ and for all $\delta>0$ and all smooth vector field $Y$ on $\M$, there exist $C=C(\delta) ,\eps_0=\eps_0(\delta)>0$ such that, for all $0<\eps<\eps_0$, we have 
\begin{align}
\label{e:estimee-ponctuelle-psi}
 |\psi_\eps(x)| +  |(Y\psi_\eps)(x)| \leq C e^{-\frac{1}{\eps} (d_{A,E}(x)-\delta)} , \quad \text{ for all } x \in \M.
\end{align}
\end{theorem}
\begin{remark}
\label{r:Cinfty-reg-bdry}
Note that the smoothness assumption $V, \qf \in C^\infty(\M)$ (as well as the smoothness assumptions in Theorem~\ref{thmtmpsminomega}) is essentially only used in~\cite{HS:84,Helffer:booksemiclassic,DS:book} to perform elliptic regularity estimates so that the pointwise estimate~\eqref{e:estimee-ponctuelle-psi} makes sense. A finer (much less demanding) regularity assumption can be formulated.
\end{remark}

As a direct corollary, we have that most of the norm of $\psi_\eps$ is near $K_E$, see~\cite[Corollary~3.3.2]{Helffer:booksemiclassic}.
\begin{corollary}[Most of the norm is in the classically allowed region]
\label{corpresx0omega}
Let $E \in V(\M) = [\min_\M V, \max_\M V]$ and assume $\psi_\eps , E_\eps$ satisfy~\eqref{e:requirements-E}.
For any open set $U$ containing $K_E$, there is $\delta, \eps_0 >0$ such that for all $0<\eps< \eps_0$, we have
\bal
\nor{\psi_\eps}{L^2(U)}^2 = 1 + \grando{e^{-\frac{\delta}{\eps}}} .
\nal
\end{corollary}

\subsection{Rough localization of eigenfunctions, and a proof of Theorem~\ref{thm-general}}
\label{sectrough}
From the decay estimates in the classically forbidden region (Theorem~\ref{t:decay-agmon}) and the rough localization of the $L^2$ mass of eigenfunctions (Corollary~\ref{corpresx0omega}), together with the existence of eigenfunctions at any energy level (Lemma~\ref{l:exist-eig}), we may now deduce a proof of Theorem~\ref{thm-general}.
Recall that $W_E$ is defined in~\eqref{e:def-WE}. We first prove the following proposition, from which Theorem~\ref{thm-general} will follow.

\begin{proposition}
\label{prop:eign-rough}
Let $E \in V(\M) = [\min_\M V, \max_\M V]$ and assume $\psi_\eps \in H^2(\M) \cap H^1_0(\M)$ and $E_\eps$ satisfy~\eqref{e:requirements-E}. 
Then for all $\delta>0$, there exists $\eps_0=\eps_0(\delta)>0$ such that, for any open set $\omega \subset \M$ and for all $\eps<\eps_0$, we have
\begin{align}
\label{e:rough-lower-bound-eign}
\nor{e^{-\frac{\f}{2\eps}} \psi_\eps}{L^2(\M)} \geq  e^{-\frac{1}{\eps}(\max_{K_E}W_E + \delta)} ,\\
\label{e:upper-bound-eign}
\nor{e^{-\frac{\f}{2\eps}} \psi_\eps}{L^2(\M)} \leq  e^{-\frac{1}{\eps}(\min_{\M} W_E  - \delta)} , \\
\label{e:upper-bound-eign-obs}
\nor{e^{-\frac{\f}{2\eps}} \psi_\eps }{L^2(\omega)} \leq  e^{-\frac{1}{\eps} \left( \min_{\ovl{\omega}} W_E -\delta \right)} .
\end{align}
Assuming also that $V, \qf \in C^\infty(\M)$ and $\Gamma\subset \partial\M$, we have
\begin{align}
\label{e:upper-bound-eign-obs-boundary}
\nor{e^{-\frac{\f}{2\eps}} \d_\nu \psi_\eps }{L^2(\Gamma)} \leq  e^{-\frac{1}{\eps} \left( \min_{\ovl{\Gamma}} W_E -\delta \right)} .
\end{align}
\end{proposition}
Note that Estimate~\eqref{e:rough-lower-bound-eign} is very rough, due to our lack of knowledge on the localization of $\psi_\eps$ in the classically allowed region $K_E$. In the one dimensional one well case, this bound can be refined, see the companion paper~\cite{LL:20-1D}. 

We first prove Proposition~\ref{prop:eign-rough}, and then deduce a proof of Theorem~\ref{thm-general}.
\bnp[Proof of Proposition~\ref{prop:eign-rough}]
First, setting $U_\delta = \{x\in \M, \f(x) < \max_{K_E} \f + \delta\} \supset K_E$, we have 
\begin{align*}
\nor{e^{-\frac{\f}{2\eps}} \psi_\eps}{L^2(\M)} 
& \geq \nor{e^{-\frac{\f}{2\eps}} \psi_\eps}{L^2(U_\delta)}  \geq e^{-\frac{1}{\eps}(\max_{K_E}\frac{\f}{2} + \frac{\delta}{2})} \nor{ \psi_\eps}{L^2(U_\delta)}  
\\
& \geq \frac12  e^{-\frac{1}{\eps}(\max_{K_E}\frac{\f}{2} + \frac{\delta}{2})},
\end{align*}
after having used Corollary~\ref{corpresx0omega} in the last inequality for $\eps < \eps_0(\delta)$. This proves~\eqref{e:rough-lower-bound-eign}, recalling that $d_{A,E} = 0$ on $K_E$.

Second, notice that~\eqref{e:upper-bound-eign} is a consequence of~\eqref{e:upper-bound-eign-obs}. Third, to prove~\eqref{e:upper-bound-eign-obs}, we use~\eqref{e:estimee-L2-psi} as follows 
\begin{align*}
 \nor{e^{-\frac{\f}{2\eps}} \psi_\eps }{L^2(\omega)}
& \leq  \nor{e^{-\frac{\f}{2\eps}} e^{-\frac{1}{\eps} d_{A,E}}  e^{\frac{1}{\eps} d_{A,E}} \psi_\eps }{L^2(\omega)}
 \leq e^{-\frac{1}{\eps} \left( \min_{\ovl{\omega}}(\frac{\f}{2} + d_{A,E}) \right)} \nor{e^{\frac{1}{\eps} d_{A,E}} \psi_\eps }{L^2(\omega)} \\
&  \leq e^{-\frac{1}{\eps} \left( \min_{\ovl{\omega}}(\frac{\f}{2} + d_{A,E}) \right)}C e^{\frac{\delta}{\eps}} .
\end{align*}
Finally, to prove~\eqref{e:upper-bound-eign-obs-boundary}, we proceed similarly using~\eqref{e:estimee-ponctuelle-psi} (instead of the sole~\eqref{e:estimee-L2-psi}) with $Y =\d_\nu$. We have
\begin{align*}
\nor{e^{-\frac{\f}{2\eps}} \d_\nu \psi_\eps }{L^2(\Gamma)} \leq \nor{e^{-\frac{\f}{2\eps}}  C e^{-\frac{1}{\eps} (d_{A,E}-\delta)} }{L^2(\Gamma)} \leq C e^{-\frac{1}{\eps} \left( \min_{\ovl{\Gamma}}(\frac{\f}{2} + d_{A,E}) -\delta\right)} ,
 \end{align*}
which concludes the proof of the proposition.
\enp

\bnp[Proof of Theorem~\ref{thm-general} from Proposition~\ref{prop:eign-rough} and Lemma~\ref{l:exist-eig}]
We use the reformulation in Lemma~\ref{lemequiveunif} and consider the observability estimate~\eqref{e:obs-heat-transp-eps-Witt} for solutions to the evolution equation~\eqref{e:heat-transp-eps-Witt}.

More precisely, we select $E \in V(\M)$, and we let $v_\eps$ be the solution to~\eqref{e:heat-transp-eps-Witt} associated to the initial condition $v_\eps(0) = \psi_\eps$, where $\psi_\eps$ is given by Lemma~\ref{l:exist-eig}. That is to say, $v_\eps(t ,x) = e^{-\frac{E_\eps}{\eps}t}\psi_\eps(x)$. We estimate both sides of~\eqref{e:obs-heat-transp-eps-Witt}.

Firstly, using~\eqref{e:rough-lower-bound-eign}, we have 
\begin{align*}
\nor{e^{-\frac{\f}{2\eps}}v_\eps(T_0)}{L^2(\M)}  = e^{-\frac{E_\eps}{\eps}T_0} \nor{e^{-\frac{\f}{2\eps}} \psi_\eps}{L^2(\M)} 
\geq \frac12 e^{-\frac{E_\eps}{\eps}T_0}e^{-\frac{1}{\eps}(\max_{K_E}W_E  + \frac{\delta}{2})},
\end{align*}
 for $\eps < \eps_0(\delta)$. Recalling that $E_\eps = E + \grando{\eps^{2/3}}$, this implies the existence of $\eps_0(\delta)>0$ such that for $\eps < \eps_0(\delta)$, 
\begin{align*}
\nor{e^{-\frac{\f}{2\eps}}v_\eps(T_0)}{L^2(\M)}  \geq \frac12 e^{-\frac{1}{\eps}(E T_0 + \max_{K_E}W_E  + \delta)},
\end{align*}

Secondly, concerning the internal case, using~\eqref{e:upper-bound-eign-obs}, we have 
\begin{align*}
\int_0^{T_0}\nor{e^{-\frac{\f}{2\eps}}v_\eps }{L^2(\omega)}^2 dt 
& = \int_0^{T_0}e^{-2 \frac{E_\eps}{\eps}t} \nor{e^{-\frac{\f}{2\eps}} \psi_\eps }{L^2(\omega)}^2 dt \\
 &= \frac{\eps}{2E_\eps} (1- e^{-2 \frac{E_\eps}{\eps}T_0})\nor{e^{-\frac{\f}{2\eps}} \psi_\eps }{L^2(\omega)}^2 
 \leq C e^{-\frac{2}{\eps} \left( \min_{\ovl{\omega}} W_E -\delta \right)} .
\end{align*}
As a consequence of these two estimates, applying~\eqref{e:obs-heat-transp-eps-Witt} implies 
\begin{align*}
 e^{-\frac{2}{\eps}(E T_0 + \max_{K_E}W_E  + \delta)} \leq C_0(T_0,\eps)^{2} C  e^{-\frac{2}{\eps} \left( \min_{\ovl{\omega}} W_E -\delta \right)} ,
\end{align*}
which concludes the proof.

In the case of observation from part of the boundary, we simply replace the above use of~\eqref{e:upper-bound-eign-obs} by that of~\eqref{e:upper-bound-eign-obs-boundary}.
\enp

\subsection{An explicit counter-example for a domain of $\R^n$}
The purpose of this section is to prove Theorem \ref{thm:counterexamplesomega}. Here $\M=\ovl{\Omega}$ where $\Omega\subset \R^n$ is an open set, endowed with the Euclidean metric. For $\delta>0$, we may assume, up to an appropriate translation of the domain $\Omega$, that there is $\eta>0$ such that 
\baln
\label{e:asspt-eta}
& \overline{B}(0,\eta)\subset \Omega .  \\
\label{e:cond-omega}
& \omega \subset  \left( \Omega \cap \{x_1>0,  x_2>0 ,  \cdots ,x_n >0\}  \right),\quad \text{ and } \quad  B(0,\eta) \cap \omega = \emptyset.
\naln
We let $\f_\lambda(x_1,\cdots,x_n)$ be defined as follows: 
\baln
\label{e:def-f-lambda}
f_\lambda(t):=\int_0^t\sqrt{\lambda^2s^2+1}ds \quad \text{ and } \quad \f_\lambda(x_1,\cdots,x_n) :=\sum_{i=1}^n f_\lambda(x_i) .
\naln
With this definition, the associated gradient vector field and potential are given respectively by
\bal
\nabla \f_\lambda(x_1,\cdots,x_n) =\sum_{i=1}^n f_\lambda'(x_i) e_i= \sum_{i=1}^n \sqrt{\lambda^2 x_i^2+1} e_i  ,\\
V_\lambda (x) =\frac{|\nabla \f_\lambda|^2}{4}=\frac{\lambda^2 |x|^2+n}{4} ,
\nal
where $(e_1,\cdots,e_n)$ denotes the canonical basis of $\R^n$.

The proof of Theorem \ref{thm:counterexamplesomega} now directly follows from the following Lemmata~\ref{lminfTOmega} and~\ref{lmsuptransOmega}, when taking $\lambda$ large enough. 

\begin{lemma}
\label{lminfTOmega}
In the above setting, recalling~\eqref{e:asspt-eta}, \eqref{e:cond-omega} and that $\f_\lambda$ is defined in~\eqref{e:def-f-lambda}, we have $T_{unif}(\omega)\geq \frac{\lambda\eta^2}{n}$ for all $\lambda>0$.
\end{lemma}

\begin{lemma}
\label{lmsuptransOmega}
In the above setting, recalling~\eqref{e:asspt-eta}, \eqref{e:cond-omega} and that $\f_\lambda$ is defined in~\eqref{e:def-f-lambda}, $(\ovl{\Omega}, \nabla \f_\lambda, \omega)$ and $(\ovl{\Omega}, \nabla \f_\lambda)$ both satisfy (FC) and we have for all $\lambda\geq0$ 
\baln
\label{e:estim-time-bdry}
T_{FC}(\ovl\Omega,\nabla \f_\lambda,\omega) & \leq T_{FC}(\ovl\Omega,\nabla \f_\lambda)
 \leq \min_{\mathsf{v}\in \R^n, |\mathsf{v}|=1, \mathsf{v}\cdot e_i \geq 0 \forall i}\left( \max_{x \in \ovl{\Omega}} x\cdot \mathsf{v} - \min_{x \in \ovl{\Omega}} x\cdot \mathsf{v}\right) \nonumber \\
& \leq \diam(\Omega) .
\naln
In particular, both are bounded by a constant uniformly in $\lambda\geq 0$.
\end{lemma}
We recall that (FC) and $T_{FC}$ are defined in this context in Definition~\ref{d:def-GCC_bdry}. 

Notice that the quantity $\max_{x \in \ovl{\Omega}} x\cdot \mathsf{v} - \min_{x \in \ovl{\Omega}} x\cdot \mathsf{v}$ represents the minimal Euclidean distance between two parallel hyperplanes (normal to $\mathsf{v}$) such that $\Omega$ is contained between the two hyperplanes.

\bnp[Proof of Lemma~\ref{lminfTOmega}]
The minimum of $V_\lambda$ is reached at $x_{\min}=0$ and $V_\lambda(x_{\min})=\frac{n}{4}$. The Agmon distance~\eqref{defAgmon} at the bottom energy $V_\lambda(x_{\min})$ can be explicitly computed for points $x \in \ovl{B}(0,\eta) \subset \Omega \setminus \omega$. Indeed, for $x \in \Omega$ we have $\sqrt{V_\lambda(x) - V_\lambda(x_{\min})} = \frac{\lambda}{2}|x|$ and thus for $x \in \ovl{B}(0,\eta)$, 
\bal
d_{A}(x) & = \frac{\lambda}{2} \inf \left\{ \int_0^1|\gamma(t)||\dot{\gamma}(t) |dt , \gamma \in W^{1,\infty} ([0,1] ; \M ), \gamma(0)=0 , \gamma(1)= x \right\} \\
& = \frac{\lambda |x|^2}{2} \int_0^1t dt=\frac{\lambda |x|^2}{4}  ,
\nal
where, by symmetry arguments, we have noticed that the straight line $\gamma(t):=xt$ reaches the infinimum for $x \in \ovl{B}(0,\eta)$.
Recalling that 
$V_\lambda(0) = \frac{n}{4}$ and $\f(0) =0$, application of~\eqref{e:estiminf-Tunifomega} in Theorem~\ref{thmtmpsminomega} implies
\baln
\label{e:app-coro}
\frac{n}{4}T_{unif}(\omega)= V_\lambda(0)T_{unif}(\omega) \geq \min_{\ovl{\omega}}\left( \frac{\f}{2} + d_{A} \right) .
\naln
By a connectedness argument, for any $x \in \ovl{\omega}$, we have from~\eqref{e:cond-omega} that  $d_A(x) \geq \min_{y\in \partial B(0,\eta)}  d_A(y)\geq \frac{\lambda \eta^2}{4}$. 
 Moreover, the condition~\eqref{e:cond-omega} together with the definition of $\f_\lambda$ in~\eqref{e:def-f-lambda} imply that $\f_\lambda(x)\geq 0$ for $x \in \omega$. We thus have
 $$
\min_{\ovl{\omega}}\left( \frac{\f_\lambda}{2} + d_{A} \right)  \geq \min_{\ovl{\omega}}\left( d_{A} \right) \geq \frac{\lambda \eta^2}{4} .
$$
Combined with~\eqref{e:app-coro}, this concludes the proof of the lemma.
\enp

\bnp[Proof of Lemma~\ref{lmsuptransOmega}]
Notice first that given $\mathsf{v}\in \R^n$ with $|\mathsf{v}|=1$ such that $e_i \cdot \mathsf{v} \geq 0$ for all $i \in \{1, \cdots , n\}$, we have for all $x\in \ovl{\Omega}$,
\baln
\label{nablaflambda}
\nabla \f_\lambda(x) \cdot \mathsf{v} = \sum_{i=1}^n \sqrt{\lambda^2 x_i^2+1} e_i \cdot \mathsf{v} \geq \sum_{i=1}^n e_i \cdot \mathsf{v} \geq \sqrt{\sum_{i=1}^n (e_i \cdot \mathsf{v})^2 }= |\mathsf{v}|= 1.
\naln
Second, following Definition~\ref{d:def-GCC_bdry}, we extend $\f_{\lambda}$ as a smooth compactly supported function $\tilde{\f}_{\lambda}$  in $\R^n$
Given $x\in \ovl{\Omega}$, we denote by $\gamma_x$ the maximal (global) solution to $\dot{\gamma}_x(t) = -\nabla\tilde{\f}_{\lambda}(\gamma_x(t))$ with $\gamma_x(0)=x$, defined in $\R^n$ for $t\in \R$. 

Given $\mathsf{v}\in \R^n$ with $|\mathsf{v}|=1$ such that $e_i \cdot \mathsf{v} \geq 0$, Estimate~\eqref{nablaflambda} thus implies that $\dot{\gamma}_x(t)\cdot \mathsf{v} =  -\nabla\tilde{\f}_{\lambda}(\gamma_x(t))\cdot \mathsf{v}\leq -1$. Integrating this between $0$ and $t$, we obtain $\gamma_x(t)\cdot \mathsf{v}-x\cdot \mathsf{v} \leq -t.$
Assuming that $T>  \max_{y \in \ovl{\Omega}} y\cdot \mathsf{v} - \min_{y \in \ovl{\Omega}} y\cdot \mathsf{v}$ thus implies 
$$
\gamma_x(T)\cdot \mathsf{v}-x\cdot \mathsf{v} \leq -T  < - \left(\max_{y \in \ovl{\Omega}} y\cdot \mathsf{v} - \min_{y \in \ovl{\Omega}} y\cdot \mathsf{v} \right) ,
$$
that is 
$$
x\cdot \mathsf{v}-\gamma_x(T)\cdot \mathsf{v}  >  \max_{y \in \ovl{\Omega}} y\cdot \mathsf{v} - \min_{y \in \ovl{\Omega}} y\cdot \mathsf{v}  .
$$
Since $x \in \ovl{\Omega}$, this implies $\gamma_x(T) \notin \ovl{\Omega}$. This holds true for any $x \in \ovl{\Omega}$.
 Recalling the definition of (FC) in Definition~\ref{d:def-GCC_bdry}, we have obtained that this condition is satisfied by both $(\ovl\Omega,\nabla \f_\lambda,T)$ and $(\ovl\Omega,\nabla \f_\lambda,\omega,T)$. Moreover, given the definition of $T_{FC}$ as an infimum, we have also obtained that $T_{FC}(\ovl\Omega,\nabla \f_\lambda,\omega) \leq T_{FC}(\ovl\Omega,\nabla \f_\lambda)\leq \max_{y \in \ovl{\Omega}} y\cdot \mathsf{v} - \min_{y \in \ovl{\Omega}} y\cdot \mathsf{v}$. Since this holds for all $\mathsf{v}$, we have proved~\eqref{e:estim-time-bdry}, which concludes the proof of the lemma.
 \enp

\section{Surfaces of revolution}
\label{s:revolution}
\subsection{General setting}
\label{s:revol}
In this section we introduce the geometric setting for the results presented in Section~\ref{s:revol-intro}. We are concerned with a revolution surface $\mathcal{S}\subset \R^3$ being either
\begin{enumerate}
\item \label{caseS}Case \ref{caseS}: diffeomorphic to a sphere $\S^2$ (in which case $\d \mathcal{S} = \emptyset$);
\item \label{caseD}Case \ref{caseD}:  diffeomorphic to a disk (in which case $\d \mathcal{S}$ is a circle embedded in $\R^3$).
\item \label{caseCY}Case \ref{caseCY}: diffeomorphic to a cylinder $[0,1] \times \S^1 \subset \R^3$ (in which case $\d \mathcal{S}$ consists in two disjoint circles embedded in $\R^3$, and belonging to two parallel hyperplanes);
\end{enumerate}
We follow~\cite{LL:18} and~\cite[Chapter~4B p95]{Besse} for the precise geometric description of such manifolds. 
At some places, we also consider the case of the torus $\T^2$, endowed with a metric invariant in one direction.This setting does not strictly speaking enter the framework of the present section, but is much simpler to describe (and we thus do not mention it in the present section).

\paragraph{Definition and differentiable structure.}
We assume that $(\calS , g)$ is an embedded 2D submanifold of $\R^3$ (endowed with the induced Euclidean structure), having $\S^1  = (\R/2\pi\Z) \sim SO(2)$ as an effective isometry group. 
The action of $\S^1$ on $\calS$, denoted by $\theta \mapsto \calR_\theta$ (such that $\calR_\theta \calS = \calS$) has:
\begin{enumerate}
\item exactly two fixed points denoted by $N, S \in \calS$ called North and South poles in Case~\ref{caseS}; we write $\d_N \calS = \{N\}$ and  $\d_S \calS = \{S\}$; 
\item exactly one fixed point denoted by $N \in \calS$ called North pole in Case~\ref{caseD}; in this case, we write $\d_N \calS = \{N\}$, and $\d_S \mathcal{S} = \d \mathcal{S}$ has a single connected component  (called ``south boundary'') which is also invariant by $\calR_\theta$;
\item no fixed point in Case~\ref{caseCY}; in this case $\d \mathcal{S}$ has two connected components denoted $\d_N\calS$ and $\d_S\calS$ (called ``north and south boundaries'') which are both invariant by $\calR_\theta$.
\end{enumerate}
We denote by $\mathcal{P}$ the set of poles, that is $\P=\{N,S\}$ in Case \ref{caseS}, $\mathcal{P}=\{N\}$ in Case \ref{caseD} and $\mathcal{P}=\emptyset$ in Case \ref{caseCY} and set 
\baln
\label{e:def-U}
U =\calS \setminus \P.
\naln
We now describe a nice parametrization of $(\calS , g)$, and, in particular, useful coordinates on the set $U$. We set $L= \dist_g(\d_N \calS, \d_S \calS)$ and denote by 
\begin{align}
\label{e:def-IL}
 I_L=(0,L) \text{  in Case \ref{caseS}}, \quad   I_L=(0,L] \text{ in Case \ref{caseD} }, \quad  I_L=[0,L] \text{ in Case \ref{caseCY}}.
\end{align}
We let $\gamma_0$ be a geodesic curve of $\calS$ joining $N$ (resp. $\d_N \calS$ in case~\ref{caseCY}) to $S$ (resp. $\d_S \calS$ in Cases~\ref{caseD} and~\ref{caseCY}). Note in particular that $\length(\gamma_0) = L$. For any $\theta \in \S^1$, the isometry $\calR_\theta$ transforms the geodesic $\gamma_0$ into $\calR_\theta (\gamma_0)$, which is another geodesic joining $N$ (resp. $\d_N \calS$) to $S$ (resp. $\d_S \calS$). For every $m \in U$ (defined in~\eqref{e:def-U}), there exists a unique $\theta \in \S^1$ such that $m$ belongs to $\calR_\theta (\gamma_0)$. The geodesic $\calR_\theta (\gamma_0)$ can be parametrized by arclength
\begin{align*}
 \rho : [0,L]\to  \calR_\theta (\gamma_0) ,  \quad \rho (0) \in \d_N \calS , \quad \rho (L)  \in \d_S \calS ,  \\
  s = \dist_g( \rho(s) , \d_N \calS) = L-  \dist_g( \rho(s) , \d_S \calS),
\end{align*}
and there exists a unique $s \in I_L$ such that $\rho(s) = m$. We use $(s,\theta)$ as a parametrization of $U \subset \calS$:
$$
\begin{array}{cccc}
\zeta : &  U = \mathcal{S} \setminus \P& \to & I_L \times \S^1\\
& m& \mapsto &\zeta(m) =  (s, \theta).
\end{array}
$$
In Case~\ref{caseCY}, $\P=\emptyset$ and thus the whole surface $\calS = U$ is diffeomorphic to the cylinder $I_L \times \S^1$ via $\zeta$. 
In Cases~\ref{caseCY} and~\ref{caseS}, we further need to describe coordinate charts around the poles.
In cases~\ref{caseS} and~\ref{caseD}, we define another exponential chart $(U_N, \zeta_N)$ centered at the pole $N$ by
\begin{align*}
U_N = \{N\} \cup \zeta^{-1} \left( \big( 0, \frac{L}{2} \big) \times \S^1 \right) = B_g\left(N, \frac{L}{2}\right)  \subset \calS, 
\\ \zeta_N : U_N \to B_{\R^2}\left(0 , \frac{L}{2}\right) , \quad \zeta_N (N) = 0 .
\end{align*}
with the transition map
$$
\begin{array}{cccc}
\zeta_N \circ \zeta^{-1} : & \zeta \big(U\cap U_N \big) = \left( 0,\frac{L}{2} \right) \times \S^1  & \to &  \zeta_N \big(U\cap U_N \big) =   B_{\R^2}\left(0 , \frac{L}{2}\right) \setminus \{ 0 \}\\
& (s ,\theta)& \mapsto &  \big(s \cos(\theta), s \sin (\theta) \big).
\end{array}
$$
In Case~\ref{caseS}, we add similarly a last exponential chart $(U_S, \zeta_S)$ centered at the pole $S$ by
\begin{align*}
U_S = \{S\} \cup \zeta^{-1} \left( \big( \frac{L}{2}, L \big) \times \S^1 \right) = B_g\left(S, \frac{L}{2}\right) \subset \calS ,
\\ \zeta_S : U_S \to B_{\R^2}\left(0 , \frac{L}{2}\right) , \quad \zeta_S (S) = 0 ,
\end{align*}
with the transition map  
$$
\begin{array}{cccc}
\zeta_S \circ \zeta^{-1} : & \zeta \big(U\cap U_S \big) = \left(\frac{L}{2} , L \right) \times \S^1  & \to &  \zeta_S \big(U\cap U_S \big) =   B_{\R^2}\left(0 , \frac{L}{2}\right) \setminus \{ 0 \}\\
& (s ,\theta)& \mapsto &  \big((L-s) \cos(\theta), (L-s) \sin (\theta) \big).
\end{array}
$$
We shall need the following notation. For a subset $J\subset I_L$, we denote by
\bal
\CO{J}=\zeta^{-1}(J\times\S^1) = \left\{m=\zeta^{-1}(s,\theta)\in U;s\in J\right\} \subset U\subset \calS 
\nal 
the $\calR_\theta$ invariant set which projects downto $J$.
We will also extend this definition to sets $J\subset [0,L]$ by adding the point $N$ if $0\in J$ (in Cases \ref{caseS} and~\ref{caseD}) and  the point $S$ if $L\in J$ (in Case \ref{caseS}).

\paragraph{Riemannian structure and operators involved.}
On the cylinder $I_L\times \S^1$, the metric $g$ is given by 
\begin{align}
\label{e:metric-revol}
(\zeta^{-1})^*g = ds^2 + R(s)^2 d\theta^2 ,
\end{align}
for some smooth function $R : I_L \to \R^+_*$ (the function $R$ can be interpreted as the Euclidean distance in $\R^3$ of the point parametrized on $\calS$ by $s$ to the symmetry axis, see e.g.~\cite[Section~3]{LL:18}). Since $g$ is a smooth metric on $\calS$, \cite[Proposition~4.6]{Besse} gives that $R$ extends to a $C^\infty$ function $[0,L] \to \R^+$ satisfying moreover
\baln
\label{e:condR}
R(0) = 0 , \quad R'(0) = 1 , \quad R^{(2p)}(0) =0 \quad \text{for any } p \in \N , \quad \text{ in Cases~\ref{caseS} and~\ref{caseD}} ,  \\ 
R(L) = 0 ,\quad R'(L) = -1 , \quad R^{(2p)}(L)=0 \quad \text{for any } p \in \N , \quad \text{ in addition in Case~\ref{caseS}} .  \nonumber
\naln

For other parametrizations of surfaces of revolution, or direct computations on the sphere $\S^2$ and the disk $\ID$, we refer to~\cite[Section~3]{LL:18}.

\begin{example}
In particular, we consider here the following three examples:
\begin{itemize}
\item the unit sphere of $\R^3$ is given by case~\ref{caseS} $L=\pi$, $s\in (0,\pi)$, $R(s) = \sin(s)$ (and the maximum of $R$ is reached at $s_0= \frac{\pi}{2}$)
\item the unit disk of $\R^2$ is given by case~\ref{caseD} with $L=1$ and $R(s)=s$;
\item flat cylinder of length $L>0$ and radius $R_0>0$ is given by case~\ref{caseCY} with $R(s)= R_0$.
\end{itemize}
\end{example}
In these coordinates, the Riemannian volume form is hence $R(s) ds d\theta$, the Riemannian gradient of a function is
\baln
\label{e:reim-gradient}
\nabla_g u = \d_s u \frac{\d}{\d s} + \frac{1}{R(s)^2} \d_\theta u \frac{\d}{\d \theta} , \quad \text{ with } \quad g(\nabla_g u , \nabla_g u ) = |\d_s u|^2 + \frac{1}{R(s)^2} |\d_\theta u|^2 ,
\naln
and the Laplace-Beltrami operator is given by
\bal
\Delta_{s,\theta} = \frac{1}{R(s)} \d_s R(s) \d_s + \frac{1}{R(s)^2} \d_\theta^2 .
\nal
We define by $L^2(\calS) := L^2(\calS, d\Vol_g)$ the space of square integrable functions, which is also invariant by the action of $(\calR_\theta)_{\theta \in \S^1}$. We will sometime also use the same definition for $L^2(\CO{J}) := L^2(\CO{J}, d\Vol_g)$ for $J\subset [0,L]$.

Another important operator is the infinitesimal generator $X_\theta$ of the group $(\calR_\theta)_{\theta \in \S^1}$, defined, for $u \in C^\infty(\calS)$, by 
\baln 
\label{e:Xtheta}
X_\theta u = \lim_{\vartheta \to 0} \vartheta^{-1} (u \circ \calR_\vartheta -u).
\naln In the chart $(U, \zeta)$, the action of $\calR_\theta$ is given by $(\zeta^{-1})^* \calR_\theta (u, \theta') = (u, \theta' +\theta)$, so that $(\zeta^{-1})^* X_\theta  = \d_\theta$. 
It is proved in~\cite[Section~3.2]{LL:18} that $X_\theta$ is a smooth vector field on $\calS$. 
Note also that $X_\theta(N) = X_\theta(S)=0$ and that its norm is given by $\sqrt{g(X_\theta, X_\theta) (s, \theta)}= R(s)$ (in the coordinates of $U$). 

Now, remark that $(\calR_\theta)_{\theta \in \S^1}$ acts as a (periodic) one-parameter unitary group on $L^2(\calS)$ by $f \mapsto f\circ \calR_\theta$. The Stone Theorem (see e.g.~\cite[Theorem~VIII-8~p266]{RS:I}) hence implies that its infinitesimal generator is $i A$, where $A$ is a selfadjoint operator on $L^2(\calS)$ with domain $D(A) \subset L^2(\calS)$. Since $i Af = X_\theta f$ for $f \in C^\infty(\calS)$ (which is dense in $D(A)$) according to~\eqref{e:Xtheta}, we have that $A$ is the selfadjoint extension of $\frac{X_\theta}{i}$. From now on, we slightly abuse the notation and still denote $\frac{X_\theta}{i}$ for its selfadjoint extension $A$.

\paragraph{Gradient vector field and conjugated operator.}
We finally introduce a function $\f : \calS \mapsto \R$, at least of class $W^{2,\infty}$ to define the gradient flow. Throughout this section, we assume that $X_\theta \f = 0$, i.e. $\f$ is invariant by rotation and the same property holds for $q$. In the coordinates of $U$, we shall thus simply write $\f=\f(s)$. These regularity assumptions can be written in the coordinate of $U$ by
\baln
\label{e:reg-f}
& s\mapsto \f(s)  \in W^{2,\infty}(0,L) , \text{ with } \nonumber\\
&  \f'(0) = 0 \text{ (in Cases~\ref{caseS} and~\ref{caseD}),} \quad
 \text{and } \f'(L) = 0 \text{ (in Case~\ref{caseS})}. 
\naln

\medskip
We may now define as in~\eqref{e:def-Peps} the conjugated operator $P_\eps$ as 
\begin{align}
\label{e:defPeps}
P_\eps & = - \eps^2 \Lap + \frac{|\nablag \f|_g^2}{4} + \eps \frac{\Lap \f}{2} -\eps \q\nonumber \\
 & = - \eps^2 \left( \frac{1}{R(s)} \d_s R(s) \d_s + \frac{1}{R(s)^2} \d_\theta^2  \right) +  \frac{|\f'(s)|^2}{4}  +  \frac{\eps}{2} \frac{1}{R(s)} \d_s (R(s)\f'(s)) -\eps \q,
\end{align}
where the second writing, in the coordinates of $U$, uses the invariance of $\f$. Note that the last term in this expression acts as a multiplication operator by a function in $L^\infty(\calS)$ with size $\eps$. We shall often consider it as a lower order term, and keep the shorter notation $\Lap \f$ in place of $\frac{1}{R(s)} \d_s (R(s)\f'(s))$.

Since both $g$ and $\f$ are invariant by the action of $\calR_\theta$, we have 
\begin{equation}
\label{e:commutation}
[X_{\theta} ,  P_\eps ] =0 .
\end{equation}
Moreover, $P_\eps$ is selfadjoint in $L^2(\calS, d\Vol_g)$ with domain $H^2(\calS) \cap H^1_0(\calS)$ ($=H^2(\calS)$ in Case~\ref{caseS}), and has compact resolvent. Therefore, the operators $P_\eps$ and $X_\theta$ share a common basis of eigenfunctions (see e.g.~\cite[Section~3.2]{LL:18} for a proof).
 If $\lambda \in\R$ is an eigenvalue of $P_\eps$, then (in the coordinates of $U$) the associated eigenfunction can be written as $e^{ik \theta}v(s)$ with $k\in \Z$, $v \in H^2_{\loc}(I_L)\cap L^2\left((0,L), R(s)ds\right)$ satisfying
\baln
\label{equationk1D}
- \frac{\eps^2 }{R(s)} \d_s\left( R(s) \d_s v \right) + \eps^2\frac{k^2}{R(s)^2} v+  \left( \frac{|\f'(s)|^2}{4}  +  \eps \qf\right) v =\lambda v ,
\naln
together with $v(L)=0$ in Case~\ref{caseD} and $v(0)=v(L)=0$ in Case~\ref{caseCY}. 

\bigskip
Restoring the dependence of the eigenelements in the parameter $\eps$, we call the normalized eigenfunctions of $P_\eps$: $\varphi_{k,n}^\eps=e^{ik \theta} v_{k,n}^\eps(s)$ with eigenvalues $\lambda_{k,n}^\eps$, where $n\in \N$.  In particular, for all $\eps>0$, we can write $L^2(\calS)=  \oplus^{\perp }_{(k,n)\in \Z\times \N}\vect (\varphi_{k,n}^\eps) $.

We further denote 
$$
L^2_k=\ker (X_{\theta}-ik)=\left\{\varphi\in L^2(\calS);\varphi_{|U} =e^{ik \theta}f(s),f \in L^2\left((0,L),R(s)ds\right)\right\},
$$ and $H^2_k=D(P_\eps) \cap L^2_k = H^2(\calS)\cap H^1_0(\calS)\cap L^2_k$. The commutation property~\eqref{e:commutation} implies that for all $\eps>0$,  $P_\eps H^2_k\subset L^2_k$, so we can define the operator 
\baln\label{defPk}P_{\eps}^{(k)} =P_{\eps\left|L^2_k\right.} , \quad \text{with domain } H^2_k,\naln which is selfadjoint. This can be seen for instance directly on the simultaneous diagonalization which implies an isometry $L^2(\calS) \approx \ell^2(\Z\times \N)$ where $L^2_k\approx \left\{(k,n)\left|n\in \N\right.\right\}$ as a closed subspace of $\ell^2(\Z\times \N)$. The fact that $P_{\eps}$ has compact resolvent implies that this is also the case for $P_{\eps}^{(k)}$.
With a slight abuse of notation, we shall still denote by $P_{\eps}^{(k)}$ the one dimensional operator $(\zeta^{-1})^*P_{\eps}^{(k)} \zeta^*$ defined on $I_L$, namely 
\baln
\label{equationk1D-bis}
P_{\eps}^{(k)} w = - \frac{\eps^2 }{R(s)} \d_s\left( R(s) \d_s w \right) +  \left(\frac{ \eps^2 k^2}{R(s)^2} + \frac{|\f'(s)|^2}{4}  + \eps \qf\right) w .
\naln

\subsection{The conditions (GCC) and (FC) on surfaces of revolution}
\label{s:minimal-time-limit-equation}
 In this section, we characterize the conditions (GCC) (see Definition~\ref{d:def-GCC} if $\d\M=\emptyset$) and (FC) (see Definition~\ref{d:def-GCC_bdry} if $\d\M\neq \emptyset$) in the above very particular geometry, and further assuming that the observation region $\omega$ is rotationally invariant as well. In case $\d\M \neq \emptyset$, we consider the two cases of internal and boundary observation, and describe the associated minimal times $T_{FC}$.

 \begin{proposition}
 \label{p:GCCFCrot}
 Let $\delta>0$ and recall that $\f$ is assumed to be $\theta$-invariant.
 \begin{enumerate}
 \item In Case~\ref{caseS}, we consider the set $\omega = B_g(N, \delta) \cup B_g(S, \delta)$; then $(\calS,\nabla_g\f,\omega)$ satisfies (GCC) if and only if $\f'(s) \neq 0$ for all $s \in [\delta, L-\delta]$ and $T_{GCC}(\calS,\nabla_g\f,\omega) = \left| \int_\delta^{L-\delta} \frac{ds}{\f'(s)} \right|$.
 \item In Case~\ref{caseD}, we consider the set $\omega = B_g(N, \delta)$; then $(\calS,\nabla_g\f,\omega)$ satisfies (FC) if and only if $\f'(s) \neq 0$ for all $s \in [\delta, L]$ and  $T_{FC}(\calS,\nabla_g\f,\omega) = \left| \int_\delta^{L} \frac{ds}{\f'(s)}  \right|$. 
 \item In Case~\ref{caseCY}, $(\calS,\nabla_g\f)$ satisfies (FC) if and only if $\f'(s) \neq 0$ for all $s \in [0, L]$ and $T_{FC}(\calS,\nabla_g\f) =  \left| \int_0^{L} \frac{ds}{\f'(s)} \right|$.
  \end{enumerate}
 \end{proposition}
 Note that in Case~\ref{caseD}, the situations $\f'>0$ and $\f'<0$ play two different roles (see the proof below). Indeed, in case $\f'>0$, all trajectories of $-\f'$ enter the controlled region $\omega$, whereas, in case $\f'<0$, all trajectories of $-\f'$ flow out of the domain $\calS$ through $\d \calS$ (without passing into $\omega$).
 However, the definition of (FC) in Definition~\ref{d:def-GCC_bdry} does not make a distinction between these two situations.
 
 \bnp
 We only prove the second item; the other two items are proved similarly.
 According to Definition~\ref{d:def-GCC_bdry},~\eqref{e:reim-gradient} and the $\theta$-invariance of $\f$, it suffices to check under which conditions the solutions to $\dot{s}(t)= - \f'(s(t))$ all enter $\omega = B_g(N, \delta)$ (resp. all exit $\calS$, that is satisfy $s(T)>L$). 
 If there is $s_0 \in [\delta , L]$ such that $\f'(s_0) = 0$, then the associated solution satisfies $s(t) = s_0\in [\delta , L]$ for all $t\in \R$, and $(\calS,\nabla_g\f,\omega)$ does not satisfy (FC). 

If $\f'>0$ on $[\delta , L]$, then $s(t)$ is decreasing, and for any $\sigma_0,\sigma_1\in \R$, one has $\int_{s(\sigma_0)}^{s(\sigma_1)} \frac{ds}{\f'(s)} =\sigma_0- \sigma_1$. The longest trajectory that does not enter $\omega$ is such that $s(0) = L$ and $s(T) =\delta$ so that $T = \int_\delta^{L} \frac{ds}{\f'(s)}$. This proves   $T_{FC}(\calS,\nabla_g\f,\omega) =  \int_\delta^{L} \frac{ds}{\f'(s)}$ in this case.

Finally, if $\f' <0$ on $[\delta , L]$, then $s(t)$ is increasing, and for any $\sigma_0,\sigma_1\in \R$, one has $\int_{s(\sigma_0)}^{s(\sigma_1)} \frac{ds}{\f'(s)} =\sigma_0- \sigma_1$. The longest trajectory that does not enter $\omega$ is such that $s(0) = \delta$ and $s(T) =L$ so that $-T = \int_\delta^{L} \frac{ds}{\f'(s)}$. This proves  $T_{FC}(\calS,\nabla_g\f,\omega) =  - \int_\delta^{L} \frac{ds}{\f'(s)}$ in this case, and hence the proposition.
  \enp

  \subsection{Existence of eigenfunctions}
  \label{s:revol-eig}
 One may consider different asymptotic regimes in the parameters $\eps\to 0^+$ and $k \to + \infty$. Note that the case $k$ bounded would correspond to the one-dimensional situation treated in the companion paper~\cite{LL:20-1D}.
Here, we shall consider the limit $k \to + \infty$ and make the following choice of the parameter $\eps$:  
\baln
\label{defe}
\e=\e_k=\ce k^{-1}
\naln 
considered as a semiclassical parameter, where $\ce> 0$ is a fixed parameter (i.e. which does not depend on $k$) that will be chosen but fixed. All constants that appear below might depend on $c$. The analysis of the asymptotic of the constant involved as $c\to0$ (low level of rotation) or $c\to\infty$ (high level of rotation) would be interesting but would require much more work.

In view of~\eqref{equationk1D-bis}, the choice~\eqref{defe} naturally leads to consider
\baln
\label{e:def-Vc}
s \mapsto V_\ce (s) := \frac{\ce^2}{R(s)^2} + \frac{| \f'(s)|^2}{4} ,
\naln
as the effective potential of the operator $P_{\eps_k}^{(k)}$ in the semiclassical limit $\eps= \eps_k = c k^{-1} \to 0^+$. In particular, the operator $P_{\eps}^{(k)} $ is now a semiclassical operator with small parameter $\e_{k}$ and \eqref{equationk1D-bis} can be rewritten
\baln
\label{equationk1D-bisbis}
P_{\eps_{k}}^{(k)} w = - \frac{\eps_k^2 }{R(s)} \d_s\left( R(s) \d_s w \right) +  \left(V_\ce (s)  + \eps_{k} \qf\right) w .
\naln

\medskip
In the present section, we recall the existence of eigenfunctions (Analogue of Lemma~\ref{l:exist-eig} above) associated to any value of the effective potential $V_\ce$. More precisely, in the chosen regime~\eqref{defe}, for any $s_0\in I_L$ (recall that $I_L$ is defined in~\eqref{e:def-IL}), we construct a sequence $\psi_k$ such that $P_{\e_k} \psi_k = (V_\ce(s_0) + r(k)) \psi_k$, with $r(k)\to 0$ as $k\to + \infty$. As in Section~\ref{s:sub-eig-schro}, the precision $r(k)$ might depend on whether $V_\ce'(s_0) \neq 0$ or  $V_\ce'(s_0) = 0$ but we will only state the worst estimate, which is sufficient for our needs. We shall later on prove localization properties of the $\psi_k$'s assuming further that $V_c(s_0) = \min_{I_L}V_\ce$ (which is a global assumption).
We recall the choice~\eqref{defe} and the definition~\eqref{e:def-Vc}.

\begin{lemma}[Existence of eigenfunctions]
\label{l:exist-mu-psibord}
For all $s_0 \in I_L$, there is $k_0>0$ such that for all $k\in\N, k\geq k_0$, there exists $\psi_k \in H^{2}(\calS)\cap L^2_k$ in case~\ref{caseS} (resp. $\psi_k \in H^{2}(\calS)\cap H^1_0(\calS)\cap L^2_k$ in cases~\ref{caseD}-\ref{caseCY}),  and $\mu_k \in \R$ such that 
\begin{align*}
\mu_k=  V_\ce(s_0) +\grando{\frac{1}{k^{2/3}}} 
=  \frac{\ce^2}{R(s_0)^2}+\frac{|\f'(s_0)|^2}{4} + \grando{\frac{1}{k^{2/3}}} , \\
P_{\e_k} \psi_k = \mu_k \psi_k , \quad \|\psi_k\|_{L^2(\calS)}=1 , \quad \psi_k(s,\theta) = e^{ik\theta} \varphi_k(s) , 
\end{align*}
with $P_\eps$ defined in~\eqref{e:defPeps} and $\e_k$ in~\eqref{defe}.
\end{lemma}
The proof is very similar to Lemma~\ref{l:exist-eig}. Indeed, the study of the operator $P_{\e_k}$ in $L^{2}_{k}$ reduces to the 1D problem for the operator $P_{\e}^{(k)}$  defined in \eqref{defPk} and \eqref{equationk1D-bisbis}. The proof consists first in constructing quasimodes exactly as in the proof of Lemma~\ref{l:exist-eig}. Deducing existence of an exact eigenfunction from a quasimode requires the use of the right selfadjoint extension $P_{\e_{k}}^{(k)}$. This issue is however treated in detail in~\cite[Lemma 3.6]{LL:18}.
Note that $s_{0} \in I_L$ implies that it cannot be a pole so that $R(s_{0})^{-1}$ is finite.

\subsection{Geometric assumptions and the Agmon distance}
The next step is to study the behavior of the eigenfunction $\psi_k$ constructed in Lemma~\ref{l:exist-mu-psibord}. This will require some global assumptions on the effective potential $V_\ce$. Recall that $V_\ce$ is defined in~\eqref{e:def-Vc}, is continuous on $I_L$, and tends to $+\infty$ near to the poles. Indeed, in Cases~\ref{caseS} and~\ref{caseD} we have for instance 
$$
V_c(s)\sim_{s\to 0^+} \frac{c^2}{R(s)^2}\sim_{s\to 0^+}  \frac{c^2}{s^2} \to_{s\to 0^+}   +\infty ,
$$ as a consequence of~\eqref{e:condR} and~\eqref{e:reg-f} (and similarly when $s\to L^-$ in Case~\ref{caseS}).
As a consequence, $V_\ce$ admits a minimum on the interval $I_L$, which we denote by
$$
V_{\min} = \min_{s \in I_L} V_\ce(s) \in \R_+^* .
$$
In the following, we make the following assumption (a precised version of~\eqref{s:asspt-uniq-min}) on the set where $V_\ce$ reaches its global minimum.
\begin{assumption}
\label{a:single-min}
The set $V_{\ce}^{-1}(V_{\min}) = \{s_{\min}\} \subset I_L$ consists in a single point.
\end{assumption}
Note that this is assumption is generic. Here it is not strictly needed to prove the main results, but simplifies the presentation and statements slightly. 
We again introduce the relevant Agmon distance at the minimal energy level $V_{\min}$, defined in the coordinates of $U$ by the eikonal equation 
\bal
 \big((d_A^{\ce})'(s) \big)^2-\left( V_\ce(s)- V_\ce(s_{\min}) \right) = 0 , \quad d_A^{\ce}(s_{\min})=0 ,  \quad \sgn((d_A^{\ce})'(s)) =\sgn(s-s_{\min}) ,
\nal
or, more explicitly, for $s\in I_L$, by~\eqref{e:defbisdA}.
In view of the $W^{2,\infty}$ regularity of $\f$ on $\calS$ and the definition of $V_\ce$ in~\eqref{e:def-Vc}, the function $d_A^{\ce}$ is of class $C^2$ away from $s_{\min}$, $0$ and $L$ and is locally Lipschitz on $I_L$. Note that this includes Lipschitz regularity up to the boundary $s=L$ in Case \ref{caseD} and to both boundaries $s=0,L$ in Case \ref{caseCY}. 
We also consider $d_A^{\ce}$ as a $\theta$-invariant function on the surface $\calS$.

\begin{lemma}[Properties of $d_A^{\ce}$]
\label{lemma-prop-dA}
Under Assumption~\ref{a:single-min}, we have $d_A^{\ce} \in C^2(I_L\setminus \{s_{\min}\})$ together with 
\baln
\label{equivdAlog}
d_A^{\ce}(s) = -\ce \log(s) + \grando{1}, \quad \text{as } s \to 0^+ , \quad \text{ in Cases~\ref{caseS} and~\ref{caseD}, }\\
\label{equivdAlog-bis}
\qquad d_A^{\ce}(s)= -\ce \log (L-s) + \grando{1}, \quad \text{as } s \to L^- , \quad \text{ in Case~\ref{caseS}}.
\naln
\end{lemma}
\bnp
We only consider the asymptotics as $s\to 0^+$, that is, prove~\eqref{equivdAlog}. The proof of~\eqref{equivdAlog-bis} follows the same.
Remark that according to \eqref{e:condR}, we have $\frac{1}{R(y)}\to + \infty$ as $y\to 0^+$ with 
$$
R(s) = s + \grando{s^3} , \text{ when } s \to 0^+ .
$$
As a consequence, since $\f \in W^{2,\infty}(I_L)$, we have 
$$
V_\ce (s) = \frac{\ce^2}{s^2} + \grando{1} , \text{ when } s \to 0^+  .
$$
With~\eqref{e:defbisdA}, we obtain $d_A^{\ce}(s) = \left| \int^s_{s_{\min}} \frac{\ce}{y}(1+\grando{y}) dy \right| = -\ce \log(s) + \grando{1}$ as $s \to 0^+$.
\enp

\subsection{Upper bounds for eigenfunctions: Agmon estimates}
\label{s:surf-Agmon}
As far as upper bounds on $\psi_k$ are concerned, we have the following Agmon type estimate.
\begin{proposition}
\label{p:Agmon-revolution}
Under Assumption~\ref{a:single-min}, assume that $\mu_k = V_\ce(s_{\min}) +r(k)$ with $r(k)\to 0$ and $\psi_k \in H^2 (\mathcal{S})\cap L^2_k(\calS)$ solves
$$
P_{\e_k}\psi_k=\mu_k \psi_k, \text{ on } \calS , \quad \psi_k|_{\d\calS} = 0 ,\quad \nor{\psi_k}{L^2(\calS)}=1 ,
$$
with $P_\eps$ defined in~\eqref{e:defPeps} and $\e_k$ in~\eqref{defe}. Then for all $0<\delta< 1$, there exist $C=C(\delta) ,k_0=k_0(\delta)>0$ such that, for all $k \in \N$, $k \geq k_0$, the following integral is well defined with the estimate
$$
\int_{\calS} e^{2(1-\delta)\frac{ d_A^{\ce}(m)}{\e_k}} |\psi_k|^2(m) d\Vol_g(m)  
\leq C e^{2\frac{\delta}{\e_k}}  .
$$
Also, if $\d\calS \neq \emptyset$ (that is, in Cases~\ref{caseD} and~\ref{caseCY}), for all $0<\delta< 1$, there exist $C=C(\delta) ,k_0=k_0(\delta)>0$ such that, for all $k \in \N$, $k \geq k_0$, 
\baln
\label{e:estimate-psik-bord}
\nor{\d_s \psi_k(L,\cdot)}{H^1(\S^1)}^2
\leq C e^{- 2(1-\delta)\frac{ d_A^{\ce}(L)}{\e_k}} , \quad 
\nor{\d_s \psi_k(0,\cdot)}{H^1(\S^1)}^2
\leq C e^{- 2(1-\delta)\frac{ d_A^{\ce}(0)}{\e_k}} ,
\naln
where the last estimate (at $s=0$) holds true in Case~\ref{caseCY} only.
\end{proposition}
Note that given the asymptotic expansion of $d_A^{\ce}$  in Lemma~\ref{lemma-prop-dA}, this estimate implies that $\psi_k$ vanishes strongly near the poles of $\calS$.
 The proof is made with classical Agmon type identity with some care with respect to the degeneracy at the poles. It is very similar to the one performed in \cite[Theorem 3.9]{LL:18} and we omit it. Note that, as opposed to \cite[Theorem 3.9]{LL:18}, we do not assume here that the minimum be non-degenerate, and only deduce an estimate with loss ($\delta>0$), which is sufficient for our needs. In the non-degenerate case, one can take $\delta=0$ in this estimate and replace the right hand-side by a polynomial bound of the type $\e_k^{-M}$, see~\cite[Theorem 3.9]{LL:18}.

The proof of the boundary estimate also requires a bootstrap argument to estimate higher $H^s$ norms and the use of a trace estimate, see~\cite[Propositions~3.3.1 and~3.3.4]{Helffer:booksemiclassic}.

We obtain the  following two direct Corollaries. 
\begin{corollary}
\label{coragmonfaible}
Under the assumptions of Proposition~\ref{p:Agmon-revolution}, for all $0<\delta< 1$, there exist $C=C(\delta) ,k_0=k_0(\delta)>0$ such that, for all $k \in \N$, $k \geq k_0$, and for all rotationally invariant set $\omega$, we have
$$
\int_{\omega} |\psi_k|^2(m) d\Vol_g(m)  
\leq C e^{-\frac{2}{\e_k} ((1-\delta)d_A^{\ce}(\omega) - \delta) }  , \quad \text{ with } d_A^{\ce}(\omega) = \inf_{s\in \omega} d_A^{\ce}(s) .
$$
\end{corollary}
\bnp 
This is a direct consequence of the following estimate:
\begin{align*}
e^{2(1-\delta)\frac{ d_A^{\ce}(\omega)}{\e_k}} \int_{\omega}  |\psi_k|^2(m) d\Vol_g(m) & \leq\int_{\omega} e^{2(1-\delta)\frac{ d_A^{\ce}(m)}{\e_k}} |\psi_k|^2(m) d\Vol_g(m) \\
& \leq \int_{\calS} e^{2(1-\delta)\frac{ d_A^{\ce}(m)}{\e_k}} |\psi_k|^2(m) d\Vol_g(m)  \leq  C e^{2\frac{\delta}{\e_k}} ,
\end{align*}
where we have used Proposition~\ref{p:Agmon-revolution} in the last inequality.
\enp
\begin{corollary}[most of the norm is close to the minimum]
\label{corpresx0}
Under the assumptions of Proposition~\ref{p:Agmon-revolution}, for any $\rho>0$, there exists $k_0\in \N$ so that 
\bal
\nor{\psi_k}{L^2(\CO{(s_{\min}-\rho , s_{\min} +\rho)})}^2\geq 1/2, \quad \text{ for all } k\in \N, k\geq k_0.
\nal
\end{corollary}
\bnp
Applying Corollary \ref{coragmonfaible} with $\omega=\mathcal{S}\setminus\CO{(s_{\min}-\rho , s_{\min} +\rho)}$ ensures that for any $\delta>0$, there is $k_0\in \N$ such that for all $k\geq k_0$, we have
$$
\int_{\omega} |\psi_k|^2(m) d\Vol_g(m)  
\leq C(\delta) e^{-\frac{2}{\e_k} ((1-\delta)d_A^{\ce}(\omega) - \delta) }  , \quad \text{ with } d_A^{\ce}(\omega) = \inf_{s\in \omega} d_A^{\ce}(s) .
$$
From Assumption~\ref{a:single-min}, we have $d_A^{\ce}(s)>0$ for all $s \in I_L\setminus \{s_{\min}\}$. That $\omega$ is closed and does not contain $s_{\min}$ implies that $d_A^{\ce}(\omega)>0$. 
 Then we fix $\delta>0$ small enough so that $(1-\delta)d_A^{\ce}(\omega) - \delta>0$. There is $k_0\in \N$ such that we have $\nor{\psi_k}{L^2(\omega)}^2 \to 0$ for $k \geq k_0$, which implies the result.
\enp

\subsection{Lower bounds for eigenfunctions: Allibert estimates}
\label{sectlower}
In Corollary~\ref{coragmonfaible}, we proved that the family of eigenfunctions $\psi_k$ decays on $\omega$ at least like $e^{-\frac{d_A^{\ce}(\omega)}{\eps_k}}$.
The purpose of this section is to prove the converse, i.e. that the $\psi_k$'s decay {\em at most} like $e^{-\frac{d_A^{\ce}(\omega)}{\eps_k}}$ up to $\delta$ loss.
This comes from the particular one-dimensional underlying context.
We follow in this section the method of Allibert~\cite{Allibert:98}. 
More precisely, we prove the following estimates. 
\begin{proposition}
\label{propinfdA}
Under Assumption~\ref{a:single-min}, assume further that 
 \baln \label{rappelmuinf}\mu_k = V_\ce(s_{\min}) +r(k)\textnormal{ with }r(k)\to 0 , 
 \naln
 and $\psi_k \in H^2 (\mathcal{S})\cap L^2_k(\calS)$ solves
$$
P_{\e_k}\psi_k=\mu_k \psi_k, \text{ on } \calS , \quad \psi_k|_{\d\calS} = 0 ,\quad \nor{\psi_k}{L^2(\calS)}=1 , \quad \psi_k(s, \theta) = e^{ik \theta}\varphi_k(s) ,
$$
with $P_\eps$ defined in~\eqref{e:defPeps} and $\e_k$ in~\eqref{defe}.

Then, for any $\eta,\delta>0$, there exist $k_0,C>0$ so that 
\bal
\nor{\psi_k}{L^2(\CO{B(s,\eta)})}\geq C e^{-\frac{1}{\e_k}(d_A^{\ce}(s)+\delta)}, \quad \text{ for all } k\geq k_0\text{ and }s \in I_L \text{ s.t. }\dist(s,\mathcal{P} )>\eta.
\nal
\end{proposition}
Note that $\eta$ is a safety distance to the set of poles $\P$ defined in Section \ref{s:revol}.
The proof of Proposition~\ref{propinfdA} relies on two lemmata, which we give in the next section. 

\subsubsection{Two preliminary lemmata}
In this section, we assume that the assumptions of Proposition~\ref{propinfdA} are satisfied. In particular, the eigenfunctions under consideration are of the form $\psi_k(s, \theta) = e^{ik \theta}\varphi_k(s)$.  
We define the following ``semiclassical energy densities'' of the eigenfunctions $\psi_k$, for $s \in I_L$ by 
\begin{align}
\label{e:defEk}\Ek(s)&:= \eps_k^2 |\d_s \psi_k|^2(s)+ (V_\ce(s) -\mu_k+1)|\psi_k|^2(s)  \\
& =  \eps_k^2 |\varphi_k'|^2(s)+ (V_\ce(s) -\mu_k+1)|\varphi_k|^2(s) ,\nonumber \\
\label{e:defEk+} \Ek^+(s)&:=\eps_k^2 |\d_s \psi_k|^2(s)+(V_\ce(s) -\mu_k)|\psi_k|^2(s) \\
& =\eps_k^2 |\varphi_k'|^2(s)+(V_\ce(s) -\mu_k)|\varphi_k|^2(s) . \nonumber
\end{align}
Note that according to elliptic regularity, we have $\varphi_k \in H^2_{\loc}(I_L)$ and, due to Sobolev embeddings,  $\varphi_k', \Ek, \Ek^+ \in C^0(I_L)$ and in particular $\Ek, \Ek^+$ are defined everywhere on $I_L$.
For $s, t \in I_L$, we define $I_{s,t}$ to be the interval between the real numbers $s$ and $t$, that is, either $[s,t]$ or $[t,s]$.
 We also set 
 \baln
 \label{e:def-E-alpha}
 E_\alpha = \{ s \in I_L ; |s-p| \geq \alpha, \text{ for all } p\in\mathcal{P}, \text{ and }|s-s_{\min}|\geq \alpha  \}.
 \naln 
 Recall that $\mathcal{P}$ is the set of poles defined in Section \ref{s:revol} (and is aimed at covering all Cases~\ref{caseS}--\ref{caseCY}
 at the same time) and $s_{\min}$ is the point at which $V_{\ce}$ reaches its minimum.
\begin{lemma}
\label{GronwallAgmon}
Assume the assumptions of Proposition~\ref{propinfdA} and recall that $\E_k^+$ is defined in~\eqref{e:defEk+}.
Then, for any $\alpha,\delta>0$, there exists $k_0>0$ so that for all $s ,t \in I_L$ such that $I_{s,t} \subset E_\alpha$, we have
\bal
\Ek^+(t)\leq e^{\frac{2}{\e_k}\left(|d_A^{\ce}(s)-d_A^{\ce}(t)|+\delta\right)}\Ek^+(s), \quad \text{for all }k\geq k_0 .
\nal
\end{lemma}
Lemma~\ref{GronwallAgmon} provides with a Gr\"onwall type estimate on the energy $\Ek^+$, with a precise description of the constant, under the condition that we remain finitely away from $s_{\min}$. It is an analogue of~\cite[Lemma~12]{Allibert:98} in our setting (with an additional uniform dependence). 

Note that $|d_A^{\ce}(s)-d_A^{\ce}(t)| = d_A^{\ce}(s,t)$ is the Agmon distance between $s$ and $t$ at the lowest energy level.

\bnp[Proof of Lemma~\ref{GronwallAgmon}]
On the interval $z \in I_{s,t}\subset E_\alpha$, we differentiate $\Ek^+$, yielding
$$(\Ek^+)'(z)=2\eps_k^2\Re( \ovl{\varphi_k}' \varphi_k'')+ V_\ce'(z)\varphi_k^2 + 2(V_\ce(z)-\mu_k) \Re(\varphi_k \ovl{\varphi_k}' ) .$$
We recall the choice of $\eps_k$ in~\eqref{defe} from the definition of $P_{\eps_k}$ in~\eqref{e:defPeps}, and the definition of $P_{\eps_k}^{(k)}$ in~\eqref{equationk1D}-\eqref{equationk1D-bisbis} that we have 
$$
\mu_k \varphi_k =P_{\eps_k} \varphi_k  = - \eps_k^2 \varphi_k'' - \eps_k^2 \frac{R'}{R}\varphi_k' + V_\ce \varphi_k  +  \eps_k\qf \varphi_k .
$$
Replacing $\eps_k^2 \varphi_k''$ in the above identity yields
\begin{align}
\nonumber (\Ek^+)'(z)& =2  \left(V_\ce(z) -\mu_k + \eps_k\qf) \right) \Re(\varphi_k \ovl{\varphi_k}' )- 2\eps_k^2 \frac{R'}{R}|\varphi_k'|^2+ V'_c(z)|\varphi_k|^2 \\
& \quad + 2(V_{\ce}(z)-\mu_k) \Re(\varphi_k \ovl{\varphi_k}' ) \nonumber\\
\label{derivNRJ}& = \big(4(V_{\ce}(z)-\mu_k)+ 2 \eps_k\qf \big)  \Re(\varphi_k \ovl{\varphi_k}' ) - 2\eps_k^2 \frac{R'}{R}|\varphi_k'|^2+ V'_c(z)|\varphi_k|^2  .
\end{align}

First, using the continuity of $V_{\ce}$ on $I_L$ and the compactness of $E_{\alpha}$ in $I_L$, we see that $V_{\ce}$ reaches its minimum on $E_{\alpha}$. This, together with the fact that  $ s_{\min}\notin E_{\alpha}$, implies that  $C_{\alpha}^{-2}\leq V_{\ce}(s)-V_{\ce}(s_{\min})\leq C_{\alpha}^2$ uniformly for $s\in E_\alpha$. Recalling~\eqref{rappelmuinf}, this yields the existence of $k_0 (\alpha)$ such that for $k \geq k_0(\alpha)$, we have $V_{\ce}-\mu_k \geq \frac{1}{2C_\alpha}>0$ on $E_{\alpha}$.

We now estimate each of the terms in~\eqref{derivNRJ}. We first have
\begin{align*}
\left| 4(V_{\ce}-\mu_k)\Re(\varphi_k \ovl{\varphi_k}' )\right| 
 & \leq 4 \e_k^{-1} \sqrt{V_{\ce}-\mu_k}  \left( \eps_k |\varphi_k'| \right) \left(\sqrt{V_{\ce}-\mu_k} |\varphi_k| \right) \\
 & \leq 2 \e_k^{-1} \sqrt{V_{\ce}-\mu_k} \left[\e_k^2 |\varphi_k'|^2 + (V_{\ce}-\mu_k) |\varphi_k|^2 \right] \\
 & =2\e_k^{-1} \sqrt{V_{\ce}-\mu_k}  \Ek^+ .
\end{align*}
Moreover, according to~\eqref{rappelmuinf}, there exists a constant $C_\alpha>0$ such that we have $\sqrt{V_{\ce}-\mu_k} \leq \sqrt{V_{\ce}-V_{\ce}(s_{\min})} + C_\alpha |r(k)|$ uniformly for $s\in E_\alpha$.
Together with the above inequality, this implies
\begin{align*}
\left| 4(V_{\ce}-\mu_k)\Re(\varphi_k \ovl{\varphi_k}' ) \right|
\leq  \e_k^{-1}(2\sqrt{V_{\ce}-V_{\ce}(s_{\min})} + C_\alpha |r(k)|)\Ek^+ , \quad \text{ on } E_\alpha.
\end{align*}
Second, we have
 \begin{align*}
\left|V_{\ce}'|\varphi_k|^2\right| 
= \frac{\left|V_{\ce}'\right|}{V_{\ce}-\mu_k}  (V_{\ce}-\mu_k) |\varphi_k|^2
 \leq 2C_\alpha\nor{V'_{\ce}}{L^\infty(E_\alpha)}  \Ek^+ , \quad \text{ on } E_\alpha.
\end{align*}
 Third, we have 
 $$
 \left|\eps_k\qf \Re(\varphi_k \ovl{\varphi_k}' ) \right|\leq C_\alpha \Ek^+ ,  \quad \text{ on } E_\alpha.
$$ 
Finally, since $\frac{R'}{R}$ is bounded on $E_{\alpha}$, we have
$$
\left|\eps_k^2 \frac{R'}{R}|\varphi_k'|^2\right|\leq C_\alpha \eps_k^2 |\varphi_k'|^2
\leq  C_\alpha \Ek^+ \quad \text{ on } E_\alpha.
$$
Combining the last four estimates in~\eqref{derivNRJ} yields for another constant $C_\alpha>0$ and for all $k\geq k_0(\alpha)$
$$
\left| (\Ek^+)'(z) \right|\leq 2\e_k^{-1}\left[\sqrt{V_{\ce}(z)-V_{\ce}(s_{\min})} + C_\alpha|r(k)|+ C_\alpha \eps_k\right]\Ek^+(z) ,\quad \text{ for all }z \in E_\alpha .
$$
Applying the Gr\"onwall Lemma on the interval $I_{s,t}$ contained in $E_{\alpha}$ yields
$$
\Ek^+(s) \leq e^{\frac{2}{\eps_k}\left(\left|\int_t^s \sqrt{V_{\ce}(z)-V_{\ce}(s_{\min})}dz \right| + C_\alpha|r(k)|+ C_\alpha \eps_k\right)}\Ek^+(t) ,
$$
which is the sought result.
\enp
The next Lemma is aimed at giving a rough Gr\"onwall type estimate for the energy $\Ek$, without precise constants. The interest of this less precise result is that it remains true close to the minimum $s_{\min}$. This allows to compensate the fact that Lemma~\ref{GronwallAgmon} is not uniform when $s$ is close to $s_{\min}$.
Similarly as $E_\alpha$ in~\eqref{e:def-E-alpha}, we define 
 \bal
  F_\alpha = \{ s \in I_L ; |s-p| \geq \alpha, \text{ for all } p\in\mathcal{P}  \}.
 \nal 
\begin{lemma}
\label{Gronwallgene}
Assume the assumptions of Proposition~\ref{propinfdA} and recall that $\E_k$ is defined in~\eqref{e:defEk}.
For any $\alpha>0$, there exist $C_{\alpha},D_{\alpha}>0$ so that for all $s,t \in  F_\alpha$, we have
\bal
\Ek(s)\leq C_{\alpha}e^{\frac{2}{\e_k}D_{\alpha}|s-t|}\Ek(t) , \quad \text{for all }k\geq k_0 .
\nal
\end{lemma}
Lemma~\ref{Gronwallgene} is an analogue of~\cite[Lemma~11]{Allibert:98} in our setting. 
 Recall that $\mathcal{P}$ is the set of poles defined in Section \ref{s:revol}.

\bnp
The proof is quite close to that of Lemma \ref{GronwallAgmon}. We only use the fact that there exists $C_{\alpha}$ so that $C_{\alpha}^{-1}\leq V_\ce -\mu_k+1\leq C_{\alpha}$ on $F_\alpha$ if $k\geq k_{0}(\alpha)$. This gives a constant $D_\alpha>0$ such that $D_{\alpha}^{-1} \left( \e_k^{2}|\varphi_k'|^{2}+|\varphi_k|^{2}\right)\leq \Ek\leq D_{\alpha} \left(\e_k^{2} |\varphi_k'|^{2}+|\varphi_k|^{2}\right)$ on $F_\alpha$. The same computation as in~\eqref{derivNRJ} gives
\begin{align*}
\Ek'& =(\Ek^+)'+2\Re(\varphi_k \ovl{\varphi_k}' )  \\
& =\left(4(V_{\ce} -\mu_k) + r_k+1\right) \Re(\varphi_k \ovl{\varphi_k}' )  - 2\eps_k^2 \frac{R'}{R}|\varphi_k'|^2+ V'_c(z)|\varphi_k|^2 ,
\end{align*}
with $\sup_{F_\alpha}r_k \to 0$ as $k \to+ \infty$. 
As a consequence, for $k\geq k_{0}(\alpha)$, we have constants $C_\alpha',C_\alpha''$ such that for all $z\in F_\alpha$
\begin{align*}
|(\Ek)'(z)|& \leq \e_k^{-1}C_{\alpha}'  \left(\e_k^{2} |\varphi_k'|^{2}+|\varphi_k|^{2}\right)\leq  \e_k^{-1} C_{\alpha}''\Ek(z),
\end{align*}
which allows to conclude as in the proof of Lemma~\ref{GronwallAgmon} above by Gr\"onwall estimates.
\enp

\subsubsection{Proof of Proposition~\ref{propinfdA} from Lemmata~\ref{GronwallAgmon} and~\ref{Gronwallgene}}
The sketch of the proof of Proposition~\ref{propinfdA} is as follows: \begin{enumerate}
\item the total mass of $\varphi_k$ is dominated by its mass near the minimum $s_{\min}$ via Corollary~\ref{corpresx0}; 
\item the mass near $s_{\min}$ is dominated by the energy at $s_{\min}-\rho$ via Lemma~\ref{Gronwallgene} (with a small loss if $\rho$ is small);
\item the energy at $s_{\min}-\rho$ is dominated by the energy near $s$ via Lemma~\ref{GronwallAgmon} (with a geometric constant $e^{\frac{2}{\eps_k}d_A^{\ce}(s)}$);
\item the energy near $s$ is dominated by the $L^2$ norm of $\varphi_k$ (or $\psi_k$) near $s$ via elliptic regularity.
\end{enumerate}
\bnp[Proof of Proposition \ref{propinfdA}]
Without loss of generality, we can assume $s<s_{\min}-\eta$. Indeed, the case $s>s_{\min}+\eta$ is treated similarly, and the case $s\in [s_{\min}-\eta, s_{\min}+\eta]$ is a direct consequence of Corollary~\ref{corpresx0} applied for $\rho = \eta/2$.

Now, since $s_{\min} \in I_L$, notice that we have necessarily $d:=\dist(s_{\min},\mathcal{P})>0$, and we may also assume that $\delta<d/4$. 
Lemma \ref{Gronwallgene} can be applied with $\alpha=d/4$ and produces some constants $C=C_{\alpha}$, $D=D_{\alpha}$. Let us now choose 
$$\rho=\rho_{\delta}=\min(\delta/(4D),d/4, \eta/2). 
$$
Application of Lemma \ref{Gronwallgene} gives for any $u\in (s_{\min}-\rho,s_{\min}+\rho) $ (so that $\dist(u,\mathcal{P})>\eta$)
\bal
\Ek(s_{\min}-\rho)\geq C^{-1}e^{-\frac{4}{\e_k}D\rho}\Ek(u)\geq C^{-1}e^{-\frac{\delta}{\e_k}}\Ek(u) .
\nal
 Integrating in $u\in (s_{\min}-\rho,s_{\min}+\rho) $ gives
\baln
\label{inegEk}
\Ek(s_{\min}-\rho)\geq \frac{C}{2\rho }e^{-\frac{\delta}{\e_k}}\int_{s_{\min}-\rho}^{s_{\min}+\rho} \Ek(u)~du \geq C_1e^{-\frac{\delta}{\e_k}} ,
\naln
where we have used $\Ek(u)\geq \frac12 |\varphi_k|^2(u)$ (for $k$ large enough) and Corollary \ref{corpresx0}.

Taking $y\in [s,s+\eta/2]$, we still have $y\leq s_{\min}-\rho$ from the definition of $\rho$. Choosing $\alpha=\rho/2\leq \eta$, we can check that $[y,s_{\min}-\rho]\subset E_{\alpha}$, so that Lemma \ref{GronwallAgmon} applies on $E_{\alpha}$ and gives 
\baln
\label{inegdaeta}
\Ek^+(y)\geq e^{ -\frac{2}{\e_k}\left(|d_A^{\ce}(y)-d_A^{\ce}(s_{\min}-\rho)|+\delta\right)}\Ek^+(s_{\min}-\rho)\geq e^{-\frac{2}{\e_k}\left(d_A^{\ce}(y)+\delta\right)}\Ek^+(s_{\min}-\rho).
\naln
where we have noticed $|d_A^{\ce}(y)-d_A^{\ce}(s_{\min}-\rho)|= d_A^{\ce}(y)-d_A^{\ce}(s_{\min}-\rho) \leq d_A^{\ce}(y)-d_A^{\ce}(s_{\min}) = d_A^{\ce}(y)$.

Since $s_{\min}$ is a strict minimum, there are constants $k_0(\rho,\delta)= k_0(\delta,\eta)>0$ and $C(\rho,\eta) = C(\delta,\eta)>0$ such that for $k\geq k_0(\delta)$,
\bal
V_{\ce}(s_{\min}-\rho)-\mu_k = V_{\ce}(s_{\min}-\rho)-V_{\ce}(s_{\min}) + r(k) \geq C(\delta,\eta)^{-1}.
\nal
This implies
$$
\Ek(s_{\min}-\delta)= \Ek^+ + |\varphi_k|^2 \leq \Ek^+ + C(\delta,\eta)(V_{\ce}-\mu_k) |\varphi_k|^2 \leq \big( 1+C(\delta,\eta) \big) \Ek^+(s_{\min}-\delta) ,
$$
where all functions are taken at the point $(s_{\min}-\delta)$.
Combining this estimate together with~\eqref{inegdaeta} and~\eqref{inegEk} yields
\begin{align*}
\Ek(y) &\geq \Ek^+(y)\geq C e^{-\frac{2}{\e_k}\left(d_A^{\ce}(y)+\delta\right)} \Ek^+(s_{\min}-\rho)
\geq C e^{-\frac{2}{\e_k}\left(d_A^{\ce}(y)+\delta\right)} \Ek (s_{\min}-\rho) \\
& \geq C e^{-\frac{2}{\e_k}\left(d_A^{\ce}(y)+2\delta\right)} 
\geq C e^{-\frac{2}{\e_k}\left(d_A^{\ce}(s)+2\delta\right)},
\end{align*}
where $C$ is a new constant depending only on $\delta$ and $\eta$. Note that in the last inequality, we have used that $d_A^{\ce}(s)\geq d_A^{\ce}(y)$. Integrating for $y\in [s,s+\eta/2]$, we get
\baln
\label{e:intermediate1}
\int_{[s,s+\eta/2]} \Ek(y)dy\geq  Ce^{-\frac{2}{\e_k}\left(d_A^{\ce}(s)+2\delta\right)}.
\naln
Now, turning to the global manifold and recalling that $\psi_k(s,\theta)=e^{ik\theta}\varphi_k(s)$, \eqref{e:reim-gradient} and \eqref{defe}, we have 
\begin{align}
\label{e:intermediate2}
\int_{\CO{[s,s+\eta/2]}} \Ek(y)R(y)dy & = \e_k^2\int_{\CO{[s,s+\eta/2]}} |\nabla_g \psi_k|^2 d\Vol_g \nonumber  \\
 & \quad+ \int_{\CO{[s,s+\eta/2]}}\left(\frac{ |\nabla_g \f|^2}{4}-\mu_k +1\right)|\psi_k|^2 d\Vol_g \nonumber \\
& \leq C \nor{\psi_k}{H^1(\CO{[s,s+\eta/2]})}^2.
\end{align}
Finally, an interpolation estimates together with $P_{\e_k} \psi_k = \mu_k \psi_k$ and the definition of $P_{\e_k}$ in~\eqref{e:defPeps} gives 
\begin{align}
\label{e:intermediate3}
\nor{\psi_k}{H^1(\CO{[s,s+\eta/2]})}^2 & \leq  C_{\eta}\left(\nor{\psi_k}{L^2(\CO{[s,s+\eta/2]})}^2+ \nor{\psi_k}{L^2(\CO{[s-\eta,s+\eta]})}\nor{\Delta_g\psi_k}{L^2(\CO{[s-\eta,s+\eta]})}\right) \nonumber \\
& \leq  C_{\eta}   \e_k^{-2}\nor{\psi_k}{L^2(\CO{[s-\eta,s+\eta]})}^2 .
\end{align}
Now combining~\eqref{e:intermediate1}-\eqref{e:intermediate2}-\eqref{e:intermediate3} gives
\bal
\nor{\psi_k}{L^2(\CO{[s-\eta,s+\eta]})}^2\geq C\e_k^{2} e^{-\frac{2}{\e_k}\left(d_A^{\ce}(s)+2\delta\right)} ,\quad k \geq k_0(\delta, \eta) .
\nal
Finally noticing that $\e_k^{2}\geq e^{-\frac{\delta}{\e_k}}$ for $k\geq k_0(\delta)$ ends the proof of Proposition \ref{propinfdA} up to replacing $\delta$ by $\delta/3$.
\enp

\subsection{Minimal time for uniform controllability in the limit $\eps \to 0^+$}
The main purpose of this section is the proof Theorem \ref{thmtmpsmin} and its corollaries.
\subsubsection{Proof of Theorem \ref{thmtmpsmin}}

Recall that we consider the following situation:
  \begin{itemize}
  \item $\mathcal{S}$ is a surface of revolution as described in Section~\ref{s:revol}, $\ce>0$ is a fixed constant. Moreover, Assumption~\eqref{s:asspt-uniq-min} is fulfilled.
 \item For $k\in \N$, $\e=\e_k=\ce k^{-1}$ and $\psi_k$ is the set of solutions of $P_{\e_k} \psi_k = \mu_k \psi_k$ defined in Lemma \ref{l:exist-mu-psibord} associated to $s_{0}=s_{\min}$ is the minimum of $V_{\ce}$.
 \item the function 
 \baln
 \label{e:def-vk}
 v_k(t,x)=e^{-\frac{\mu_k}{\e_{k}} t}\psi_k(x)
 \naln is the solution to \eqref{e:heat-transp-eps-Witt}, namely $(\e_k \d_t + P_{\eps_k})v_k = 0$, $v_k|_{\d\calS}=0$, and $v_k|_{t=0} = \psi_k$. Here, $\psi_k$ denotes the eigenfunction constructed in Section~\ref{s:revol-eig} above (in particular $\nor{\psi_k}{L^2(\calS)}=1$) and studied in Sections~\ref{s:surf-Agmon}-\ref{sectlower}.
 \end{itemize}
 We now want to test Inequality~\eqref{e:obs-heat-transp-eps-Witt} on $v_k$, and thus estimate both sides of this inequality. This is achieved in Lemmata~\ref{lminfnorexp} and~\ref{lmsupobserv}.
 Theorem~\ref{thmtmpsmin} is then a direct consequence of these two lemmata.
 We recall that $V_\ce$ is defined in~\eqref{e:def-Vce} and $W^{\ce}$ in~\eqref{e:def-Wce}.

\begin{lemma}\label{lminfnorexp}
For any $\delta>0$, there exist $C,k_0>0$ such that for all $k\geq k_{0}$ and all $0\leq T_{0}\leq \delta^{-1}$, we have
\bal
\nor{e^{-\frac{\f}{2\eps}}v_{k}(T_0)}{L^2(\calS)} \geq C e^{-\frac{V_{c}(s_{\min})T_{0}+W^{\ce}_m+\delta }{\e_k}}, \quad W^{\ce}_m = \min_{I_L}W^{\ce} ,
\nal
with $v_k$ defined in~\eqref{e:def-vk} and $W^\ce(s) = d_A^{\ce}(s) +\frac{\f(s)}{2}$.
\end{lemma}
\bnp
Note first that the function $W^{\ce} = d_A^{\ce} +\frac{\f}{2}$ is continuous on $I_L$ and converges to $+\infty$ close to the poles $p\in \mathcal{P}$, according to the asymptotics of $d_A^{\ce}$ in Lemma~\ref{lemma-prop-dA}. Hence, it reaches its minimum in (at least) one point that we denote $s_1\in I_L$, that is $W^{\ce}(s_1) = W^{\ce}_m=\min_{I_L}W^{\ce}$.
We take $0<\eta<\dist(s_{1},\mathcal{P})$ small enough so that $|\f(s)-\f(s_1)|\leq \delta$ for $|s-s_{1}|\leq \eta$. We have for $k\geq k_{0}$ large enough
\begin{align*}
\nor{e^{-\frac{\f}{2\eps}}v_{k}( T_0)}{L^2(\calS)} &=e^{-\frac{\mu_k}{\e_{k}} T_{0}}\nor{e^{-\frac{\f}{2\eps}}\psi_k}{L^2(\calS)}\geq e^{-\frac{\mu_k}{\e_{k}} T_{0}} \nor{e^{-\frac{\f}{2\eps}}\psi_k}{L^2(\CO{(s_{1}-\eta, s_1 +\eta)})}\\
&\geq e^{-\frac{\mu_k}{\e_{k}} T_{0}} e^{-\frac{\f(s_1)+\delta }{2\e_{k}} }\nor{\psi_k}{L^2(\CO{(s_{1}-\eta,s_1+\eta)})} \\
& \geq Ce^{-\frac{\mu_k}{\e_{k}} T_{0}} e^{-\frac{\f(s_1)+\delta }{2\e_{k}} } e^{-\frac{1}{\e_k}(d_A^{\ce}(s_{1})+\delta)} ,
\end{align*}
where we have used Proposition \ref{propinfdA} for the last estimate. 
Since $\mu_{k}\to V_{c}(s_{\min})$, we have $\mu_{k}T_{0}\leq V_{c}(s_{\min})T_{0}+\delta$ for $k$ large enough, which, together with the above estimate, concludes the proof (up to changing $3\delta$ into $\delta$).
\enp
\begin{lemma}
\label{lmsupobserv}
For any $\omega \subset \calS$ and $\delta>0$, there exist $C,k_0>0$ such that for all $k\geq k_{0}$,  and all $0\leq T_{0}\leq \delta^{-1}$, we have
\bal
\int_0^{T_0}\nor{e^{-\frac{\f}{2\eps}}v_k(t,\cdot)}{L^2(\omega)}^2 dt
& \leq C e^{-2\frac{W^{\ce}_{\omega} - \delta }{\e_k}} , \quad W^{\ce}_\omega = \min_{\bar{\omega}}W^{\ce} , \\
\int_0^{T_0}\nor{e^{-\frac{\f}{2\eps}} \d_s v_k(t)|_{s=L} }{H^1(\S^1)}^2 dt & \leq C e^{-2\frac{W^{\ce}(L) - \delta }{\e_k}}, \quad \text{in Cases~\ref{caseD} and~\ref{caseCY}}, \\
\int_0^{T_0}\nor{e^{-\frac{\f}{2\eps}} \d_s v_k(t)|_{s=0} }{H^1(\S^1)}^2 dt & \leq C e^{-2\frac{W^{\ce}(L) - \delta }{\e_k}}, \quad \text{in Case~\ref{caseCY}},
\nal
with $v_k$ defined in~\eqref{e:def-vk}.
\end{lemma}
\bnp
Since we need an upper bound, we can assume without loss of generality that $\omega$ is invariant by rotation. Also, $W^{\ce}_\omega$ is finite except in the trivial case $\omega\subset \mathcal{P}$.
Let $\delta>0$. 

We first estimate the contribution close to the poles in case $\mathcal{P}\neq \emptyset$. There, the function $\psi_k$ (hence $v_k$) is supposed to be very small since $d_A^{\ce}$ is large. More precisely, using the asymptotics of $d_A^{\ce}$ close to $\mathcal{P}$ given by Lemma~\ref{lemma-prop-dA}, there exists $\widetilde{\eta}>0$ so that 
\baln
\label{infdAP}
\frac{d_A^{\ce}(s)}{2}\geq W^{\ce}_{\omega}+\frac{\nor{\f}{L^{\infty}}}{2} +1;\quad \text{ for all } s\in  \calN_{\tilde{\eta}}=\{ s\in  I_L \textnormal{ such that } \dist(s,\mathcal{P}) < \widetilde{\eta}\} .
\naln
We start with the estimate
\baln
\label{e:interm-noref}
\int_0^{T_0}\nor{e^{-\frac{\f}{2\eps}}v_{k} }{L^2(\omega\cap\CO{\calN_{\tilde{\eta}}})}^2 dt\leq e^{\frac{\nor{\f}{L^{\infty}}}{\eps}} \nor{\psi_k }{L^2(\CO{\calN_{\tilde{\eta}}})}^2\int_0^{T_0}e^{-2\frac{\mu_k}{\e_{k}} t} dt .
\naln
 We have $\int_0^{T_0}e^{-2\frac{\mu_k}{\e_{k}} t} dt=\frac{\e_{k}}{2\mu_k}\int_0^{ T_0\frac{2\mu_k}{\e_{k}}}e^{- s}  ds\leq C\frac{\e_{k}}{2\mu_k}\leq 1$ for $k$ large. Corollary \ref{coragmonfaible} applied to the set $\CO{\calN_{\tilde{\eta}}}$ and the constant $\delta=1/2$, together with~\eqref{infdAP} implies 
 \bal
 \nor{\psi_k }{L^2(\CO{\calN_{\tilde{\eta}}})}\leq C e^{-\frac{1}{\e_k}  (W^{\ce}_{\omega} +\frac{\nor{\f}{L^{\infty}}}{2})}  .
 \nal
 With~\eqref{e:interm-noref}, this gives
 $$
 \int_0^{T_0}\nor{e^{-\frac{\f}{2\eps}}v_{k} }{L^2(\omega\cap \CO{\calN_{\tilde{\eta}}})}^2 dt
 \leq C e^{\frac{\nor{\f}{L^{\infty}}}{\eps}}e^{-\frac{2}{\e_k}  (W^{\ce}_{\omega} +\frac{\nor{\f}{L^{\infty}}}{2})}  = Ce^{-\frac{2}{\e_k}  W^{\ce}_{\omega}},
 $$
which is the expected bound for this part.

 Let us now treat the contribution of the norm away from the poles (which is the whole $I_L$ in Case~\ref{caseCY}).  
Since $\f$ is uniformly continuous on $[0,L]$ and $d_A^{\ce}$ is uniformly continuous on the compact set $\calN_{\tilde{\eta}}^{\mathsf{c}} :=I_L \setminus \calN_{\tilde{\eta}}$, there exists $\eta>0$ so that 
\begin{align}
\label{e:unif-cont}
s,s'\in \calN_{\tilde{\eta}}^{\mathsf{c}} \  \text{ and } \ |s-s'|\leq \eta \quad \Longrightarrow \quad  |\f(s)-\f(s')|\leq \delta \text{ and } |d_A^{\ce}(s)-d_A^{\ce}(s')|\leq \delta .
\end{align}
We now select a finite sequence $s_i\in  \calN_{\tilde{\eta}}^{\mathsf{c}}, i=1,\cdots,N$ so that $\left\{|s-s_i|\leq \eta\right\}$ is a finite covering of $\calN_{\tilde{\eta}}^{\mathsf{c}}$. This property gives the estimate
\bal
\int_0^{T_0}\nor{e^{-\frac{\f}{2\eps}}v_{k} }{L^2(\omega\cap \CO{\calN_{\tilde{\eta}}^{\mathsf{c}}})}^2 dt\leq \int_0^{T_0}e^{-2\frac{\mu_k}{\e_{k}} t} dt\sum_{i\in \mathcal{J}} \nor{e^{-\frac{\f}{2\eps}}\psi_k }{L^2(\CO{(s_i-\eta, s_i +\eta)})}^2
\nal
where $\mathcal{J}=\left\{i=1, \cdots ,N; \ (s_i-\eta, s_i +\eta)\cap\omega \neq \emptyset\right\}$.
Using $|f(s)-f(s_{i})|\leq \delta$ and $|d_A^{\ce}(s)-d_A^{\ce}(s_{i})|\leq \delta$ for $s\in (s_i-\eta, s_i +\eta)$ and then Corollary \ref{coragmonfaible} with some $0<\delta'\leq \delta/(1+d_A^{\ce}(s_{i}))$ so that $(1-\delta')d_A^{\ce}(s_{i}) - \delta' \geq d_A^{\ce}(s_{i}) - \delta $, we obtain
\begin{align*}
\nor{e^{-\frac{\f}{2\eps}}\psi_k }{L^2(\CO{(s_i-\eta, s_i +\eta)})}& \leq e^{-\frac{f(s_i)+\delta}{2\eps}} \nor{\psi_k }{L^2(\CO{(s_i-\eta, s_i +\eta)})}\\
& \leq C e^{-\frac{f(s_i)+\delta}{2\eps}}e^{-\frac{1}{\e_k} ((1-\delta')d_A^{\ce}(s_{i}) - \delta') } e^{\frac{\delta}{\eps}} 
\leq C e^{-\frac{W^{\ce}(s_i)-4\delta}{\e_k} } .
\end{align*}
We finally obtain 
\bal
\int_0^{T_0}\nor{e^{-\frac{\f}{2\eps}}v_{k} }{L^2(\omega\cap \CO{\calN_{\tilde{\eta}}^{\mathsf{c}}})}^2 dt\leq  C e^{\frac{8\delta}{\e_k} }\max_{i\in \mathcal{J}}e^{-2\frac{W^{\ce}(s_i)}{\e_k} }. 
\nal
We remark from~\eqref{e:unif-cont} that $\min_{i\in \mathcal{J}}W^{\ce}(s_i)\geq W^{\ce}_{\omega}-2\delta $. This finishes the proof since $\delta$ is arbitrary and $k$ can be chosen large enough.

Finally, the proof of the boundary estimates simply consists in replacing the use of Corollary~\ref{coragmonfaible} by that of the Inequality~\eqref{e:estimate-psik-bord}.
\enp

We may now conclude the proof of Theorem~\ref{thmtmpsmin} from Lemmata~\ref{lminfnorexp} and~\ref{lmsupobserv}
\bnp[Proof of Theorem~\ref{thmtmpsmin}]
Using Lemma \ref{lemequiveunif}, if uniform observability holds for $T_{0}$, then, we have the inequality
\bal
\nor{e^{-\frac{\f}{2\eps}}v(T_0)}{L^2(\calS)}^2 \leq C_0^2\int_0^{T_0}\nor{e^{-\frac{\f}{2\eps}}v }{L^2(\omega)}^2 dt .
\nal
for any solution $v$ of \eqref{e:heat-transp-eps-Witt}. In particular, this inequality holds true for the sequence $v_{k}$ described above. So, combining Lemmata~\ref{lminfnorexp} and~\ref{lmsupobserv}, we obtain
$$
e^{-\frac{V_{c}(s_{\min})T_{0}+W^{\ce}_m+\delta }{\e_k}}\leq C_{0}e^{-\frac{W^{\ce}_{\omega} - \delta }{\e_k}  }, \quad k\geq k_{0}(\delta), \quad \text{ for all } k\geq k_0 .
$$
This implies $V_{c}(s_{\min})T_0\geq W^{\ce}_{\omega}-W^{\ce}_{m}-2\delta$ when letting $k\to +\infty$, which gives the expected result since $\delta$ is arbitrary. 
The proof of the boundary observability estimate~\eqref{e:estiminf-Tunif-bord} follows the same. We also notice that all the previous results also apply in the case with the alternative Definition~\ref{e:defbisdAT2} of the Agmon distance which is also Lipschitz with the same properties that we used.
\enp
\subsubsection{Proof of Corollaries of Theorem~\ref{thmtmpsmin}}

Corollaries~\ref{c:sphere-disk},~\ref{c:example-torus} and~\ref{c:example-cylinder}, stated in the introduction, are significant examples of application of Theorem~\ref{thmtmpsmin}.
In this section, we prove these three results.

\bnp[Proof of Corollary~\ref{c:sphere-disk}]
We first consider the case $\calS$ diffeomorphic to $\S^2$, i.e. Case~\ref{caseS}. The case $\calS$ diffeomorphic to $\ID$, i.e. Case~\ref{caseD}, is discussed at the end of the proof.
We define $f_\delta (s) = \int_0^s \chi_\delta (t) dt$ where $\chi_\delta \in C^\infty_c((0,L);[0,1])$, $\chi_\delta(s) =1$ in a neighborhood of $[\delta , L-\delta]$. Such a function $f_\delta$ is constant near $0$ and $L$, and hence can be extended by continuity as a $C^\infty$ function on $\calS$ (see e.g.~\cite[Proposition~4.6]{Besse}). We notice that $f_\delta' (s)=1$ for $s\in [\delta,L-\delta]$, so the statements about $T_{GCC}$ are direct consequences of Proposition \ref{p:GCCFCrot}.

Notice now that we have
$$
V_\ce (s) := \frac{\ce^2}{R(s)^2} + \frac{|\chi_\delta(s)|^2}{4} .
$$
Let us call $s_{\min} \in (0,L)$ the unique point such that $R(s_{\min}) = \max R$, that is $\frac{1}{R(s_{\min})^2}= \min \frac{1}{R^2}$.

\paragraph{Claim:} For all $\ce>0$, there is $\delta_0>0$ such that for all $\delta \in (0,\delta_0)$, we have $V_\ce(s_{\min})= \min V_\ce$ and $V_\ce^{-1}(V_\ce(s_{\min})) = \{ s_{\min}\}$.

 To prove the claim, we let $\delta_0>0$ be such that $$
s_{\min} \in [\delta_0 , L-\delta_0] , \quad \text{ and }
 \quad  \frac{\ce^2}{R^2(s)} > \frac{\ce^2}{R(s_{\min})^2} +\frac14 \text{ for } s \notin [\delta_0 , L-\delta_0]
$$
(note that $\delta_0$ thus depends on $\ce$).
This is possible since $R(s)\to 0$ as $s\to 0^+$ and $s\to L^-$.
Hence, recalling the definitions of $\chi_\delta$ and $V_\ce$, for $\delta < \delta_0$ we have 
$$
V_\ce(s_{\min})=  \frac{\ce^2}{R(s_{\min})^2} + \frac14 , \quad V_\ce(s) > \frac{\ce^2}{R(s_{\min})^2} + \frac14  \text{ for all }s \in [\delta_0 , L-\delta_0] \setminus \{s_{\min}\} ,
$$
together with 
$$
V_\ce(s) \geq  \frac{\ce^2}{R^2(s)} > \frac{\ce^2}{R(s_{\min})^2} + \frac14  \text{ for }s \notin [\delta_0 , L-\delta_0] .
$$ 
As a consequence, $V_\ce$ reaches its minimum at $s_{\min}$ only, which proves the claim.

Now Assumption~\ref{a:single-min} is satisfied and $s_{\min}$ does not depend on $\delta$. We compute the Agmon distance~\eqref{e:defbisdA} 
\begin{align*}
d_A^{\ce}(s) & = \left| \int_{s_{\min}}^s \sqrt{ V_\ce(y)- V_\ce(s_{\min})} dy \right| \\
&=  \left| \int_{s_{\min}}^s \sqrt{ \frac{\ce^2}{R(y)^2} + \frac{|\chi_\delta(y)|^2}{4}- \frac{\ce^2}{R(s_{\min})^2} - \frac14} dy \right| .
\end{align*}
Notice that $0\geq \frac{|\chi_\delta(y)|^2}{4}- \frac14 \geq - \frac14$ uniformly with respect to $\delta$ so that the asymptotic expansion of Lemma~\ref{lemma-prop-dA} is valid uniformly in $\delta$: there are $C,\gamma>0$ such that for all $\delta \in (0,\delta_0)$, 
\begin{align*}
|d_A^{\ce}(s) + \ce \log(s)| \leq C  \quad \text{for } s \in (0,  \gamma ] ,  \\  
|d_A^{\ce}(s) + \ce \log (L-s) | \leq C, \quad \text{for } s \in [L-\gamma, L) .  
\end{align*}
Now we recall the definition of $W^{\ce} = d_A^{\ce} +\frac{\f_{\delta}}{2}$, notice that $0\leq \f'_{\delta}\leq 1$ so that $0\leq \f_{\delta}(s)\leq L$ for $s \in I_L$. As a consequence, using that $d_A^{\ce}\geq0$ and $d_A^{\ce}(s_{\min}) = 0$, we have for $\delta \in (0,  \gamma ]$,
\begin{align*}
W^{\ce}_m & = \min_{I_L}W^{\ce}  = \min_{I_L} \left( d_A^{\ce} +\frac{\f_{\delta}}{2} \right) \leq \frac{L}{2} + \min_{I_L} d_A^{\ce} \leq \frac{L}{2}  \\
W^{\ce}_{\omega} & = \min_{\ovl{\omega}}W^{\ce} \geq  \min_{\ovl{\omega}} d_A^{\ce} \geq  - \ce \log(\delta) - C \quad \text{ for }\omega = \ B_g(N, \delta) \cup B_g(S, \delta) .
\end{align*}
The bound~\eqref{e:estiminf-Tunif} of Theorem~\ref{thmtmpsmin} then yields
$$
\left( \frac{\ce^2}{R(s_{\min})^2} + \frac14 \right) T_{unif}(\omega)\geq  - \ce \log(\delta) - C -\frac{L}{2} ,
$$
and hence concludes the proof in the case $\calS$ diffeomorphic to $\S^2$, i.e. Case~\ref{caseS}. 

In the case $\calS$ diffeomorphic to $\ID$, i.e. Case~\ref{caseD}, we instead define $f_\delta (s) = \int_0^s \chi_\delta (t) dt$ where $\chi_\delta \in C^\infty_c((0,L];[0,1])$, $\chi_\delta(s) =1$ in a neighborhood of $[\delta , L]$. Then, the remainder of the proof is the same except that the minimum can be achieved at $s=L$, and all sets of the form $[\delta_0 , L-\delta_0]$ have to be replaced by $[\delta_0 , L]$ (i.e., only a neighborhood of zero is avoided, and not a neighborhood of $L$).
\enp

\bnp[Proof of Corollary~\ref{c:example-torus}]
For the sake of simplicity, we may identify $\S^1_L = [-L/2,L/2]$ and $I_\omega= (-\alpha,\alpha)$ for $\alpha \in (0,L/2)$. 
We first choose $\chi_\delta \in C^\infty(\S^1_L)$ even in this identification, and $s_{\min} = \pm L/2 \notin (-\alpha,\alpha)$ such that 
\begin{align*}
&\chi_\delta =\frac{1}{\delta} \text{ on } \left[-\frac12\big(\alpha+\frac{L}{2}\big),\frac12\big(\alpha+\frac{L}{2}\big) \right],  \quad  \chi_\delta(s_{\min}) = 1 , \\ 
& \chi_\delta'(s)  = 0 \text{ if and only if } s\in \left[-\frac12\big(\alpha+\frac{L}{2}\big),\frac12\big(\alpha+\frac{L}{2}\big) \right] \text{ or } s= s_{\min} . 
\end{align*}
Note that $\chi_\delta=\frac{1}{\delta}$ in a neighborhood of $(-\alpha, \alpha)$, and that for $\delta<1$, $\chi_\delta$ reaches at $s_{\min}=\pm L/2$ a unique global minimum and in particular $\chi_\delta \geq 1$ on  $\S^1_L$.
We then set
$$
V^\delta(s) := \chi_\delta +M, \quad M:= \max_{\S^1_L}\frac{|\f'|^2}{4} ,
$$
so that $V^\delta(s)>\frac{|\f'(s)|^2}{4}$ on $\S^1_L$, and 
$$
 R_\delta(s) := \left(V^\delta(s) - \frac{|\f'(s)|^2}{4} \right)^{-\frac12}  , \quad R_\delta \in C^\infty(\S^1_L; \R^+_*) ,
$$
where we have used that $\f\in C^{\infty}(\S^1_L)$.
Notice that with these definitions, we have (where $V_1$ denotes $V_\ce$ with $\ce=1$)
$$
V_1 (s) = \frac{1}{R_\delta(s)^2} + \frac{|\f'(s)|^2}{4} = V^\delta(s) .
$$
Notice then that $V^\delta$ admits a unique global minimum at $s_{\min} =\pm L/2$. 

On the one hand, we have 
$$
(\chi_\delta(s) +M)^{-\frac12}\leq R_\delta(s) \leq \chi_\delta^{-\frac12}(s).
$$
On the other hand, since $\chi_\delta$ is even and with the appropriate definition \eqref{e:defbisdAT2}, we have for $s\in [0, L/2]$, 
$$
d_A^{\ce}(s) = - \int_{L/2}^s \sqrt{ V^\delta(y)- V^\delta(L/2)} dy 
 =  - \int_{L/2}^s \sqrt{\chi_\delta(y) - 1} dy .
$$
As a consequence, we have (recall $\omega = (-\alpha ,\alpha) \times \S^1$)
\begin{align*}
W^{\ce}_m & = \min_{[0,L]}W^{\ce}  = \min_{[0,L]} \left( d_A^{\ce} +\frac{\f}{2} \right) \leq \left( d_A^{\ce} +\frac{\f}{2} \right) (L/2) = \frac{\f(L/2)}{2} , \\
W^{\ce}_{\omega} & \geq \min_{\overline{\omega}} d_{A}^{\ce} +\min_{\overline{\omega}} \frac{\f}{2}  = d_{A}^{\ce}(\alpha) +\min_{\overline{\omega}} \frac{\f}{2},
\end{align*}
where
\begin{align*}
d_A^{\ce} (\alpha) &=\int_\alpha^{L/2} \sqrt{\chi_\delta(y) - 1} dy \geq \int_\alpha^{\frac12(\alpha+\frac{L}{2})} \sqrt{\chi_\delta(y) - 1} dy
 =\int_\alpha^{\frac12(\alpha+\frac{L}{2})} \sqrt{\frac{1}{\delta} - 1} dy \\
 & \geq \frac12(\frac{L}{2}-\alpha) (\delta^{-1/2}-1)
\end{align*}
for $\delta<1$. Applying Theorem~\ref{thmtmpsmin} concludes the proof of the corollary.
\enp

We next prove Corollary~\ref{c:example-cylinder}, stated in the introduction. The proof is close to that of Corollary~\ref{c:example-torus}.

\bnp[Proof of Corollary~\ref{c:example-cylinder}]
We first take $\chi \in C^\infty(\R_+;[0,1])$ such that $\chi = 1$ on $[0, L/4]$ and $\supp\chi \subset [0,L/2)$. Next, with 
$$
M:= \max_{[0,L]}\frac{|\f'|^2}{4} ,  
$$ we set,
\bal
&V^\delta(s) := \frac{\chi(s)}{(s+\delta)^\gamma} + \left( 1-\chi(s) \right) \left(s-\frac{L}{2} \right)^2+M, \quad  \text{ for }s \in [0,\frac{L}{2}] ,\\
&V^\delta(s) := V^\delta(L - s) \quad \text{ for }s \in [\frac{L}{2}, L] . 
\nal
This function is symmetric about $L/2$, smooth on $[0,L]$, and satisfies $V^\delta(s)>\frac{|\f'(s)|^2}{4}$ on $[0,L]$. Hence, defining  $R_\delta$ as 
$$
 R_\delta(s) := \left(V^\delta(s) - \frac{|\f'(s)|^2}{4} \right)^{-\frac12}  ,
$$
and using that $\f\in C^{\infty}([0,L])$, we deduce that $R_\delta \in C^\infty([0,L]; \R^+_*)$.
Notice that with these definitions, we have (where $V_1$ denotes $V_\ce$ with $\ce=1$)
$$
V_1 (s) = \frac{1}{R_\delta(s)^2} + \frac{|\f'(s)|^2}{4} = V^\delta(s) .
$$
Notice then that $V^\delta$ admits a unique global minimum at $s_{\min} = L/2$. 

We have on the one hand that for all $s\in [0,L/4]$ and $\delta \in (0,\delta_0]$, 
$$
(s+\delta)^{\frac{\gamma}{2}} \frac{1}{\sqrt{1+M(s+\delta)^\gamma}}= \left(V^\delta(s)  \right)^{-\frac12} 
\leq  R_\delta(s) \leq \left(V^\delta(s) -M \right)^{-\frac12} = (s+\delta)^{\frac{\gamma}{2}} ,
$$
which proves the Item \ref{itemdelta} (on account to the symmetry $V^\delta(s) := V^\delta(L - s)$). 

On the other hand, recalling that $s_{\min}=L/2$, we have for $s\in [0,L/2]$
\bal
V_1(s_{\min}) & = V^\delta(L/2)  = M ,\\  
d_A^{\ce}(s) & = \left| \int_{L/2}^s \sqrt{ V^\delta(y)- V^\delta(s_{\min})} dy \right| \\
&  =  \left| \int_{L/2}^s \sqrt{\frac{\chi(s)}{(s+\delta)^\gamma} + \left( 1-\chi(s) \right) \left(s-\frac{L}{2} \right)^2} dy \right| .
\nal
As a consequence, we have 
\begin{align*}
W^{\ce}_m & = \min_{[0,L]}W^{\ce}  = \min_{[0,L]} \left( d_A^{\ce} +\frac{\f}{2} \right) \leq \left( d_A^{\ce} +\frac{\f}{2} \right) (L/2) = \frac{\f(L/2)}{2} , \\
W^{\ce}_{\{0\}} & = W^{\ce}(0) \geq  d_A^{\ce} (0) +\frac{\f(0)}{2} ,
\end{align*}
where
\begin{align*}
d_A^{\ce} (0) &=   \int_0^{L/2} \sqrt{\frac{\chi(s)}{(s+\delta)^\gamma} + \left( 1-\chi(s) \right) \left(s-\frac{L}{2} \right)^2}dy 
\geq  \int_0^{L/4} \sqrt{\frac{1}{(s+\delta)^\gamma}}dy \\
&=  \frac{\delta^{1-\gamma/2}}{\gamma/2-1} - \frac{(\delta+L/4)^{1-\gamma/2}}{\gamma/2-1}  
=  \frac{\delta^{1-\gamma/2}}{\gamma/2-1} + \grando{1} , \quad \text{for } \gamma >2 .
\end{align*}
By symmetry of $V^\delta$ about $L/2$, we also have $d_A^{\ce} (L) = d_A^{\ce} (0)$.
Applying Theorem~\ref{thmtmpsmin} concludes the proof of the corollary.
\enp

\section{Uniform time of observability for positive solutions}
\label{s:positive}

The proofs of Theorem \ref{t:positive-intro} and Proposition~\ref{proplowerobse+} rely on fine estimates on the semiclassical heat kernel, which we borrow from~\cite{LY:86}. The latter are first presented in Section~\ref{s:Li-Yau}. Then, in Section~\ref{s:L1obs}, we deduce $L^1$ observability statements and finally conclude the proofs of Theorem \ref{t:positive-intro} and Proposition~\ref{proplowerobse+} in Section~\ref{s:L1toL2}.
Throughout this section, we assume $\d\M=\emptyset$.

\subsection{Estimates on the semiclassical heat kernel}
\label{s:Li-Yau}
The main tool we use to estimate the heat kernel in the semiclassical limit is the following theorem taken from Li-Yau~\cite{LY:86} (see also \cite{Simon:83}, for a similar result on $\R^{n}$).
\begin{theorem}[Theorem 6.1 of \cite{LY:86}]
\label{t:thm-LiYau}
Let $\M$ be a compact manifold without boundary. Suppose $V_\eps = V + \eps q$ with $V,q \in C^2(\M)$. For any $\e>0$, we consider $H_{\e}$, the fundamental solution of 
\bal
\partial_t w-\Delta_g w+\frac{1}{\e^2}V_\eps (x)w=0 , \quad \text{ on }(0,+\infty)\times \M .
\nal
Then, we have 
\baln
\label{e:Li-Yau}
\underset{\e\to 0}{\lim}\e \log H_{\e}(x,y,\e t)=-\rho(x,y,t)
\naln
where
\baln
\label{e:def-rho-V}
\rho(x,y,t)=  \inf \left\{ \int_0^t \frac{1}{4}|\dot{\gamma}(s)|_g^2 +V \big( \gamma(s) \big)ds , \gamma \in W^{1,\infty} ([0,t] ; \M ), \gamma(0)=x , \gamma(t)= y \right\} .
\naln 
Moreover, the limit in~\eqref{e:Li-Yau} is uniform on any compact set of $\M^2\times (0,+\infty)$. 
\end{theorem}
We recall that $H_\eps(x,y,t)$ is defined to be the unique solution to 
  \begin{equation}
  \label{e:heat-kernel-def}
 \left\{
 \begin{aligned}
&\Big(\d_t -\Delta_g+\frac{1}{\e^2} V_\eps(x) \Big) H_{\e}(x,y,t)= 0 , & \text{ for } (t,x) \in \R^+_*\times \M ,  \\
& H_{\e}(x,y,t) |_{t=0} = \delta_{x=y}, & \text{ for } x \in  \M ,
 \end{aligned}
 \right.
 \end{equation}
 where $y \in \M$ is fixed, and the differential operator $ -\Delta_g+\frac{1}{\e^2} V_\eps(x)$ acts in the $x$-variable.
We also recall that $H_{\e}(x,y,t)$ is (well-defined and) continuous in $\M^2\times (0,+\infty)$ as soon as $V_\eps\in L^\infty(\M)$, see e.g.~\cite[Theorem~B.7.1~(a$'''$)]{Simon:82}, so that pointwise estimates like~\eqref{e:Li-Yau} make sense.

\medskip
The statement of Theorem~\ref{t:thm-LiYau} is not strictly speaking a consequence of~\cite[Theorem 6.1]{LY:86} for the following two reasons: 
\begin{itemize}
\item the potential $V_\eps$ is assumed independent of $\eps$ in~\cite{LY:86}; 
\item the uniformity of the limit on any compact subset is not explicitly written in~\cite{LY:86}. 
\end{itemize}
However, let us explain why the proof of~\cite[Theorem 6.1]{LY:86}  actually contains these two points.
The limit~\eqref{e:Li-Yau} is proved in two steps, a lower bound and an upper bound.
The lower bound
$$
\underset{\e\to 0}{\lim}\e \log H_{\e}(x,y,\e t) \geq -\rho(x,y,t)
$$ is obtained as a consequence of the explicit estimate (\cite[Theorem 2.1]{LY:86})
$$
\e \log H_{\e}(x,y,\e t_1) \leq \e \log H_{\e}(x,y,\e t)  
+ \eps^2 A_{1/\eps} (t-t_1) + \rho_{\alpha,R}(x,y,t-t_1) , 
$$
for all $0<t_1<t$, $\alpha>1$ and
where the constant $A_{\lambda}$ ($\lambda=\e^{-1}$ in our context) only depends on $\nor{\Delta_g V_\eps}{L^\infty(\M)} , \nor{\nabla_g V_\eps}{L^\infty(\M)}$ and $\alpha$. In particular, the limit $\lim_{\eps\to 0} \eps A_{1/\eps} = 0$ holds uniformly.
The proof finally proceeds by taking the limit $R\to + \infty$, and then $\alpha \to 1^+$ (in which case $\rho^{V_\eps}_{\alpha,R}(x,y,t-t_1)) \to \rho^{V_\eps}(x,y,t-t_1)$, uniformly on compact sets; here $\rho^{V_\eps}$ denotes the function $\rho$ defined by~\eqref{e:def-rho-V} with $V$ replaced by $V_\eps$). Taking finally the limit $\eps \to 0^+$, and then $t_1 \to 0^+$, and noticing that
$$\lim_{t_1\to 0^+}\lim_{\eps \to 0^+} \e \log H_{\e}(x,y,\e t_1) \geq 0$$ concludes the proof of the lower bound.
This last argument relies only on a uniform upper bound on $\nor{V_\eps}{L^\infty(\M)}$ together with a comparison argument, namely Inequality \cite[(6.3)]{LY:86}. It can be checked that the convergence is uniform on any compact set since it involves the asymptotics on the diagonal of kernels of heat equations on large balls and without potential, which are known to be uniform on compact sets, see for instance \cite[Theorem 4.6]{Va:67b}.
It only remains to notice that the limit $\lim_{t_1\to 0^+} \lim_{\eps \to 0^+} \rho^{V_\eps}(x,y,t-t_1) = \rho(x,y, t)$ is uniform on compact sets.

The upper bound 
$$
\underset{\e\to 0}{\lim}\e \log H_{\e}(x,y,\e t) \leq -\rho(x,y,t)
$$
follows from \cite[Theorem 3.3]{LY:86}. As for the lower bound, this result also furnishes an explicit and uniform bound involving another constant $A_{\lambda}$ enjoying the same type of convergence properties as for the lower bound. 

\medskip
Note that we have chosen to use the estimate of the semiclassical limit of the Kernel of~\cite{LY:86} but it is likely that we could have obtained  the observability inequality starting directly from Harnack type inequalities like~\cite[Theorem 2.1]{LY:86}, as we did in \cite{LL:18}.

Here, we are mostly interested in the case $V(x) = \frac14 |\nabla_g \f(x)|_g^2$. In this situation, we can reformulate the result in terms of the transport equation with vanishing viscosity~\eqref{e:transport-viscous-gradient}.

We next define $K_{\e}(x,y,t)$, the fundamental solution of~\eqref{e:transport-viscous-gradient} on $\M^2\times(0,+\infty)$ by the unique solution to 
  \begin{equation}
  \label{e:transport-viscous-kernel}
 \left\{
 \begin{aligned}
&(\d_t - \nablag \f \cdot \nablag -\q - \eps \Delta_g ) K_{\e}(x,y,t)= 0 , & \text{ for } (t,x) \in \R^+_*\times \M ,  \\
& K_{\e}(x,y,t) |_{t=0} = \delta_{x=y}, & \text{ for } x \in  \M ,
 \end{aligned}
 \right.
 \end{equation}
 where $y \in \M$ is fixed, and the differential operator $- \nablag \f \cdot \nablag -\q - \eps \Delta_g$ acts in the $x$-variable.
 Recall that for $u_0 \in L^1(\M)$, the function
 $$
 u(t,x) = \int_\M K_\eps(x,y,t)u_0(y) d\Vol_g(y)
 $$
 is the unique solution of~\eqref{e:transport-viscous-gradient} on $(0,+\infty)\times \M$ issued from $u|_{t=0} = u_0$ (This uses the choice of volume form $d\Vol_g$ in the embedding $L^1(\M) \hookrightarrow \mathcal{D}'(\M)$).

\begin{corollary}
\label{corLiYauconj}
Let $\M$ be a compact manifold without boundary. Suppose $X=\nabla_g \f$ where $\f$ is a $C^3$ function defined on $\M$. For any $\e>0$, we consider $K_{\e}$, the fundamental solution of~\eqref{e:transport-viscous-gradient} on $(0,+\infty)\times \M$. Then, we have
\bal
\underset{\e\to 0}{\lim}\e \log K_{\e}(x,y, t)=-d_{\nablag \f}(x,y,t) ,
\nal
with  
\begin{align}
\label{defdX}
d_{\nablag \f}(x,y,t)& := \rho(x,y,t)+ \frac{\f(x)-\f(y)}{2} ,
\end{align}
where $\rho(x,y,t)$ is defined by~\eqref{e:def-rho-V} with $V(x) = \frac{|\nabla_g \f|_g^2}{4}$.
Moreover, the limit is uniform on any compact set of $\M^2 \times (0,+\infty)$.

In particular, for any $\delta>0$ and any compact subset $I \Subset (0,+\infty)$, there exists $\e_0>0$ such that   
\baln
\label{inegkernel}
e^{-\frac{d_{\nablag \f}(x,y,t)+\delta}{\e}}\leq K_{\e}(x,y,t)\leq e^{-\frac{d_{\nablag \f}(x,y,t)-\delta}{\e}}
\naln 
for any $(x,y,t)\in \M^2 \times I$ and $0<\e\leq \e_0$.
\end{corollary}
Note that the definition of $d_{\nablag \f}$ in~\eqref{defdX}  is not the same as that given in~\eqref{defdXintro} in the introduction. Equivalence between these two definitions is proved in Lemma~\ref{lmdist}.

Note that although the kernel $H_\e(x,y,t)$ is symmetric (with respect to the Riemannian volume measure $d\Vol_g$) since $-\Delta_g+\frac{1}{\e^2} V_\eps$ is, this is no longer the case for the kernel $K_\e(x,y,t)$ (since the operator in~\eqref{e:transport-viscous-gradient} is not symmetric in $L^2(\M,d\Vol_g)$). Similarly, $\rho(x,y,t)$ is symmetric whereas $d_{\nablag \f}(x,y,t)$ is not.

\bnp
Setting $u_y(t,x) = K_{\e}(x,y, t)$ and $w_y(t,x) = e^{\f(x)/2\eps} u_y(t/\eps,x)$, we have from Equation~\eqref{e:transport-viscous-kernel} and Lemma~\ref{lemequiveunif} that $w_y(t,x)$ solves~\eqref{e:heat-transp-eps-Witt-eps}. Moreover, we have $$w_y(0,x)= e^{\f(x)/2\eps} u_y(0,x)= e^{\f(x)/2\eps} \delta_{x=y} = e^{\f(y)/2\eps}\delta_{x=y}.$$
This implies that $w_y (t,x) = e^{\f(y)/2\eps} H_\eps (x,y,t)$ where $H_\eps (x,y,t)$ is defined in~\eqref{e:heat-kernel-def} with $V_\eps = \frac{|\nablag \f|_g^2}{4} + \eps \left( \frac{\Lap \f}{2}- \q \right)$.
Finally, we have proved that 
$$
H_\eps (x,y,\eps t) =e^{\f(x)/2\eps} e^{-\f(y)/2\eps} K_\eps (x,y,t) ,
$$ 
and hence 
$$
\underset{\e\to 0}{\lim}\e \log K_{\e}(x,y, t)= \frac{\f(y)-\f(x)}{2} +\underset{\e\to 0}{\lim}\e \log H_{\e}(x,y,\eps t)
 = \frac{\f(y)-\f(x)}{2} -\rho(x,y,t) ,
 $$
after having applied~\eqref{e:Li-Yau}.
\enp

\subsection{$L^1$ observability estimates for positive solutions}
\label{s:L1obs}
We first prove intermediate observability statements in $L^1$. 
The following elementary abstract lemma shows that concerning positive solutions, observability in a (possibly weighted) $L^1$ norm is equivalent to the ``observability of the Kernel''. 
\begin{proposition}
\label{propequivobsKernel}
Let $\M$ be a compact Borel space (on which we denote by $dx$ a distinguished measure) and $T>0$. Take $K =K(x,y,t) \in C^0( \M \times \M\times(0,T] )$ be a nonnegative kernel. Assume further that $(t,y) \mapsto \nor{K(\cdot,y,t)}{L^1(\M)} = \int_\M K(x,y,t) dx$ is uniformly bounded for $(t,y)\in (0,T] \times \M$. 
Define for $t\in (0,T]$ the operator 
$$
S(t) : \Meas(\M) \to C^0(\M) , \quad 
\left[S(t)\mu \right](x)=\int_{\M} K(x,y,t)d\mu(y) = \langle \mu , K(x,\cdot,t) \rangle ,
$$ 
where $\Meas(\M)$ denotes the space of Radon measures on $\M$.
Let $w_1\in L^{\infty}([0,T]\times \omega)$ and $w_2\in L^{\infty}(\M)$ be two nonnegative weight functions. 

 Then, for all $T, s,C_0>0$, the following statements are equivalent
\begin{enumerate}
\item \label{enumObsmeasure}Observability of positive measures:
\baln
\label{e:meas-mu-obs}
\nor{w_2 \left[S(s)\mu\right]}{L^1(\M)}\leq C_0 \nor{w_1 [S(\cdot)\mu]}{L^1((0,T]\times\omega)}
\naln
for any $\mu \in  \Meas_+(\M)$ nonnegative Radon measure.
\item \label{enumObsL1} Observability of positive $L^1$ functions:
\bal
\nor{w_2 \left[S(s)u_0\right]}{L^1(\M)}\leq C_0 \nor{w_1\left[S(\cdot)u_0\right]}{L^1([0,T]\times \omega)}
\nal
for any $u_0\in L^1(\M, dx)$ with nonnegative value.
\item \label{enumObsDirac}Observability of Dirac distributions:
\bal
\nor{w_2 \left[S(s) \delta_y \right]}{L^1(\M)}\leq C_0 \nor{w_1 \left[S(\cdot) \delta_y \right]}{L^1([0,T]\times \omega)}
\nal
for any $y\in \M$. 
\item \label{enumObsKernel}Observability of the Kernel: 
\bal
I_s (y)\leq C_0 O_T(y) , \quad \text{ for all } y \in \M ,
\nal 
where
\bal
O_T(y)=\int_0^T \int_{\omega} w_1(t,y)K(x,y,t)dx dt , \quad \text{and} \quad 
I_s (y)=\int_{\M} w_2(x)K(x,y,s)dx . 
\nal
\end{enumerate}
\end{proposition}
Note first that under the assumption of the theorem, both terms in~\eqref{e:meas-mu-obs} are well-defined. Indeed, the Tonelli theorem (all functions/measures are nonnegative) implies 
\begin{align*}
 \nor{w_1 S(\cdot)\mu}{L^1([0,T]\times \omega)} & = \int_0^T \int_\omega w_1(t,x) [S(t)\mu](x) dx dt \\
 &  =  \int_0^T \int_\omega w_1(t,x) \int_{\M} K(x,y,t)d\mu(y) dx dt \\
 & =  \int_0^T   \int_{\M} \left(\int_\omega w_1(t, x) K(x,y,t)dx \right)d\mu(y)  dt \\
& \leq T \nor{w_1}{L^\infty([0,T]\times \M)}  \|\mu\|_{TV} \sup_{t\in [0,T]}  \sup_{y\in \M} \int_\omega K(x,y,t)dx 
< + \infty ,
\end{align*}
by assumption. Here $\|\mu\|_{TV}$ denotes the total variation of the measure $\mu$.
Note also that $u(t) = [S(t)(\mu)]$ is a continuous nonnegative function for any nonnegative measure $\mu$ and $t>0$.
Remark that $O_T(y)$ is essentially the observation of solutions starting from $\delta_y$ while $I_s(y)$ is the weighted norm of this solution at time $s$.

\bnp
Remark first that 
\begin{align*}
O_T(y)& =\int_0^T \int_{\omega} w_1(t,y)K(x,y,t)dx dt   = \int_0^T \int_{\omega} w_1(t,y)\langle \delta_y, K(x,\cdot,t ) \rangle dx dt \\
& = \nor{w_1 \left[S(\cdot) \delta_y \right]}{L^1([0,T]\times \omega)}\\
I_s (y)& =\int_{\M} w_2(x)K(x,y, s)dx = \int_{\M}w_2(x)\langle \delta_y ,  K(x,\cdot ,s)\rangle dx = \nor{w_2 \left[S(s) \delta_y \right]}{L^1(\M)} ,
\end{align*}
so that Item~\ref{enumObsDirac} $\Leftrightarrow$ Item~\ref{enumObsKernel}. Moreover, applying the Fubini Theorem, we get 
\begin{align*}
\nor{w_1 [S(\cdot)\mu]}{L^1([0,T]\times \omega)} &=\int_0^T \int_{x\in \omega}\int_{y\in \M} w_1(t,y)K(x,y,t)d\mu(y)dxdt=\langle \mu , O_T \rangle\\
\nor{w_2 \left[S(s)\mu\right]}{L^1(\M)}&=\int_{y\in \M}\int_{x\in \M} w_1(t,x)K(x,y,t)d\mu(y)dxdt=\langle \mu , I_s \rangle.
\end{align*}
Therefore, Item~\ref{enumObsmeasure} is equivalent to $\langle \mu , C_0O_T-I_s \rangle \geq 0$ for any $\mu \in \Meas_+(\M)$ and Item~\ref{enumObsL1} to $\int_\M u_0(y) (C_0O_T-I_s )(y) dy \geq 0$ for any $u_0 \in L^1_+(\M,dx)$. That Item~\ref{enumObsKernel} $\Leftrightarrow$ Item~\ref{enumObsmeasure} $\Leftrightarrow$ Item~\ref{enumObsL1} follows from the general fact that if $f \in C^0(\M)$, one has 
\begin{align*}
f \geq 0 \text{ on }\M & \Leftrightarrow \ \langle \mu , f \rangle \geq 0 \text{ for all }\mu \in \Meas_+(\M) \\
&  \Leftrightarrow \  \int_\M u_0(y) f(y) dy \geq 0\text{ for all }u_0 \in L^1_+(\M,dx).
\end{align*}
\enp
We now give the $L^1$ observability estimate for positive solutions to~\eqref{e:transport-viscous-gradient}.
\begin{proposition}
\label{propobsL1}
Assume that $T>T_{GCC}(\M, \nabla_g\f,\omega)$. Then, for any $\delta, s >0$, there exists $\e_0>0$ so that we have 
\baln
\label{estimL1+}
\nor{u(s)}{L^1(\M)}\leq e^{\frac{\delta}{\e}}\nor{u}{L^1([0,T]\times \omega)}  , \quad \text{ for all } \eps \in (0,\eps_0) ,
\naln
for any $u_0\in L^1(\M)$ with non-negative values and $u$ solution of \eqref{e:transport-viscous-gradient}.
\end{proposition}
Remark that we only use the case $s=T$ below. It is however remarkable that the stronger result for $s >0$ small holds as well. This is linked to the $L^1$ setting here.
Note also that in $L^1$, we have a ``converse inequality'', which we state for the sake of the comparison. Proposition~\ref{propobsL1} is proved afterwards.
\begin{lemma} 
\label{l:reverse-heat}
Assume $\d \M = \emptyset$. For all $T>0$, there is $C_T>0$ such that for all $u_0\in L^1(\M; \R^+)$ and $u$ the associated solution of~\eqref{e:transport-viscous-gradient}, we have 
\begin{equation}
\label{e:gronwall}
\nor{u(t)}{L^1(\M)}\leq C_T \nor{u(T)}{L^1(\M)}, \quad \text{ for all } t \in [0,T] \text{ and }\eps >0. 
\end{equation}
\end{lemma}
In particular, this implies that one cannot hope to replace the loss $e^{\frac{\delta}{\e}}$ by a gain $e^{-\frac{\delta}{\e}}$ in~\eqref{estimL1+}.

\bnp[Proof of Lemma~\ref{l:reverse-heat}]
Assume first that $u_0 \in W^{2,1}(\M)$ with $u_0\geq 0$ a.e. on $\M$. Then, notice that $u(t,x)\geq 0$ for a.e. $(t,x) \in (0,T) \times \M$.
Integrating~\eqref{e:transport-viscous-gradient} on $\M$, we obtain after an integration by parts (using that $\d \M = \emptyset$)
$$
\frac{d}{dt} \int_\M u  - \int_\M (\div_gX - q ) u = 0.
$$
Since $u\geq 0$, this implies $\frac{d}{dt} \nor{u(t,\cdot)}{L^1} \geq \gamma \nor{u(t,\cdot)}{L^1}$ for all $t>0$, with $\gamma = \inf_\M (\div_gX - q )$. The Gr\"onwall inequality yields~\eqref{e:gronwall}. The conclusion for a general $u_0 \in L^1$ follows from a density argument.  
\enp

We now turn to the proof of the $L^1$ observability estimate of Proposition~\ref{propobsL1}, which will use the following lemma.
\begin{lemma}
\label{l:pourri}
Assume $\big( \M, \nablag\f,\omega ,(\delta',T)\big)$ satisfies (GCC) (see Definition \ref{d:def-GCC}) for some $\delta'>0$.Then, for all $\delta>0$ there is $C_\delta>0$ such that for all $y \in \M$ there is an open set $U_y \subset \omega \times  (\delta',T)$ with $|U_y| \geq C_\delta$ and for all $(x,t)\in U_y$, we have $d_{\nablag \f}(x,y,t)\leq \delta$.
\end{lemma}
\bnp[Proof of Lemma~\ref{l:pourri}]
The assumption (GCC) implies that for any $y_0 \in \M$, there is $t_{y_0} \in (\delta',T)$ and $x_{y_0}=\phi_{-t_{y_0}}(y_0) \in \omega$ where $(\phi_t)_{t\in \R}$ is the flow of $\nablag\f$. The trajectory $\gamma(s)=\phi_s(x_{y_0})$ satisfies $\dot{\gamma}(s)=\nablag\f(s)$ with $\gamma(0)=x_{y_0}$ and $\gamma(t_{y_0})=y_0$ so that Proposition \ref{propdX} implies $d_{\nablag\f}(x_{y_0},y_0,t_{y_0})=0$. In particular, we obtain that for any $y_0 \in \M$, there is $t_{y_0} \in (\delta',T)$ and $x_{y_0} \in \omega$ such that $d_{\nablag\f}(x_{y_0},y_0,t_{y_0})=0$.

By uniform continuity of $d_{\nablag\f}$ (on the compact set $\M^2\times [\delta'/2,T+\delta']$) together with the fact that $\omega \times (\delta',T)$ is open, there exists $\nu_{y_0}>0$ so that $B_g(x_{y_0},\nu_{y_0})\times [t_{y_0}-\nu_{y_0}, t_{y_0} +\nu_{y_0}]\subset \omega\times (\delta',T)$
 and for any $y\in B_g(y_0,\nu_{y_0})$, $x\in B_g(x_{y_0},\nu_{y_0}))$ and $t\in [t_{y_0}-\nu_{y_0}, t_{y_0} +\nu_{y_0}]$, we have $d_{\nablag \f}(x,y,t)\leq \delta$. By compactness, we can cover $\M$ by $\M=\bigcup_{i \in I} B_g(y_i , \nu_{y_i})$ where $I$ is finite. Then, for any $y \in \M$, there is $i\in I$ such that $y \in B_g(y_i , \nu_{y_i})$, and the set $$U_y := B_g(x_{y_i},\nu_{y_i})\times [t_{y_i}-\nu_{y_i}, t_{y_i} +\nu_{y_i}]$$
satisfies the sought properties. 
\enp

As a consequence of this lemma together with Corollary~\ref{corLiYauconj} and Proposition \ref{propequivobsKernel}, we may now deduce a proof of the $L^1$ observability estimate of Proposition~\ref{propobsL1}.
\bnp[Proof of Proposition~\ref{propobsL1}]
Without any loss of generality, we can assume $0<\delta<s$.
According to Proposition~\ref{propequivobsKernel}, it is enough to study the ``observability of the Kernel''.
Using Corollary \ref{corLiYauconj}, for any $\delta>0$, there exists $\e_0$ so that~\eqref{inegkernel}
holds for all $(x,y)\in \M^2$, $t\in [\delta,\delta^{-1}]$, $0<\e\leq \e_0$.
Proposition \ref{propequivobsKernel} leads to compare
\bal
O_{T}(y): =\int_0^T \int_{x\in \omega} K_{\e}(x,y,t)dxdt
\quad \text{with} \quad 
I_{s}(y): =\int_{x\in \M} K_{\e}(x,y,s)dx
\nal
From~\eqref{inegkernel}, the fact that $K_{\e}(x,y,t)\geq0$ for $(x,y,t) \in \M^2 \times (0,\infty)$, and $s >\delta$, we deduce 
\bal
O_{T}(y) \geq \int_{\delta'}^T \int_{x\in \omega}e^{-\frac{d_{\nablag \f}(x,y,t)+\delta}{\e}} \quad \text{and} \quad 
I_{s}(y) \leq \int_{x\in \M} e^{-\frac{d_{\nablag \f}(x,y,s)-\delta}{\e}} ,
\nal
where $\delta'>0$ is chosen sufficiently small so that $\big(\M , \nablag\f,\omega,(\delta',T)\big)$ still satisfies (GCC), which is possible since $T>T_{GCC}(\M, \nabla_g\f,\omega)$.
Using Lemma~\ref{l:pourri} (where $U_y$ and $C_\delta$ are defined), we now have 
\baln
\label{uperO}
O_{T}(y)  \geq \int_{(x,t)\in U_y} e^{-\frac{d_{\nablag \f}(x,y,t)+\delta}{\eps}} dt dx 
\geq \int_{(x,t)\in U_y} e^{-\frac{2\delta}{\eps}} dt dx 
\geq C_\delta e^{-\frac{2\delta }{\e}}.
\naln
Also, for any $s >\delta$, using that $d_{\nablag \f}\geq 0$ (see Proposition \ref{propdX}), we have 
\bal
I_{s}(y)\leq \Vol_g (\M) e^{\frac{\delta}{\e}}. 
\nal
When combined with \eqref{uperO}, we obtain $O_{T}(y) \geq C e^{-3\frac{\delta}{\e}} I_{s}(y)$. By Proposition \ref{propequivobsKernel}, this gives \eqref{estimL1+} which concludes the proof of the proposition (up to changing $4\delta$ into $\delta$).
\enp

\subsection{From $L^1$ to $L^2$ observability estimates for positive solutions}
\label{s:L1toL2}

In this section, we conclude the proofs of Proposition~\ref{proplowerobse+} and Theorem~\ref{t:positive-intro}.
We first prove the negative result of Proposition~\ref{proplowerobse+} (uniform observability of positive solutions does not hold for $T<T_{GCC}(\M, \nabla_g\f,\omega)$, with an exponential lower bound of the cost).
\bnp[Proof of Proposition \ref{proplowerobse+}]Let us check the first part of the proposition (geometric statement). Since $(\M,\nablag \f,\overline{\omega},T)$ does not satisfy (GCC), there is $y_0\in \M$, so that for all $t \in [0,T]$, $\phi_{-t}(y_0)\notin \overline{\omega}$. In particular, for any $(t,x)\in [0,T]\times \overline{\omega}$, we have $\phi_t(x)\neq y_0$, which implies $d_{\nablag \f}(x,y_0,t)>0$ by Proposition~\ref{propdX}. By compactness of $[0,T]\times \overline{\omega}$, $\inf_{x\in\overline{\omega},t\in[0,T]}d_{\nablag \f}(x,y_0,t)>0$. Therefore, $d_{([0,T], \overline{\omega})}>0$ as expected.  

For the second part, for any $\delta>0$, select $y_0\in\M$ so that $$\inf_{(t,x)\in[0,T]\times \overline{\omega}}d_{\nablag \f}(x,y_0,t)\leq d_{([0,T], \overline{\omega})}+\delta/8.$$ By uniform continuity of $d_{\nablag \f}(x,y_0,t)$ defined on $[0,2T]\times \overline{\omega}$, we can also find $\eta>0$ so that $\inf_{(t,x)\in [0,T+\eta]\times \overline{\omega}}d_{\nablag \f}(x,y_0,t)\geq d_{([0,T], \overline{\omega})}-\delta/4$ and $d_{\nablag \f}(x,y_0,\eta )\leq \delta/4$ for $x\in B(y_0,\eta)$. We take as initial datum $u_{0}(x)=K_{\e}(x,y_0,\eta)$, yielding $u(t,x) = K_{\e}(x,y_0,t+ \eta)$ (see the definition of $K_\eps$ in~\eqref{e:transport-viscous-kernel}) as the associated solution of~\eqref{e:transport-viscous-gradient}. We have $u_{0}\in L^2(\M)$ together with a lower bound coming from \eqref{inegkernel} with a sufficiently small $\delta$ (replaced by $\delta/4$)
\bal
\nor{u_0}{L^2(\M)}\geq \nor{u_0}{L^2(B(y_0,\eta))}\geq C(\eta) e^{-\frac{\delta}{2\e}}.
\nal
Concerning the observation term, we deduce from the upper bound in \eqref{inegkernel} that for $\e$ small enough,
\bal
 \int_0^T\int_\omega |u(t,x)|^2 dx dt & =  \int_0^T\int_{\overline{\omega}} K_{\e}(x,y_0,t+ \eta)^2  dx dt \\
 & \leq \Vol_g(\omega)T\sup_{x\in\overline{\omega},t\in[\eta,T+\eta]}e^{-2\frac{d_{\nablag \f}(x,y_0,t)-\delta/4}{\e}}
  \leq C e^{-2\frac{d_{([0,T], \overline{\omega})}-\delta/2}{\e}}.
\nal
Applying the observability inequality~\eqref{e:transport-viscous-obsplus} to $u$ thus implies that $$C_0^+(T,\eps)^2 e^{-2\frac{d_{([0,T], \overline{\omega})}-\delta/4}{\e}}\geq C, $$ uniformly in $\eps \in (0,\eps_0]$, which concludes the proof of the proposition.
\enp

To conclude the proof of Theorem~\ref{t:positive-intro}, we need the following dissipation result taken from Guerrero-Lebeau \cite{GL:07}. In that reference, it is written on an open subset $\Omega \subset \R^n$ with the flat metric; however, it can be checked that the result also applies to the case of a Riemannian manifold $(\M,g)$ without boundary and with an additional potential $q$.
\begin{proposition}[Proposition 3 of \cite{GL:07}, $m=1$]
\label{p:guerrero-lebeau}
Assume that $\d\M=\emptyset$ and that $(\M,X,\omega,T)$ satisfies (GCC). Then, there exist $C,C_0>0$ such that
\baln
\label{dissip}
\nor{u(T)}{L^{2}(\M)}^2 \leq C \left(\int_{0}^{T}\int_{\omega}|u|^{2}dtdx+e^{-\frac{C_{0}}{\e}}  \nor{u(0)}{L^{2}(\M)}^2 \right)
\naln
for all $\eps \in (0,1]$ and any (not necessarily positive) solution $u$ to \eqref{e:transport-viscous} (and a fortiori for all solutions $u$ to \eqref{e:heat-transp-eps}).
The same statement holds true if $\d\M\neq\emptyset$ and $(\M,X,\omega,T)$ satisfies (FC).
\end{proposition}

We shall also need the following lemma in the proof of Theorem \ref{t:positive-intro}.
\begin{lemma}
\label{lmregular}
For any $\delta'>\delta>0$, there exists $\e_{0}>0$ so that
\baln
\label{regulL1L2}
\nor{u(\delta')}{L^{2}(\M)}\leq C \nor{u(\delta')}{L^{\infty}(\M)}\leq Ce^{\frac{\delta}{\e}}  \nor{u(\delta)}{L^{1}(\M)} 
\naln
for any solution $u$ to \eqref{e:heat-transp-eps} and $0<\e\leq \e_{0}$. 
\end{lemma}
\bnp
Since the manifold is compact, we only need to prove the $L^{\infty}$ bound which follows from the bound $\nor{K_{\e}(\cdot,\cdot,\delta' -\delta)}{L^{\infty}(\M \times \M)}\leq Ce^{\frac{\delta}{\e}} $ on the Kernel. This estimate follows from Corollary \ref{corLiYauconj} (e.g.~\eqref{inegkernel} together with the fact that $d_{\nablag \f} \geq0$, see Proposition~\ref{propdX}).
\enp
\bnp[Proof of Theorem \ref{t:positive-intro}]
Inequality \eqref{lowercost+} of Proposition \ref{proplowerobse+} directly yields $T_{unif}^{+}(\omega)\geq T_{GCC}(\M, \nabla_g\f, \omega)$. Note that it was mostly proved in Guerrero-Lebeau \cite[Theorem 1]{GL:07} since one can check that the counterexample they build is a nonnegative solution. 

Now, we prove $T_{unif}^{+}(\omega)\leq T_{GCC}(\M, \nabla_g\f,\omega)$. 
For any $\delta>0$ (we will later need $2\delta<C_{0}$ where $C_{0}$ is the constant in \eqref{dissip}), and for $T\geq T_{GCC}(\M, \nabla_g\f,\omega)+2\delta$, we prove the observability inequality for positive solution
\baln
\label{observproof}
\nor{u(T)}{L^{2}(\M)}^2\leq C \int_{0}^{T}\int_{\omega}|u|^{2}dtdx.
\naln
The combination of \eqref{dissip} on the time interval $(2\delta,T)$ together with~\eqref{regulL1L2} on the time interval $(\delta, 2\delta)$ implies
\baln
\label{estimintermL2obs}
\nor{u(T)}{L^{2}(\M)}^2 \leq C \left(\int_{2\delta}^{T}\int_{\omega}u^{2}dtdx+e^{-\frac{C_{0}-\delta}{\e}}  \nor{u(\delta)}{L^{1}(\M)}^2  \right) .
\naln
Now, applying \eqref{estimL1+} with $\lambda>0$ such that $\lambda T=\delta$, we obtain 
\baln
\label{estimdisspidelta}
\nor{u(\delta)}{L^{1}(\M)} \leq e^{\frac{\delta}{\e}}\int_{0}^{T}\int_{\omega}u dtdx .
\naln
Combining \eqref{estimintermL2obs}, \eqref{estimdisspidelta} together with the H\"older inequality, we deduce
\bal
\nor{u(T)}{L^{2}(\M)}\leq \left(C+ Ce^{-\frac{C_{0}-2\delta}{\e}}\right)\int_{0}^{T}\int_{\omega}u^{2}dtdx .
\nal
Choosing $\delta\in (0,C_0/2)$ implies~\eqref{observproof} uniformly for $\eps \in (0,\e_{0}(\delta)]$, and hence concludes the proof of the theorem.
\enp

\subsection{From observability of positive solutions to a controllability statement}
\label{s:cont-pos}
This section is devoted to the proof of the controllability result of Corollary~\ref{corcontrol+} from the observability of positive solutions.
It relies on the following lemma. The result and its proof follow~\cite[Theorem 4.1]{L:18}. 

\begin{lemma}
\label{lm:dualiteconvex}
Let $\mathcal{V}$ be a closed convex set of $L^2(\M)$ with $0\in \mathcal{V}$ and $\widetilde{\mathcal{V}}\subset \mathcal{V}$ so that $\widetilde{\mathcal{V}}\subset\mathcal{V}-v $ for any $v\in \mathcal{V}$. Assume moreover that there exists $C_{\mathcal{V}}>0$ so that
 \begin{multline}
 \label{e:observV}
C_{\mathcal{V}}(T,\eps)^2 \int_0^T\int_\omega |u(t)|^2 ds(x)dt \geq \|u (T)\|_{L^2(\M)}^2 , \\  \text{ for all } u_0 \in \mathcal{V}\text{ and } u \text{ solution of~\eqref{e:transport-viscous}}.
 \end{multline}
Then, for any $y_{0}\in L^{2}(\M)$ and $0<\e\leq \e_{0}$, there exists a control $h\in L^2([0,T],L^2(\M))$ with $$\nor{h}{L^2([0,T],L^2(\M))}\leq C_{\mathcal{V}}(T,\eps)\nor{y_{0}}{L^{2}(\M)}$$ so that the solution of \eqref{e:transport-viscous-control-internal} satisfies $\left( y(T),u_0\right)_{L^2(\M)}\geq 0$ for any $u_0\in \widetilde{\mathcal{V}}$.
\end{lemma}
For the proof of Corollary~\ref{corcontrol+}, we apply this lemma to the sets $\widetilde{\mathcal{V}}= \mathcal{V}=L^2(\M; \R^+)$. Notice that  Lemma~\ref{lm:dualiteconvex} also contains one implication (namely Observability $\implies$ Controllability) in Corollary~\ref{c:obscontrol} when applied to $\widetilde{\mathcal{V}}= \mathcal{V}=L^2(\M; \R)$. 

\bnp
For any $\alpha>0$, we consider the functional $J_{\alpha}$ defined for any $u_0\in \mathcal{V}$ by
\bal
J_{\alpha}(u_0)=\frac12 \int_0^T\int_\omega |u(t,x)|^2 dtdx + \frac{\alpha}{2} \nor{u_0}{L^2(\M)}^2+\left( u (T),y_0\right)_{L^2(\M)} ,
\nal
where $u$ is the solution of~\eqref{e:transport-viscous}. The functional $J_{\alpha}$ is continuous, convex and coercive. Therefore, $J_{\alpha}$ admits a minimum $u_{0,\alpha}\in \mathcal{V}$ (se e.g.~\cite[Chapter II, Proposition 1.2]{ET:74}). The minimality condition gives (see e.g.~\cite[Chapter II, Proposition 2.1]{ET:74}) for any $p_0\in L^{2}(\M)$ that can be written $p_0=v_0-u_{0,\alpha}$, with $v_0\in \mathcal{V}$,
\baln
\label{e:Lagrange}
\int_0^T\int_\omega p(t,x) u_{\alpha}(t,x)dtdx+ \alpha\left( p_0,u_{0,\alpha}\right)_{L^2(\M)}+ \left( p (T),y_0\right)_{L^2(\M)} \geq 0.
\naln
where we have denoted $p$ (resp. $u_\alpha$) the solution of~\eqref{e:transport-viscous} with $p(0)=p_0$ (resp. $u(0)=u_{0,\alpha}$). 

Now, let $y_{\alpha}$ be the solution of \eqref{e:transport-viscous-control-internal} with control function $h_{\alpha}(t,x)=u_{\alpha}(T-t,x)$ and initial datum $y_{\alpha}(0) =y_0$. The duality equation \eqref{e:duality-equation-1} gives for any $p_0\in L^2(\M)$ initial datum for $p$ solution of~\eqref{e:transport-viscous}
\bal
\int_0^T\int_\omega p(t,x) u_{\alpha}(t,x)dtdx= \left( p_0,y_{\alpha}(T)\right)_{L^2(\M)}-\left( p(T),y_0\right)_{L^2(\M)}.
\nal
Combined with~\eqref{e:Lagrange}, this implies 
\baln
\label{e:dualitealpha}
\left( p_0,y_{\alpha}(T)\right)_{L^2(\M)}+ \alpha\left( p_0,u_{0,\alpha}\right)_{L^2(\M)}\geq 0 ,
\naln
 for every $p_0\in  \mathcal{V}-u_{0,\alpha}$. This also holds for any $p_0\in \widetilde{\mathcal{V}}$ since $\widetilde{\mathcal{V}}\subset\mathcal{V}-u_{0,\alpha}$ by assumption. 
 
 To obtain an estimate of the control, we apply \eqref{e:Lagrange} to $p_0=0-u_{0,\alpha}$. After an application of the Cauchy-Schwarz inequality, we have
\bal
\int_0^T\int_\omega |u_{\alpha}(t,x)|^2dtdx+\alpha\nor{u_{0,\alpha}}{L^2(\M)}^2\leq \nor{u_{\alpha}(T)}{L^2(\M)}\nor{y_0}{L^2(\M)}.
\nal
The observability inequality~\eqref{e:observV} applies to $u_{0,\alpha}\in \mathcal{V}$, so that
\begin{multline*}
\int_0^T\int_\omega |u_{\alpha}(t,x)|^2dtdx+\alpha\nor{u_{0,\alpha}}{L^2(\M)}^2 \\
\leq C_{\mathcal{V}}(T,\eps)\left(\int_0^T\int_\omega |u_{\alpha}(t,x)|^2dtdx\right)^{1/2}\nor{y_0}{L^2(\M)}.
\end{multline*}
We obtain successively
\baln
\label{estimuobservalpha}\int_0^T\int_\omega |u_{\alpha}(t,x)|^2dtdx&\leq C_{\mathcal{V}}(T,\eps)^2\nor{y_0}{L^2(\M)}^2,\\
\label{estimu0alpha}\alpha\nor{u_{0,\alpha}}{L^2(\M)}^2&\leq C_{\mathcal{V}}(T,\eps)^2\nor{y_0}{L^2(\M)}^2 .
\naln
We obtain that $\mathds{1}_{\omega}u_{\alpha}$ is bounded in $L^2([0,T]\times \M)$ uniformly in $\alpha>0$. Take a sequence $\alpha_n\to 0$ so that $\mathds{1}_{\omega}u_{\alpha_n}\rightharpoonup \mathds{1}_{\omega}u$ in $L^2([0,T]\times \M)$. The associated solutions $y_{\alpha_n}$ with control $\mathds{1}_{\omega}u_{\alpha_n}(T-t,x)$ is therefore bounded in $L^{\infty}([0,T],L^2(\M))$ and, again up to a subsequence, converges weakly-$*$ to a solution $y$ of~\eqref{e:transport-viscous-control-internal} with control $h(t,x)=\mathds{1}_{\omega}u(T-t,x)$ and initial datum $y_0$. Moreover, up to a subsequence, we can impose $y_{\alpha_n}(T)\rightharpoonup y(T)$ in $L^2(\M)$. Passing to the limit in \eqref{e:dualitealpha} using \eqref{estimu0alpha}, we finally obtain
\bal
\left( p_0,y(T)\right)_{L^2(\M)}\geq 0
\nal
for any $p_0\in \widetilde{\mathcal{V}}$. We finally get the expected estimate on $h(t,x)=u(T-t,x)$ passing to the limit in \eqref{estimuobservalpha}.
\enp

We may now conclude the proof of Corollary \ref{corcontrol+} from Lemma \ref{lm:dualiteconvex}.
\bnp[Proof of Corollary \ref{corcontrol+}]
We apply Lemma \ref{lm:dualiteconvex}, with $\widetilde{\mathcal{V}}= \mathcal{V}=L^2(\M; \R^+)$. Note that the Lemma applies because for any $v\in L^2(\M; \R^+)$, $L^2(\M; \R^+)\subset L^2(\M; \R^+)-v$. Indeed, any $u \in L^2(\M; \R^+)$ can be written $u=(u+v)-v\in L^2(\M; \R^+)-v$ since $u+v\geq 0$. This gives a control $h$ with the expected uniform bound and so that $\left( y(T),u_0\right)_{L^2(\M)}\geq 0$ for any $u_0\in L^2(\M; \R^+)$. This implies $y(T)\geq 0$.
\enp

\appendix

\section{About the distances}
\label{sectdist}
In this section, $(\M,g)$ is a compact Riemannian manifold without boundary.
\subsection{A general lemma}
We start with a general lemma.
\begin{lemma}
\label{lmparamgeod}
Let $\mathsf{V} \in W^{1,\infty}(\M)$ with nonnegative value. Then, for all $x,y\in \M^2$, we have 
\bal
&\frac{1}{2} \inf_{\gamma, t} \left\{ \int_0^t |\dot{\gamma}(s)|_g^2 +\mathsf{V}(\gamma(s)) ds , \ t>0,  \gamma \in \mathbf{U}_t(x,y)\right\}\\
&=\inf_{\gamma, t} \left\{ \int_0^t |\dot{\gamma}(s)|_g\sqrt{\mathsf{V}(\gamma(s))} ds , \ t>0,   \gamma \in \mathbf{U}_t(x,y) \right\}\\
&=\inf_{\gamma} \left\{ \int_0^1 |\dot{\gamma}(s)|_g\sqrt{\mathsf{V}(\gamma(s))} ds , \ \gamma \in \mathbf{U}_1 (x,y) \right\}.
\nal
where $\mathbf{U}_t(x,y)= \left\{ \gamma \in W^{1,\infty} ([0,t] ; \M ), \gamma(0)=x , \gamma(t)= y \right\}$ for $t>0$.
\end{lemma}
This lemma is particularly useful for $\mathsf{V} = (V-E)_+$ in which case the (pseudo-) distance defined is the Agmon distance at energy level $E$. 

\bnp
We denote by $d_1$, $d_2$, $d_3$ respectively the three (pseudo-)distances defined in the statement of the lemma.
Then, we notice that the last two quantities are invariant by reparametrization, so that $d_2=d_3$ after a change of variable in the integral. Then, the inequality $ab\leq \frac{1}{2}(a^2+b^2)$ directly yields $d_2\leq d_1$. Let us now prove the converse inequality, namely $d_3 \geq d_1$. For $\e>0$ there exist $\delta>0$ and a path $\gamma:[0,1]\to \M$ such that $ \int_0^1 |\dot{\gamma}(s)|_g\sqrt{\mathsf{V}(\gamma(s))+\delta} \, ds\leq d_{3}+\e$. We can further assume $|\dot{\gamma}(s)|_g>0$ with the same estimate (indeed, defining a new parametrization $\zeta$ by $\gamma(t) = \zeta(\phi(t))$ with $\phi(t) = \int_0^t |\dot{\gamma}(s)|_g ds$ even yields a Lipschitz reparametrization with constant positive speed, see e.g.~\cite[Proof of Lemma~3.16]{BarilariAgrachevBoscainBook:20}). Using an approximation argument, we can further assume that $\gamma$ is smooth up to replacing $\e$ by $2\eps$.
We now define the following reparametrization $\widetilde{\gamma}(s)=\gamma(\varphi(s))$ where $\varphi$ solves $\dot{\varphi}(s)=\frac{\sqrt{\mathsf{V}(\gamma(\varphi(s)))+\delta}}{|\dot{\gamma}(\varphi(s))|_g}\geq \frac{\sqrt{\delta}}{\max_{[0,1]} |\dot{\gamma}|_g}>0$, $\varphi(0)=0$ so that $|\dot{\widetilde{\gamma}}(s)|_g=\sqrt{\mathsf{V}(\widetilde{\gamma}(s))+\delta} $ for any $s\in [0,\varphi^{-1}(1)]$. In particular, $|\dot{\widetilde{\gamma}}(s)|_g\sqrt{\mathsf{V}(\widetilde{\gamma}(s))+\delta}=\frac12\left(|\dot{\widetilde{\gamma}}(s)|_g^2+V(\widetilde{\gamma}(s))+\delta\right) $ and
\begin{align*}
 \frac12 \int_0^{\varphi^{-1}(1)}\left(|\dot{\widetilde{\gamma}}(s)|_g^2+\mathsf{V}(\widetilde{\gamma}(s))+\delta\right) ds
 &  = \int_0^{\varphi^{-1}(1)} |\dot{\widetilde{\gamma}}(s)|_g\sqrt{\mathsf{V}(\widetilde{\gamma}(s))+\delta}  ds \\
 &=\int_0^1 |\dot{\gamma}(s)|_g\sqrt{\mathsf{V}(\gamma(s))+\delta} ds\leq d_{3}+2\e ,
\end{align*}
which gives $d_1\leq d_3+2\e$, and concludes the proof of the lemma.
\enp

\subsection{Equivalence between the two definitions of $d_{\nablag \f}$}

In this section, we prove equivalence between the two definitions of $d_{\nablag \f}$, respectively given in~\eqref{defdXintro} in the introduction and in~\eqref{defdX}. 
We further give an equivalent quantity in terms of the Riemannian distance and the flow $\phi_t$ of $\nablag\f$. The function $\f$ is assumed to be $C^{3}$ throughout the section.

We recall that $\rho$ is defined in~\eqref{e:def-rho-V} with $V= \frac{|\nabla_g \f|_g^2}{4}$, that is to say
\baln
\label{e:def-rho}
\rho(x,y,t)&:=\frac{1}{4} \inf \left\{ \int_0^t |\dot{\gamma}(s)|_g^2 +|\nabla_g\f(\gamma(s))|_g^2 ds , \ \gamma \in \mathbf{U}_t(x,y)  \right\} ,
\naln
with $\mathbf{U}_t(x,y)= \left\{ \gamma \in W^{1,\infty} ([0,t] ; \M ), \gamma(0)=x , \gamma(t)= y \right\}$.
Note that it is proved in~\cite[Appendix]{LY:86} that $\rho$ is continuous on $\M^2 \times (0,+\infty)$ and, for all $t>0$ fixed, Lipschitz continuous as a function of $(x,y)\in \M^2$.
These quantities are related to the Agmon distance but in finite time, see Section~\ref{s:app-further-links} below. Note that the quantity $\rho$ is symmetric, $\rho(x,y,t) = \rho(y,x,t)$, and remains unchanged under the change of $\f$ by $-\f$. 
This is not the case for $d_{\nablag \f}$. 

\begin{lemma}
\label{lmdist}
The function $d_{\nablag \f}$ defined as
\baln
\label{e:def-dX}
d_{\nablag \f}(x,y,t) &:= \rho(x,y,t)+\frac{\f(x)-\f(y)}{2} 
\naln
is continuous on $\M^2 \times (0,+\infty)$ and, for all $t>0$ fixed, Lipschitz continuous as a function of $(x,y)\in \M^2$. Moreover, we have
\begin{align}
\label{e:equiv-def-dx}
d_{\nablag \f}(x,y,t)&=\frac{1}{4} \inf \left\{ \int_0^t \left|\dot{\gamma}(s)+\nablag\f(\gamma(s))\right|_g^{2} ds ,  \ \gamma \in \mathbf{U}_t(x,y) \right\} , \\
\label{defdXappend}&=\frac{1}{4} \inf \left\{ \int_0^t \left|\dot{\gamma}(s)-\nablag\f(\gamma(s))\right|_g^{2} ds , \ \gamma \in \mathbf{U}_t(x,y) \right\} ,\\
\label{defdXbisappend}&=\frac{1}{4} \inf \left\{ \int_0^t \left|\dot{\gamma}(s)\right|_{g_s}^{2} ds , \ \gamma \in W^{1,\infty} ([0,t] ; \M ), \gamma(0)=x , \gamma(t)= \phi_{-t}(y) \right\} ,
\end{align}
where $\mathbf{U}_t(x,y)= \left\{ \gamma \in W^{1,\infty} ([0,t] ; \M ), \gamma(0)=x , \gamma(t)= y \right\}$, and
 $g_s$ is the time varying metric defined by $ \left|Y\right|_{g_s}= \left|D\phi_s(Y)\right|_{g}$.

In particular, for any $T\geq 0$, there exists some constant $C_T>0$ so that
 $$
 C_T^{-1}d(x,\phi_{-t}(y))\leq d_{\nablag \f}(x,y,t)\leq C_T d(x,\phi_{-t}(y)), \quad \text{ for all } t\in [0,T], 
 $$
 where $d$ denotes the Riemannian distance (associated to $g$).
\end{lemma}

\bnp
The continuity property directly follows from that of $\rho$ proved in~\cite[Appendix]{LY:86}.
To prove~\eqref{e:equiv-def-dx}, we remark that for any path $\gamma$ so that $\gamma(0)=y , \gamma(t)= x$, we have $\f(x)-\f(y)=\int_0^t \frac{d}{ds}(f\circ \gamma)(t) dt =\int_0^t \nablag\f (\gamma(t))\cdot \dot{\gamma}(t)dt$. In particular, from the definition of $d_{\nablag \f}$ in~\eqref{e:def-dX} and $\rho$ in~\eqref{e:def-rho}, we can rewrite $d_{\nablag \f}(x,y,t)$ as
\bal
d_{\nablag \f}(x,y,t)&=\frac{1}{4} \inf \bigg\{ \int_0^t |\dot{\gamma}(s)|^2 +|\nablag\f(\gamma(s))|^2 +2\nablag\f(\gamma(s))\cdot \dot{\gamma}(s)ds  ,\  \gamma \in \mathbf{U}_t(x,y) \bigg\}\\
&=\frac{1}{4} \inf \left\{ \int_0^t \left|\dot{\gamma}(s)+\nablag\f(\gamma(s))\right|^{2} ds ,\  \gamma \in  \mathbf{U}_t(x,y) \right\}.
\nal
The statement~\eqref{defdXappend} is obtained thanks to the change of path $\widetilde{\gamma}(s)=\gamma(t-s)$. 

Now, we compute $d_{\nablag \f}(x,\phi_t (y),t)$ according to formula \eqref{defdXappend}. To this aim, let $\gamma \in W^{1,\infty} ([0,t] ; \M )$ so that $\gamma(0)=x$, $\gamma(t)=\phi_t(y)$. Let $\widetilde{\gamma}(s)=\phi_{-s}(\gamma(s))$ so that $\dot{\gamma}(s)=\nablag\f(\gamma(s))+D\phi_s(\dot{\widetilde{\gamma}}(s))$. In particular, $\left|\dot{\gamma}(s)-\nablag\f(\gamma(s))\right|_g^{2}=\left|D\phi_s(\dot{\widetilde{\gamma}}(s))\right|_g^{2}=\left|\dot{\widetilde{\gamma}}(s)\right|_{g_s}^{2}$. It gives \eqref{defdXbisappend} since any path $\gamma \in W^{1,\infty} ([0,t] ; \M )$ so that $\gamma(0)=x$, $\gamma(t)=y$ can be written $\gamma(s)=\phi_{s}(\widetilde{\gamma}(s))$ with $\widetilde{\gamma}(0)=x$, $\widetilde{\gamma}(t)=\phi_{-t}(y)$, and conversely.
\enp

\subsection{Further links between the different distances}
\label{s:app-further-links}
In this section, we relate the above quantities $\rho(x,y,t), d_{\nablag \f}(x,y,t)$ with the Agmon distance to the bottom energy (see~\ref{defAgmon-0} for $V =  \frac{| \nabla_g \f(x)|_g^2}{4}$ and $E_0 = \min_\M V = 0$), that is to say
\baln
\label{e:def-dA-Appendix}
d_{A}(x,y) &= \frac12 \inf \left\{ \int_0^1|\nabla_g \f(\gamma(s))|_g  |\dot{\gamma}(s)|_g  ds, \  \gamma \in  \mathbf{U}_1(x,y) \right\} ,
\naln
with $\mathbf{U}_1(x,y)= \left\{ \gamma \in W^{1,\infty} ([0,1] ; \M ), \gamma(0)=x , \gamma(1)= y \right\}$ and the associated quantity (compare with the definition of $d_{\nablag \f}$ in terms of $\rho$ in~\eqref{e:def-dX})
\baln
\label{e:def-W-Appendix}
W(x,y) &:= d_{A}(x,y)+\frac{\f(x)-\f(y)}{2}  .
\naln
The results of this section are not explicitly used in the proofs of the main part of the paper; however we believe these links are interesting and enlightening. Indeed, they relate the quantity $d_{A}(x,y)$ appearing in all general bounds of Section~\ref{s:general-bounds} together with the quantities $\rho(x,y,t), d_{\nablag \f}(x,y,t)$ appearing in results of Section~\ref{s:positive} concerning positive solutions.
\begin{lemma}
For all $(x,y)\in \M^2$, we have
\baln
\label{dAinf}
d_{A}(x,y)&=\inf_{t>0}\rho(x,y,t) ,\\ 
\label{Winf}W(x,y)&=\inf_{t>0}d_{\nablag \f}(x,y,t) .  
\naln
Moreover, if $\nabla_g\f(y)=0$, then we have $d_{A}(x,y) =\underset{t\to+\infty}{\lim}\rho(x,y,t)$.

\end{lemma}
\bnp
Equality in~\eqref{dAinf} is a consequence of Lemma \ref{lmparamgeod} applied to $\mathsf{V}=|\nabla_g\f|_g^2$. 
Then, \eqref{Winf} is a direct consequence of the expression of $W$ and $d_{\nablag \f}$ in terms of $d_A$ and $\rho$ in \eqref{e:def-W-Appendix}, \eqref{e:def-dX}, together with~\eqref{dAinf}.

Finally, if $\nabla_g\f(y)=0$, then the function $t\mapsto \rho(x,y,t)$ is non-increasing. Indeed, taking $t_1\leq t_2$, from a path $\gamma_1: [0,t_1]\to \M$ such that $\gamma_1(0)=x$ and $\gamma_1(t_1)=y$, we can construct the path $\gamma_2:[0,t_2]\to \M$ by $\gamma_2(s)=\gamma_1(s)$ for $s\in [0,t_1]$ and $\gamma_2(s) = y$ for $s\in[t_1,t_2]$. This yields a path in $W^{1,\infty}([0,t_2];\M)$ if $\gamma_1 \in W^{1,\infty}([0,t_1];\M)$, and thus the set of admissible paths on $[0,t_2]$ is larger than the set of admissible paths on $[0,t_1]$. Since the contribution $\int_{t_1}^{t_2} |\dot{\gamma}_2(s)|^2 + |\nabla_g\gamma_2(s)|_g^2 ds= 0$ we deduce that $\rho(x,y,t_2)\leq \rho(x,y,t_1)$. This proves that the $\inf$ is actually a $\underset{t\to+\infty}{\lim}$ in this case.
\enp

Note that related properties are proved in the Appendix of \cite{HS:85}. For instance, \cite[Lemma A2.2]{HS:85} with our notations can be loosely stated as follows:
If $W(x,y)=0$, then every minimizing geodesic of $d_{A}$ is a generalized integral curve of $\nabla_g \f$.

\bigskip 
Finally, we state a last result that explains that $d_{\nablag \f}(x,y,t)$ measures how far $x$ is the final state of a path of the vector field at time $t$ and starting at $y$.
Part of this result is contained in the last statement of Lemma~\ref{lmdist}; we here give a different proof, which, we believe, is interesting in itself.
\begin{proposition}
\label{propdX}
With $d_{\nablag \f}(x,y,t)$ defined in \eqref{defdX}, we have for all $(x,y,t) \in \M^2\times (0,+\infty)$,
\begin{enumerate}
\item \label{enum1dX}$d_{\nablag \f}(x,y,t)\geq 0$;
\item \label{enum2dX}$d_{\nablag \f}(x,y,t)=0$ if and only if there exists a trajectory of $\dot{\gamma}(s)=\nabla_g\f(\gamma(s))$ with $\gamma(0)=x , \gamma(t)= y$, that is if and only if $y=\phi_t(x)$.
\end{enumerate}
In particular, $(\M,\nablag\f,\omega,T)$ satisfies (GCC) if and only if for any $x\in \M$, there exist $y\in \omega$ and $t\in (0,T)$ so that $d_{\nablag \f}(x,y,t)=0$.
\end{proposition}
Recall that the flow $(\phi_t)_{t\in \R}$ is defined in~\eqref{e:def-flot} and the Geometric Control Condition (GCC) is defined in Definition~\ref{d:def-GCC}.
\bnp
Statement \ref{enum1dX} follows from the definition of $d_{\nablag \f}$ in~\eqref{defdXappend}. 
Let us now consider Statement~\ref{enum2dX}. 
Assume first that there exists a trajectory of $\dot{\gamma}(s)=\nablag\f(\gamma(s))$ with $\gamma(0)=x$ and $\gamma(t)=y$. Then, by definition of the infimum in~\eqref{defdXappend}, this yields $d_{\nablag \f}(x,y,t)\leq 0$ and hence $d_{\nablag \f}(x,y,t)=0$. Conversely, assume $d_{\nablag \f}(x,y,t)=0$. Take a minimizing sequence in~\eqref{defdXappend}, that is to say that $\gamma_n \in W^{1,\infty}([0,t]; \M)$ such that $\gamma_n(0)=x$, $\gamma_n(t)=y$ and $\dot{\gamma}_n-\nablag\f (\gamma_n)=R_n$ (bounded continuous with values in the tangent bundle of $\M$) with
\baln
\label{e:Rn-to-zer}
\int_0^t |R_n(s)|_g^2 ds \to 0.
\naln
Since $\nablag\f$ is bounded on $\M$, the sequence $\int_0^t |\dot{\gamma}_n(s)|_g^2ds$ is then uniformly bounded in $\R$. As a consequence, the sequence of paths $(\gamma_n)_{n\in \N}$ is equicontinuous. From Ascoli's theorem, we may extract a subsequence (which we do not relabel) $(\gamma_n)_{n\in \N}$ which converges strongly for the topology $C^0([0,t]; \M)$ to a limit $\gamma \in C^0([0,t]; \M)$. The latter thus has the same and endpoints $\gamma(0)=x$ and $\gamma(t)=y$. It is solution of $\dot{\gamma}=\nabla_g \f (\gamma)$ in the distributional sense according to~\eqref{e:Rn-to-zer}. Bootstrapping in the differential equation implies $\gamma \in W^{1,\infty}([0,1];\M)$ and $\gamma$ is a strong solution to $\dot{\gamma}=\nabla_g \f (\gamma)$. This concludes the proof of Statement~\ref{enum2dX}. 

 According to Lemma~\ref{l:equiv-GCC} Item~\ref{i:-GCC-X},  that $(\M,X,\omega,T)$ satisfies (GCC) is equivalent to the fact that for any $x\in \M$, there exist $t\in (0,T)$ such that $y:= \phi_t(x)\in \omega$. As a consequence of Item~\ref{enum2dX}, this is equivalent to having $d_{\nablag \f}(x,y,t)=0$.
\enp

Note that an analogue statement for the function $W$ in~\eqref{Winf} is proved in~\cite[Lemma A2.2]{HS:85}, and could also be deduced from Proposition~\ref{propdX}.

\small
\bibliographystyle{alpha}
\bibliography{bibli}
\end{document}